\documentclass[12pt]{amsart}
\newtheorem{theorem}{Theorem}
\newtheorem{lemma}{Lemma}

\newtheorem{corollary}{Corollary}
\newtheorem{proposition}{Proposition}
\newtheorem{definition}{Definition}
\newtheorem{remark}{Remark}

\usepackage{appendix}
\numberwithin{equation}{section}
\numberwithin{theorem}{section}
\numberwithin{corollary}{section}
\numberwithin{proposition}{section}
\numberwithin{definition}{section}
\numberwithin{lemma}{section}
\numberwithin{sublemma}{section}
\numberwithin{remark}{section}
\numberwithin{convention}{section}
\usepackage{mathrsfs}
\usepackage{diagbox}

\begin{document}
\title{The Isoparametric Story, a Heritage of \'{E}lie Cartan}
\author{Quo-Shin Chi}
%\address{Department of Mathematics, Washington University, St. Louis, MO 63130}
%\email{chi@math.wustl.edu}
%%\thanks{}

%%\date{}

\begin{abstract} In this article, we survey along the historical route the classification of isoparametric hypersurfaces in the sphere, paying attention to the employed techniques in the case of four principal curvatures.
\end{abstract}
\maketitle
%\pagestyle{myheadings}
%\markboth{QUO-SHIN CHI}{ISOPARAMETRIC HYPERSURFACES}
%*************
%%\footnote{Keywords: Isoparametric hypersurfaces; Subjclass: Primary, 53C40}

\section{Prologue}  When writing this survey article, I kept in mind that there had been an extensive body of research papers on the subject of isoparametric submanifolds and beyond. With the addition of the comprehensive books~\cite{BCO,CR} that wrapped up much of what had been known up to the time of their publications, and the survey articles~\cite{Th,Th1} that went in-depth beyond the isoparametric category, I would therefore devote this article primarily to the classification part of the hypersurface case. 

As in~\cite{Chi6}, I once more followed the historical development of the extensive studies to let the flow of presentation as motivated and seamless as possible. The difference from~\cite{Chi6} is that I included much more detailed expositions and proofs in the present article. Since the subject is so all-encompassing that I was compelled to assume background knowledge in the first place, a comfortable understanding of symmetric spaces is preferred, for which I would refer to the two volumes by Loos~\cite{Lo}, Volume 2 of Kobayashi and Nomizu~\cite{KN}, and Helgason~\cite{He}. Meanwhile, in~\cite{Chi4}, I wrote a fairly detailed exposition on the comprehensive commutative algebra engaged in the isoparametric story, to which I would thus refer without dwelling more on it than is necessary.
Also, I would only report on the methods entailed in the classification, in the case of four principal curvatures, done by my collaborators and myself~\cite{CCJ}, and subsequently by myself~\cite{Ch1, Ch2, Chi3}, leaving the classification in the case of six principal curvatures~\cite{DN,Mi2,Mi4} to the expertise of the authors themselves.

The codimension 2 estimates prevalent in the classification derived and developed from the powerful criterion of Serre on prime ideals is stressed in this article, whose unifying capacity in conjunction with the underlying isoparametric geometry lifts us from the jungle of intertwined components of tensors to the canopy of ideals of polynomial rings, to enable us to see the entire landscape of classification.

I would like to thank the referees for many valuable comments to better the exposition, and Zhenxiao Xie for his careful reading through the manuscript during his visit at Washington University.

\section{The dawn, 1918-1940} When we stroll along a beach, the last arriving wavefront gently brushes our feet to a halt, where we often see cusps forming of the wavefront due to the different speeds at front points. Waves form singularity. The same phenomenon, if applied to the lenses of our eyes shone upon by lights refracted through different media with different resulting speeds, most of us will feel disoriented after some exposure time, because the formation of wave singularity plays the trick. It would then be interesting to understand the wavefronts that travel at a constant speed each moment. This was investigated by Laura~\cite{La} in 1918. He concluded that such wavefronts were either planar, cylindrical, or spherical, as our daily experiences would almost certainly convince us that this is the case, by seeing laser beams of planar wavefronts, fluorescent tube light of cylindrical wavefronts, and candle light of spherical wavefronts. 

To see the equations that govern the wavefronts traveling at a constant speed each moment, let us start with the wave equation,
$$
\Delta \phi = \frac{\partial^2 \phi}{\partial t^2}.
$$
Wavefronts are {\em level surfaces} of $\phi$, at each moment $t$,
which propagate along the normal directions of the level surfaces. That the wavefront speed remains constant on each level surface means $ds/dt=b(\phi)$ and
\begin{equation}\label{nabla}
|\nabla \phi|=\; {\rm change\; per\; unit\; length\; of}\;\phi\;{\rm along\; the\;normal}=a(\phi),
\end{equation}
for some smooth functions $a$ and $b$, where $s$ is the distance by which a wavefront travels.
Therefore, we derive
$$
\frac{\partial \phi}{\partial t}=\frac{\partial \phi}{\partial s}\frac{ds}{dt}=a(\phi)b(\phi):=c(\phi),
$$
and, as a result,
\begin{equation}\label{laplace}
\Delta \phi=\frac{\partial^2 \phi}{\partial t^2}=c'(\phi)c(\phi).
\end{equation}
\begin{definition}\label{hyper} A smooth function $f$ over a Riemannian manifold is transnormal
if
$$
|\nabla f|=A(f)
$$
for some smooth function $A$. A transnormal function is isoparametric if, for some smooth function $B$,
$$
\Delta f=B(f).
$$
Let $c$ be a regular value of an isoparametric function $f$.
The level surface $f^{-1}(c)$ is called
an isoparametric hypersurface. 
\end{definition}

Somigliana's paper (1918-1919) brought mean curvature to the foreground.

\vspace{2mm}

\begin{theorem} {\rm (Somigliana,~\cite{So})} Over ${\mathbb R}^3$, a transnormal
function $f$ is isoparametric
if and only if each regular level surface of $f$ has constant mean curvature.
\end{theorem}
\begin{proof} For each regular level surface
$M:=f^{-1}(c)$
of a transnormal function $f$,
$$
{\bf n}=\nabla f/|\nabla f|=\nabla f/A(f)
$$
\vspace{1mm}

\noindent is a unit normal field to $M$. The {\em shape operator} $S$ of the surface $M$ is
\vspace{2mm}
$S(X):=-d{\bf n}(X)$
for a tangent vector $X$ of $M$. However,
\begin{eqnarray}\nonumber
\aligned
&d(\nabla f)(X)=d(A(f){\bf n})(X)=A'(f)\,df(X)\,{\bf n}+A(f)\,d{\bf n}(X)\\
&=A(f)\,d{\bf n}(X)=-A(f)\,S(X).
\endaligned
\end{eqnarray}
On the other hand, as a linear operator,
$$
d(\nabla f): X\mapsto \text{Hessian}(f)\, X.
$$
%where $\text{Hessian}(f)=(\partial^2_{ij} f)$. 
Taking trace, we obtain
\begin{eqnarray}\nonumber
\aligned
&\Delta f=trace(d(\nabla f))=-A(f)\, trace(S)+<d(\nabla f)({\bf n}),{\bf n}>\\
&=-2A(f)\, H +A'(f)\,A(f), 
\endaligned
\end{eqnarray}
where $H$ is the mean curvature of $M$. In other words,
$$
H=-(B(f)-A'(f)A(f))/2A(f)
$$
is a constant along $M$ if the transnormal $f$ is also isoparametric.

Conversely, if $H$ is constant along
regular level surfaces of the transnormal $f$, then $H$ is a function of $f$ and so
$\Delta f$ is a function of $f$ so that $f$ is isoparametric.
\end{proof}
In particular, he arrived at the same conclusion that over ${\mathbb R}^3$, the regular level surfaces of an isoparametric
function are either all spheres, all cylinders or all planes.

This theorem was rediscovered later by Segre in {\rm 1924}~\cite{Se1} and Levi-Civita
in {\rm 1937}~\cite{Le}. The approach Levi-Civita gave is what we will look at next.
\begin{theorem} {\rm (Levi-Civita,~\cite{Le})} Over ${\mathbb R}^3$, a transnormal $f$ is isoparametric if and only if
the two principal curvatures of each regular level surface are constant.
\end{theorem}

\begin{proof} Observe first that the integral curves of
the unit normal field ${\bf n}=\nabla f/|\nabla f|$ are just line segments.
In fact, an integral curve $c$ of ${\bf n}$ from $f=a$ to $f=b$ assumes
the length
$$
{\rm Length\; of}\;c=\int_a^b \frac{df}{|\nabla f|}=\int_a^b \frac{df}{A(f)}.
$$
On the other hand, for any curve $\gamma(t), 0\leq t\leq 1,$ beginning and ending
at the two end points of the given integral curve, we have
$$
|\frac{df(\gamma(t))}{dt}| = |\langle\nabla f(\gamma(t)),\gamma'(t)\rangle|\leq A(f(\gamma(t))|\gamma'(t)|,
$$
so that
$$
{\rm Length\; of}\;\gamma=\int_0^1 |\gamma'(t)|dt\geq \int_0^1 \frac{1}{A(f)}\frac{df}{dt}dt=\int_a^b \frac{df}{A(f)}.
$$
In other words, the given integral curve $c$ assumes the shortest distance
among all curves beginning and ending at its end points, i.e., the integral curve
is a line segment. In view of this observation, instead of using $f$ to parametrize the level surfaces,
we might as well use the arc length $s$ of the normal lines of an initial level
surface to parametrize other level surfaces, so that 
$$
M_s:= M+s{\bf n}
$$
is now the 1-parameter family of level surfaces of the transnormal $f$, where
$M$ is the initial level surface with unit normal field ${\bf n}$.

Let us calculate the mean curvature $H_s$ of $M_s$ by using the fact that
${\bf n}$ is still the unit normal to $M_s$. The upshot is~\cite[p. 209]{On}
$$
H_s= \frac{H-sK}{1-2sH+s^2K},
$$
where $k_1$ and $k_2$ are the eigenvalues (the principal curvatures)
of the shape operator
$S$ of $M$, so that $H=(k_1+k_2)/2$, and $K=k_1k_2$ is the Gaussian curvature
of $M$.
Therefore, the mean curvature $H_s$ is constant on $M_s$ for all $s$, i.e.,
the transnormal $f$ is isoparametric, if and only if
$H,K$ are constant on $M$,
if and only if the principal curvatures $k_1,k_2$ are constant.

Case 1. $k_1=k_2\neq 0$. $M$ is a sphere.

Case 2. $k_1=k_2=0$. $M$ is a plane.

Case 3. $k_1\neq k_2$. One employs $dk_1=dk_2=0$ and a bit more surface geometry
to conclude $k_1k_2=0$~\cite[p. 255]{On}, so that one of $k_1,k_2$ is zero. Then $M$ is a cylinder.
\end{proof}

Segre then took up the investigation of generalizing the question to ${\mathbb R}^n$ in 1938, followed by Cartan's look into the hyperbolic space $H^n$ in the same year. 
\begin{theorem}\label{S} {\rm (Segre,~\cite{Se2})} The same conclusion holds on ${\mathbb R}^n$, namely,
an isoparametric hypersurface, which is a regular level hypersurface of an isoparametric
function $f$ over ${\mathbb R}^n$ satisfying
$$
|\nabla f|=A(f),\quad \Delta f=B(f),
$$
is either a hypersphere, a hyperplane, both are totally umbilic {\rm(}with one principal curvature{\rm )},
or a cylinder $S^k\times{\mathbb R}^{n-1-k}$.
\end{theorem}

\begin{theorem}\label{C} {\rm (Cartan,~\cite{Ca1})} The same conclusion holds on
the hyperbolic space $H^n$ of constant curvature $-1$, namely, an isoparametric hypersurface
in $H^n$ is either a sphere, a hyperbolic $H^{n-1}$, a Euclidean ${\mathbb R}^{n-1}$
{\rm(}i.e., a horosphere{\rm )}, all three being totally umbilic,
or a cylinder $S^k\times H^{n-k-1}$.
\end{theorem}
      
\begin{proof} (sketch) We show again that there are at most two (constant) principal curvatures of the shape operator. 
Indeed, let $\lambda_1,\cdots,\lambda_{n-1}$
be the principal curvatures of an isoparametric hypersurface in a standard
space form of dimension $n$ with constant curvature $C$. Then we have 
\begin{equation}\label{fund}
\sum_{j\neq k}m_j\frac{C+\lambda_k\lambda_j}{\lambda_k-\lambda_j}=0,\quad {\rm summed\; on}\; j,
\end{equation}
referred to by Cartan as the ``Fundamental Formula'', which was proven by Segre in the Euclidean case
and by Cartan in general~\cite[p. 84]{BCO}. Here, $m_j$ is the multiplicity
of $\lambda_j$
and $\lambda_i\neq\lambda_j$ if $i\neq j$.

\vspace{2mm}

Case 1, $C=0$. Let $\lambda_k$ be the smallest positive principal curvature. Note that
each term, if nontrivial, in the fundamental formula must be negative, which is
a contradiction. Therefore, there are at most two principal curvatures, one of
them is zero if there are two.

\vspace{2mm}

Case 2, $C=-1$. It is easy to see that 

\begin{equation}\label{-1}
\frac{C+\lambda_k\lambda_j}{\lambda_k-\lambda_j}<0
\end{equation}
if $\lambda_j\leq 0$ and $\lambda_k>0$. Consider those positive principal values. If there is a $0<\lambda_l\leq 1$
such that $(\lambda_j)^{-1}\leq\lambda_l$ for all $\lambda_j>1$, we let $\lambda_k$
be the largest positive principal value $\leq 1$.
It follows that~\eqref{-1} is negative for all those $0<\lambda_j<1$ and for those
$\lambda_j>1$ not reciprocal to $\lambda_k$. We conclude that none of the positive $\lambda_j$
other than the reciprocal of $\lambda_k$ exist. 
Otherwise, there exists some $\lambda_j>1$
such that its reciprocal is greater than the above $\lambda_k$, which
we replace by the smallest principal value $>1$. Once more,~\eqref{-1} is negative
for all positive $\lambda_j$ not reciprocal to $\lambda_k$. We arrive at the
same conclusion as in the preceding case. So,
we have at most two principal curvatures reciprocal to each other.

In the case of two distinct principal curvatures $\lambda$ and $\mu$, the isoparametric hypersurface is the product of two simply connected
space forms of constant curvatures $\lambda^2-1$ and $\mu^2-1$.  All these hypersurfaces are homogeneous.

See~\cite[Theorem 2.5, p. 373]{Ry} on the product structure in both cases.
\end{proof}

The story now takes a fascinating turn when Cartan directed his attention to the spherical case. Let us denote by $g$ the number of principal curvatures of an isoparametric hypersurface in the sphere $S^n$. He quickly settled the cases when $g\leq 2$. For $g=1$, the hypersurface is any sphere in the 1-parameter family of such spheres perpendicular to the axis of the North and South poles of $S^n$; note that the 1-parameter family degenerates to the South and North poles. For $g=2$,  the hypersurface is any one in the 1-parameter family  of the product of two subspheres of the form $S^k\times S^{n-k-1}$ whose points are given in coordinates as $(x_0,\cdots,x_{k},x_{k+1},\cdots,x_n)$, where
$$
x_0^2+\cdots+x_k^2=r^2,\quad
\;x_{k+1}^2+\cdots+x_n^2=s^2,\quad r^2+s^2=1,
$$
which degenerates to the two manifolds $S^k$ and $S^{n-k-1}$ of radius $1$
as $r$ approaches $0$ or $1$. We refer to~\cite{Ry} for a proof. All these hypersurfaces are homogeneous.

Then Cartan worked on the case $g=3$. {\em A priori}, the three principal values could potentially carry distinct multiplicities, which was ruled out by him. He went on to classify and was amazed by the beauty of such hypersurfaces, as one can tell from the title of his 1939 papers~\cite{Ca2},~\cite{Ca3}:
\begin{theorem}\label{CC} Let $g=3$.
\begin{itemize}
\item[I.] Over the ambient Euclidean space ${\mathbb R}^{n+1}\supset S^n$, there
is a homogeneous polynomial $F$ of degree $3$, satisfying
\begin{eqnarray}\label{CMP}
\aligned
&|\nabla F|^2=9r^2,\;r=|x| \; \text{for}\;x\in{\mathbb R}^{n+1},\\
&\Delta F=0,
\endaligned
\end{eqnarray}
whose restriction to $S^n$ is exactly the isoparametric function $f$.
The range of $f$
is $[-1,1]$ with $\pm 1$ the only critical values. Thus $f^{-1}(c),-1<c<1,$
form a $1$-parameter family of isoparametric hypersurfaces that degenerates to the
two critical sets $f^{-1}(1)$ and $f^{-1}(-1)$.
\item[II.] The three principal values have equal multiplicity $m=1,2,4$,
or $8$.
\item[III.] The two critical sets of $f$ are the real, complex, quaternionic,
or octonion projective plane corresponding to the principal multiplicity $m=1,2,4,$ or
$8$. Each isoparametric hypersurface in the family is a tube around the projective plane.
\item[IV.] Let ${\mathbb F}$ be one of the normed algebras ${\mathbb R},
{\mathbb C},{\mathbb H}$, and ${\mathbb O}$. Let $X,Y,Z\in{\mathbb F}$ and
$a,b\in{\mathbb R}$. Then

\begin{eqnarray}\label{Cartan}
\aligned
&F=a^3-3ab^2+\frac{3a}{2}(X\overline{X}+Y\overline{Y}-2Z\overline{Z})\\
&+\frac{3\sqrt{3}b}{2}(X\overline{X}-Y\overline{Y})
+\frac{3\sqrt{3}}{2}((XY)Z+\overline{(XY)Z}).
\endaligned
\end{eqnarray}
\end{itemize}
\end{theorem}
In particular, all these hypersurfaces are homogeneous. %He was very pleased with the geometric realization of the Cayley projective plane through the isoparametric way~\cite{Ca2},~\cite{Ca3}. 

His proof was an algebraic analysis of the homogeneous polynomial $F$ of degree 3, nowadays called the Cartan-M\"{u}nzner polynomial (see the next section for more expositions on it). He expanded the polynomial at a point $p$ of, say, the critical (or, focal) manifold $f^{-1}(1)$ with the coordinate
$u$ parametrizing ${\mathbb R}p$ and $x_1,\cdots,x_n$ parametrizing $({\mathbb R}p)^\perp$. He observed that the equations~\eqref{CMP} dictated that
$$
F=u^3 +P(x)u +Q(x),
$$
where $P(x)$ is quadratic and $Q(x)$ is cubic homogeneous in $x$. Substituting this into the equations in~\eqref{CMP}, he found that there was an integer $m$ such that
$$
P=\frac{3}{2}(x_1^2+\cdots+x_{2m}^2) -3(x_{2m+1}^2+ \cdots +x_{3m+1}^2), \quad n=3m+1.
$$
In fact, 
$$
w_0:=x_{2m+1},\cdots,w_m:=x_{3m+1}
$$ 
constitute the normal coordinates to the focal manifold in $S^n$ at $p$. It follows that the isoparametric hypersurface is of equal multiplicity $m$ by identifying any of the three principal spaces of its shape operator
with the normal space to a point on the associated focal manifold. 

He then turned to $Q(x)$, which he expanded into
 $$
 Q=A_3+A_2+A_1+A_0,
 $$
 where $A_i$ is of degree $i$ in $x_1,\cdots,x_{2m}$ and degree $3-i$ in $w_0,\cdots,w_m$ and deduced that $A_3=A_1=A_0=0$ so that
 $$
 Q=Q_0z_0+\cdots+Q_m z_m
 $$ 
 for some homogeneous quadratic polynomials $Q_0, \cdots, Q_m$ in $x_1,\cdots,x_{2m}$ that are in fact the components of the second fundamental form of the focal manifold at $p$. Since the principal curvatures of the focal manifold are $\cot(\pi/3)$ and $\cot(2\pi/3)$ (see~\eqref{focalcurvature} below for details), only differing by a sign, he diagonalized $Q_0$ so that
 $$
 Q_0=\frac{3\sqrt{3}}{2}(x_1^2+\cdots+x_m^2-y_1^2-\cdots -y_m^2),
 $$
by reindexing the coordinates $y_1:=x_{m+1},\cdots,y_m:=x_{2m}$. He then set $H_a:=Q_a/3\sqrt{3}$ for $1\leq a\leq m$ and concluded
$$
H_1^2+\cdots+H_m^2 = (x_1^2+\cdots+x_m^2)(y_1^2+\cdots+y_m^2),
$$
or, in other words, the map 
$$
{\mathcal L}:(x,y)\in {\mathbb R}^m\times{\mathbb R}^m\mapsto (H_1(x,y), \cdots,H_m(x,y))\in{\mathbb R}^m
$$
satisfies $|{\mathcal L}(x,y)|=|x||y|$, i.e., it is an orthogonal multiplication~\cite{Hu}, from which he deduced that $m=1, 2, 4,$ or, $8$. 

With the expressions of $P$ and $Q$ pinned down, he set $v:=w_0$ and wrote down
$$
\aligned
&F=u^3-3uv^2+\frac{3}{2}u\sum_{i=1}^m (x_i^2+y_i^2)-3u\sum_{i=1}^m w_i^2\\
&+\frac{3\sqrt{3}}{2}v\sum_{i=1}^m (x_i^2-y_i^2)+\sum_{i=1}^m w_i\, Q_i(x,y).
\endaligned
$$
This is exactly~\eqref{Cartan} when we employ the orthogonal multiplications ${\mathcal L}$, which give rise to the products of the four normed algebras.

Thanks to~\cite{Ka},~\cite{KK},~\cite{CO}, we have geometric alternatives to the proof, which were inspired, directly or indirectly, by M\"{u}nzner's fundamental papers~\cite{Mu} that entered the stage in the early 1970s. 

We will derive~\eqref{Cartan} in Subsection 4.3 for the real case, using a Lie-theoretic approach.

Next, Cartan worked on the case $g=4$~\cite{Ca4}. The situation got much more complicated. So, he assumed equal multiplicities $m$ to the four principal values and pointed out that the hypersurface satisfies the equations
$$
\aligned
&|\nabla F|^2=16r^2,\;r=|x| \; \text{for}\;x\in{\mathbb R}^{n+1},\\
&\Delta F=0.
\endaligned
$$
He indicated the classification without proof in the cases when $m=$ 1, or 2, by writing down the respective 4th degree Cartan-M\"{u}nzner polynomials. When $m=1$, the polynomial is~\eqref{nomizu} below for $k=3$, whereas
when $m=2$, the polynomial is defined over the Euclidean space $so(5,{\mathbb R})$, where an element $Z\in so(5,{\mathbb R})$ is written as
$$\begin{pmatrix} Z_1&Z_2&Z_3&Z_4&Z_5\end{pmatrix}$$
by five column vectors, for which the polynomial is
\begin{equation}\label{F(Z)}
F(Z):= \frac{5}{4}\sum_i\,|Z_i|^4-\frac{3}{2}\sum_{i<j}\,|Z_i|^2|Z_j|^2+4\sum_{i<j}\,(\langle Z_i,Z_j\rangle)^2,
\end{equation}
with respect to which he studied geometric properties of these two hypersurfaces. Both are homogeneous.
It is quite evident that in the paper he had in mind the isotropy representations of the rank-2 symmetric spaces $SO(5)\times SO(5)/\Delta(SO(5)),$ $\Delta$ denoting the diagonal map, and $SO(10)/U(5)$, respectively.

He left us with three natural questions:

\vspace{2mm}

\noindent  {\rm ({Cartan, 1940~\cite{Ca3})}
\begin{itemize}
\item[(i)] What are the possible $g${\rm ?}
\item[(ii)] Is equal multiplicity of principal curvatures always true{\rm ?}
\item[(iii)] Are all isoparametric hypersurfaces homogeneous{\rm ?}
\end{itemize}

Before proceeding, let us make the following slightly more general definition of an isoparametric hypersurface than the one in Definition~\ref{hyper}, though they turn out to be equivalent through M\"{u}nzner's work.

\begin{definition} 
A hypersurface in $S^n$ is isoparametric if all its principal curvatures are constant with fixed multiplicities. Let $M$ be an isoparametric hypersurface in $S^n$. The number
of principal curvatures of $M$ is denoted by $g$.
\end{definition}

\section{The dormancy, 1941-1969} When I was a student at Stanford, I had a few times of working as a TA for David Gilbarg. During a chat he said that he had been trained as an algebraist working under Emil Artin. But war had a better dictation over everything. Upon graduation he was designated to work on partial differential equations in the wartime, from that point on it was a point of no return for him. Mathematics more on the ``pure'' side essentially ground to a halt.

Diversified activities awakened and resumed after the war. One of the most notable was the vibrancy in the field of algebraic and differential topology, which led to the profound discoveries in the 1950-60s, from which the stage was set for the revival of the isoparametric story. At the same time, algebraic geometry also experienced a gigantic transformation in which the theory of schemes via ideals changed its landscape.

\section{The Renaissance, 1970-1979} 
\subsection{Nomizu, Takagi and Takahashi's work on homogeneous isoparametric hypersurfaces} Nomizu wrote two papers to start this period in the
early 1970s, in which he constructed examples~\cite{No1} that answered the second question of Cartan in the negative and reviewed Cartan's work~\cite{No2}. 

Indeed, consider ${\mathbb C}^k={\mathbb R}^k\oplus{\mathbb R}^k$ and write
$z\in{\mathbb C}^k$
as $z=x+\sqrt{-1}y$ accordingly. Define a homogeneous polynomial of degree $4$
on ${\mathbb C}^k$ by
$$
\tilde{F}=(|x|^2-|y|^2)^2+4(\langle x,y\rangle)^2.
$$
Then $\tilde{f}:=F|_{S^{2k-1}}$ is an isoparametric function whose regular level sets form a 1-parameter family of isoparametric hypersurfaces with four principal values and multiplicities $\{1, k-2\}$. 

The isoparametric hypersurfaces (respectively, all the above isoparametric hypersurfaces by Cartan) are the principal isotropy orbits of the symmetric spaces
$SO(k+2)/S(2)\times SO(k)$ (respectively, of appropriate symmetric spaces of rank 2 to be seen in Subsection 4.3). 

Note that $\tilde{f}$ has range $[0,1]$. So, we normalize
it by defining $f:=1-2\tilde{f}$, or rather, by setting
\begin{equation}\label{nomizu}
F:=(|x|^2+|y|^2)^2-2\tilde{F}.
\end{equation}
$F$ is an isoparametric function such that $f$ has range $[-1,1]$. 

Note also that 
$$
f^{-1}(1):=\{(x,y): |x|^2=|y|^2=1/2, \langle x,y\rangle=0\},
$$
the Stiefel manifold of oriented 2-frames, which is orientable, whereas 
$$
f^{-1}(-1)=\{z\in {\mathbb C}^k: z=e^{\sqrt{-1}\theta} \, v,\; v\in S^{k-1}, \theta\in [0,2\pi)\},
$$
which is doubly covered by $S^1\times S^{k-1}$ and is not orientable when $k$ is odd, as was pointed out by Cartan in~\cite{Ca4} in the case of multiplicity pair $(1,1)$ for $k=3$.

%In hind sight, one cannot but wonder why Cartan, being the father of symmetric spaces, did not look at the isotropy representations for more examples of isoparametric hypersurfaces.

At about the same time Takagi and Takahashi~\cite{TT} classified homogeneous isoparametric hypersurfaces
in spheres, which are the principal orbits of isotropy representations of simply connected symmetric spaces of rank 2. They also calculated the number $g$ of principal curvatures to be 1, 2, 3, 4, or 6, and, moreover, verified that there were at most two distinct multiplicities, in the homogeneous category. We begin with a definition.

\begin{definition} A connected hypersurface $M$ in a smooth Riemannian manifold $X$
is homogeneous if 
$I(X,M)$, the group of isometries of $X$ leaving $M$ invariant, acts transitively
on $M$.
\end{definition}

It is clear that, for such a hypersurface, the principal curvatures of its shape
operator are everywhere constant, counting multiplicities. Hence, they are isoparametric.

%\begin{definition} A hypersurfaces in ${\mathbb R}^n,S^n$ or $H^n,$ is
%called {\bf isoparametric} if its principal curvatures are everywhere constant, counting
%multiplicities.
%\end{definition}

Theorems~\ref{S} and~\ref{C} classify all isoparametric hypersurfaces in ${\mathbb R}^n$
and $H^n$ to be exactly the homogeneous hypersurfaces in these space forms.
What is interesting is then the spherical case. 

%\begin{question} Classify all isoparametric hypersurface

Recall a representation $\rho: G\hookrightarrow SO(n+1)$ acting
on ${\mathbb R}^{n+1}$ is {\em effective} if every nontrivial element in $G$ displaces some vector in ${\mathbb R}^{n+1}$. Now,
let $I(M)$ be the group of isometries
of $M$ and let $\iota:I(S^n,M)\rightarrow I(M)$ be the restriction map. Let
$I_0(S^n,M)$ be the connected component of the identity of $I(S^n,M)$ and let
$G:=\iota(I_0(S^n,M))$.

\begin{proposition}\label{kn}~\cite[II, p. 15]{OT}
$\iota:I_0(S^n,M)\rightarrow G$
is an isomorphism, so that $\iota^{-1}:G\hookrightarrow SO(n+1)$ is an effective
representation on ${\mathbb R}^{n+1}$ with $M$ an orbit.
Furthermore, $M$ is compact, and so in particular $G$ is compact and hence is a Lie group.
\end{proposition}

\begin{definition} An effective representation $\rho: G\hookrightarrow SO(n+1)$ acting
on ${\mathbb R}^{n+1}$ is of cohomogeneity $r$ if the smallest codimension
of all orbits of $\rho$ is $r$ in ${\mathbb R}^{n+1}$.
\end{definition}

In particular, the representation $\iota$ above of a homogeneous
hypersurface in $S^n$ is of cohomogeneity 2.

\begin{definition} Given an effective representation
$\rho:G\hookrightarrow SO(n+1)$ of cohomogeneity $r$,
$\rho$ is maximal if there is no effective representation
$\rho_1:G_1:\hookrightarrow
SO(n+1)$ of cohomogeneity $r$ such that $G$ is a proper subgroup of $G_1$
with $\rho(g)=\rho_1(g)$ for all $g\in G$.
\end{definition}

\vspace{1mm}

\begin{proposition}~\cite[p. 16]{OT}
\begin{itemize}
\item[(1)] The effective representation $\iota:G\hookrightarrow SO(n+1)$
in Proposition{\rm ~\ref{kn}}
is a maximal effective representation of cohomogeneity $2$.
\item[(2)] Let $\rho:G\hookrightarrow SO(n+1)$ be a
maximal effective representation of cohomogeneity $2$. Let $M$ be a $G$-orbit of
codimension $2$ in ${\mathbb R}^{n+1}$. Then $\rho(G)=I_0(S^n,M)$.
\item[(3)] In particular, any maximal effective representation
$\rho:G\hookrightarrow SO(n+1)$
is obtained as the representation of a homogeneous hypersurface in $S^n$.
\item[(4)] Two homogeneous hypersurfaces $M$ and $N$ in $S^n$ are equivalent,
i.e., $N=f(M)$ for an $f\in O(n+1)$, if and only if $I_0(S^n,M)\simeq I_0(S^n,N)$
through the isomorphism $g\mapsto fgf^{-1}$.
\end{itemize}
\end{proposition}

Therefore, the classification of homogeneous hypersurfaces in $S^n$ is equivalent to
first classifying maximal effective orthogonal representations $\rho:G\hookrightarrow SO(n+1)$
of cohomogeneity 2 and then classifying their orbits of codimension 2.
Hsiang and Lawson classified all maximal orthogonal representations in~\cite{HL}.
They are closely tied with what are called the $s$-representations of symmetric spaces. We will return to this in Subsection 4.3.
 
%We next report on Takagi and Takahashi's work, which is based on the comprehensive work
%of Cartan on the classification of symmetric spaces.

\subsection{M\"{u}nzner's work on the general case} M\"{u}nzner~\cite{Mu} (preprint circulating in 1973)
established a breakthrough result that developed Cartan's work, recorded
in Theorem~\ref{CC}, in a far-reaching manner:

\begin{theorem} %{\rm (M\"{u}nzner's structure theory on isoparametric hypersurfaces)}
Let $M$ be any isoparametric hypersurface with $g$ principal curvatures in
$S^n$. 
Then we have the following.
\begin{itemize}
\item[(1)] There is a homogeneous polynomial $F$, called the Cartan-M\"{u}nzner
polynomial, of degree $g$ over ${\mathbb R}^{n+1}$ satisfying
\begin{equation}\label{CM}
|\nabla F|^2=g^2r^{2g-2},\quad \Delta F=\frac{m_{-}-m_{+}}{2}g^2r^{g-2},
\end{equation}
where $r$ is the radial function over ${\mathbb R}^{n+1}$.
\item[(2)] Let $f:=F|_{S^n}$. Then the range of $f$ is $[-1,1]$. The only critical
values of $f$ are $\pm 1$. Moreover, $M_{\pm}:=f^{-1}(\pm 1)$ are connected submanifolds
of codimension $m_{\pm}+1$ in $S^n$, called focal manifolds, whose principal
curvatures are $\cot(k\pi/g),1\leq k\leq g-1$.
\item[(3)] For any $c\in(-1,1)$, $f^{-1}(c)$ is an isoparametric hypersurface
with at most two multiplicities $m_{\pm}$ associated with the principal curvatures.
In fact, if we order the principal curvatures $\lambda_1>\cdots>\lambda_g$
with multiplicities $m_1,\cdots,m_g$, then $m_i=m_{i+2}$ with index modulo $g$;
in particular, all multiplicities are equal when $g$ is odd, and
when $g$ is even, there are at most
two multiplicities equal to $m_{\pm}$.
\item[(4)]
Each of the $1$-parameter isoparametric hypersurfaces is a tube around each of the
two focal manifolds, so that $S^n$ is obtained by gluing two disk
bundles over $M_{\pm}$ along the isoparametric hypersurface
$M_0:=f^{-1}(0)$. As a consequence, algebraic topology implies that 
the only possible values of $g$ are $1,2,3,4,$ or $6$.
\end{itemize}
\end{theorem}

Indeed, start with an isoparametric hypersurface
$$
x:M\hookrightarrow S^n
$$
whose principal curvatures are set to be
$$
\lambda_j=\cot(\theta_j), \quad 0<\theta_1<\cdots<\theta_g<\pi,
$$
with respect to the outward normal field $n$. Let us consider
the parallel transport of $M$,
\begin{equation}\label{xt}
x_t:=\cos(t) x+\sin(t) n,
\end{equation}
which is the counterpart to the Euclidean parallel transport
along the normal direction. {\em A priori}, $M_t:=x_t(M)$ is an embedding
for small $t$. Since
$$
n_t:=-\sin(t) x+\cos(t) n
$$
is normal to $M_t$, a straightforward calculation derives that the principal
curvatures of $M_t$, with respect to the chosen normal field $n_t$, are

\begin{eqnarray}\label{symm}
\lambda_j(t) =\cot (\theta_j-t)
\end{eqnarray}
with the same eigenspace and multiplicity as $\lambda_j$. On the other hand,
for a fixed $l$, the eigenspace of $\lambda_l$ from point to point defines an integrable
distribution, called the $l$-th curvature distribution, on $M$
with spheres of radius $|\sin(\theta_l)|$ as leaves. This can be directly
checked by differentiating
$$
f_l(x):=x+v_l(x)/|v_l(x)|^2, \quad v_l(x):= -x+\cot(\theta_l) n,
$$
to see that $f_l(x)$ is a constant $c_l$ on the $l$-th curvature leaf through $x$; we have
\begin{equation}\label{c}
c_l=\cos(\theta_l)(\cos(\theta_l)x+\sin(\theta_l)n),
\end{equation}
i.e., the unit vector pointing in the same direction as $c_l$ assumes the angle $\theta_l$ on the unit circle
oriented from $x$ to $n$. 
Now that the curvature leaf through $x$ is a sphere of
radius $|\sin(\theta_l)|$ 
centered at $c_l$, the antipodal point to $x$ on this leaf gives the reflection
map $\phi_l$ about $c_l$:
$$
\phi_l(x):=x+2v_l(x)/|v_l(x)|^2=\cos(2\theta_l)x+\sin(2\theta_l) n,
$$
i.e., $\phi(x)$ is the point of reflection of $x$ about the line spanned by $c_l$
on the $(x,n)$-plane.
Therefore, by~\eqref{symm}, the principal curvatures of $M$ at $\phi_l(x)$ are
\begin{equation}\label{anti}
-\cot(\theta_j-2\theta_l),\quad 1\leq j\leq g,
\end{equation}
with the same eigenspaces and multiplicities
as $x$. Note that the sign in~\eqref{anti} differs from that in~\eqref{symm},
because the circle $x_t$ leaves $M$ at $x$ and enters $M$ at $\phi(x)$, so that
$n_{2\theta_l}$ at $\phi(x)$ is negative of the chosen outward normal field $n$
of $M$ at $\phi(x)$. 
Since $M$ has constant principal curvatures, counting multiplicities,
we conclude that the following sets
\begin{equation}\label{reflection}
\{\cot(\theta_j)\}, \quad \{\cot(2\theta_l-\theta_j)\}
\end{equation}
are identical for all $j,l$, and two numbers, one from each set, having the same index $j$
have the same principal multiplicity $m_j$, regardless of what $l$ is.

Now,~\eqref{reflection} means that the lines $L_j$ spanned by $c_j$ on the
$(x,n)$-plane, all through the origin, satisfies the property that the reflection of $L_j$ about
any $L_l$ is another $L_k$. It follows that these lines $L_j,1\leq j\leq g,$
are equally spaced in the $(x,n)$-plane so that
$$
\theta_j=(j-1)\pi/g+\theta_1.
$$
Thus, the reflections about the lines $L_j$ result in $m_i=m_{i+2}$ with
index modulo $g$. Accordingly, we denote $m_1$ and $m_2$ by $m_{+}$ and
$m_{-}$, respectively. (This is reminiscent of a root system and its Weyl chambers.)

Having done so, M\"{u}nzner went on to construct a
local isoparametric function, which is nothing but an appropriate distance function,
in a neighborhood of $M$, already observed by Cartan for $g=3$, as follows. Any $p$ in a tubular neighborhood $U$
of $M$ can be written uniquely as
$$
p= \cos(t)x +\sin(t) n
$$
for some small $t$. Define
$$
\mu(p):=\theta_1-t,\quad V(p)=\cos(g\,\mu(p)).
$$
Extend $V(p)$ to a neighborhood of $M$ in the ambient Euclidean space by defining
$$
F(rp)=r^gV(p),
$$
where $r$ is the Euclidean radial function. 

\begin{theorem}\label{thm} $F$ is in fact a homogeneous polynomial of degree $g$ satisfying
$$
|\nabla F|^2=g^2r^{2g-2},\quad \Delta F=g^2\frac{m_{-}-m_{+}}{2}.
$$
\end{theorem}

\begin{proof} (Sketch) Define
$$
G:=F-ar^{g},
$$
where
$$
a:=\frac{g}{g+n-1}\frac{m_{-}-m_{+}}{2}.
$$
Then verify that
$$
\Delta G=0.
$$
In general, it is true that for a harmonic function $G$
over ${\mathbb R}^{n+1}$, we have
\begin{equation}\label{ng}
\Delta^g|\nabla G|^2=\sum (\partial^{g+1} G/\partial x_{i_1}\cdots\partial x_{i_{g+1}})^2.
\end{equation}
On the other hand, for the $G$ engaged in our consideration, a calculation gives
$$
|\nabla G|^2=g^2r^{2g-2}(1+a^2)-2ag^2r^{g-2}F.
$$
We therefore find
$$
\Delta^{g-1}|\nabla G|^2=c
$$
with $c$ an appropriate constant. $F$ is thus a homogeneous polynomial
by~\eqref{ng}.
\end{proof}

Now that $F$ is globally analytic over ${\mathbb R}^{n+1}$,
we set $f:=F|_{S^n}$.
A calculation by the formulae
$$
|\nabla F|^2=(\frac{\partial F}{\partial r})^2+|\nabla f|^2,\quad
\Delta f=\Delta F-\frac{\partial^2 F}{\partial r^2} -n\frac{\partial F}{\partial r}
$$
derives, by Theorem~\ref{thm}, that
$$
|\nabla f|^2=A(f),\quad \Delta f=B(f),
$$
where
$$
A(f)=g^2(1-f^2),\quad B(f)=-g(n+g-1)f+\frac{m_{-}-m_{+}}{2}g^2.
$$
So, $f$ is an isoparametric function on $S^n$. Note that $A(f)=0$ only at $f=\pm 1$,
so that the range of $f$ is $[-1,1]$ and $\pm 1$ are the only critical values.
Let $M_{\pm}:=f^{-1}(\pm 1)$ be the singular set. $S^n\setminus (M_{+}\cup M_{-})$
is open and dense and is diffeomorphic to $M_c\times (-1,1)$ for any fixed $c$, where
$M_c:=f^{-1}(c)$ for $c\in(-1,1)$.

{\em A priori}, $M_c$ might not be
connected. We claim that this is not the case.
Define
$$
d:M\times (0,\pi/g)\rightarrow S^n,\quad d(x,\mu)=\cos(\theta_1-\mu) x+\sin(\theta_1-\mu) n.
$$
Then
$$
f(d(x,\mu))=\cos(g\mu)
$$
by the analytic nature of $f$ because $f|_U=V$ and the identity holds
on $U$; in particular, $M$ is contained in $M_c$ with $c=\cos(g\theta_1)$. But then
the map
\begin{equation}\label{diffe}
d_c:M_c\times (0,\pi/g)\rightarrow S^n,\quad d_c(x,\mu)=\cos(\theta_1-\mu)x+\sin(\theta_1-\mu)n
\end{equation}
also satisfies $f(d_c(x,\mu))=\cos(g\mu)$ and
$$
d_c:M_c\times (0,\pi/g)\rightarrow S^n\setminus (M_{+}\cup M_{-})
$$
is a diffeomorphism. From this we see that the map 
$$
h:M_c\rightarrow S^n,\quad x\mapsto \cos(\theta_1)x+\sin(\theta_1)n
$$
maps $M_c$ to $M_{+}$. Observe that $h(x)$ points in the same direction as $c_1$ for the curvature leaf
through $x$ whose tangent space is the eigenspace with
principal value $\cot(\theta_1)$, where $c_1$ is
defined in~\eqref{c}.
It follows that $h:M_c\rightarrow M_{+}$ is a sphere bundle whose fiber is a
curvature leaf diffeomorphic to $S^{m_{+}}$.  Meanwhile, it is easy to check
that $dh$ has kernel dimension $m_{+}$; at $x$, the derivative $dh$
preserves eigenspaces of all principal values other than that of
$\cot(\theta_1)$. Therefore, $M_{+}$ is a manifold of dimension $\dim(M)-m_{+}$,
which is of codimension at least $2$ in $S^n$. Likewise, the codimension
of $M_{-}$ is at least $2$ in $S^n$.

Returning to the map~\eqref{diffe}, we see now $S^n\setminus(M_{+}\cup M_{-})$
is connected as $M_{+}$ and $M_{-}$ are of codimension at least $2$ in $S^n$.
Therefore, that $d_c$ is a diffeomorphism ensures that $M_c$ is connected, for all $c$.
As a consequence, $M_{\pm}$ are also connected via the map $h$.

Lastly, since $h(x)=x_{\theta_1}$ defined in~\eqref{xt}, we see by~\eqref{symm}
that the principal values of $M_{+}$, in any normal direction, are
\begin{equation}\label{focalcurvature}
\cot(\theta_j-\theta_1)=\cot((j-1)\pi/g), \quad 2\leq j\leq g.
\end{equation}
This also holds true for $M_{-}$.

We remark that the sphere bundle property via the map $h$ holds true over any complete Riemannian manifold when we only assume transnormality~\cite{Wa1}, though the focal manifolds need not be connected in general.
\begin{corollary} $M_{\pm}$ are minimal submanifolds of $S^n$. The minimality
condition is exactly equation~\eqref{fund}, the fundamental formula of Segre and Cartan, when $C=1$.
\end{corollary}

\begin{proof} By the preceding formula, the mean curvature of $M_{+}$ in any normal direction is
$$
\sum_{j=1}^{g-1}\cot(j\pi/g)=0,
$$
which is the fundamental formula~\eqref{fund}.
\end{proof}

\begin{corollary}\label{minimal} There is a unique minimal isoparametric hypersurface in the $1$-parameter
family $M_t$.
\end{corollary}

\begin{proof} By~\eqref{symm}, the mean curvature of $M_t$ is
$$
H:=\sum_{j=1}^g\cot(\theta_j-t)
$$
for $t\in(0,\pi/g)$. $H$ is strictly increasing as the derivative is $>0$.
Near $t=\theta_1<\pi/g$ the function is $>0$ whereas near $\theta=\pi/g-\theta_1>0$
the function is $<0$. Therefore, there is a unique $t\in(0,\pi/g)$ at which
$H=0$.
\end{proof}
Now that $S^n$ is obtained by gluing two disk bundles over the focal manifolds
$M_{\pm}$
along an isoparametric hypersurface $M$, M\"{u}nzner used algebraic topology
to express the cohomology ring of $M$, with ${\mathbb Z}_2$ coefficients, as
modules of those of $M_{\pm}$, whose intertwining module structures via Steenrod squares
then give the remarkably short list that $g=1,2,3,4,$ or $6$! Thus, it answered the first question of Cartan.

\subsection{The homogeneous case in retrospect} Now that we have a grand view of the structure of isoparametric hypersurfaces thanks to  M\"{u}nzner's theorem, let us return to the homogeneous case, this time with more geometric insight.
\begin{definition} An {\bf $s$-representation} of rank $r$ is the isotropy
representation of a connected, simply connected semisimple
Riemannian symmetric space
of rank $r$. Here, if the symmetric space is decomposed into its irreducible components,
the rank is the sum of the ranks of the components.
\end{definition}

An $s$-representation of rank 2 is either the isotropy representation
of two irreducible symmetric spaces of rank 1,
or of an irreducible symmetric space of rank 2. Note that
${\mathbb R}\times M,$ where $M$ is irreducible symmetric of rank 1, is
also of rank 2, although its isotropy representation is not an
$s$-representation.

In connection with classifying homogeneous hypersurfaces in $S^n$, we are particularly
interested in the isotropy representations of simply connected noncompact
symmetric spaces of rank 2, because of the theorem of Hsiang and Lawson~\cite{HL} on
the classification of
all maximal effective orthogonal representations
$\rho:G\hookrightarrow SO(n+1)$ of cohomogeneity 2:

\begin{theorem}\label{hl} Up to equivalence,
the maximal effective orthogonal representations
$\rho:G\hookrightarrow SO(n+1)$ of cohomogeneity $2$ are exactly the isotropy
representations of the  simply connected noncompact symmetric spaces of
rank $2$, i.e., the isotropy representations of
\begin{itemize}
\item[(1)] ${\mathbb R}\times H^n$, where the principal orbits
are spheres $S^{n-1}\subset S^n\subset{\mathbb R}^{n+1}$,
\item[(2)] $H^p\times H^q$, where the principal orbits are
$S^{p-1}\times S^{q-1}\subset
S^{p+q-1}\subset{\mathbb R}^{p+q}$, and
\item[(3)] the noncompact irreducible symmetric spaces of rank $2$,
where principal orbits are those of $s$-representations.
\end{itemize}
\end{theorem}

More generally, Dadok's classification~\cite{D} shows that any polar representation is orbit equivalent to an $s$-representation. See also~\cite{EH},~\cite{EH1} for a more conceptual proof.

The principal orbits of the first two items in Theorem~\ref{hl} are easy to visualize.
Item (1) is the isotropy representation
of $O_0(1,n)$ on the hyperbolic space $H^n$ identified with $O_0(1,n)/SO(n)$ as a symmetric space, where $O_0(1,n)$ is the identity component of $O(1,n)$ and $SO(n)$ is identified with
$$
K=\{\begin{pmatrix}1&0\\0&A\end{pmatrix},\quad A\in SO(n)\}.
$$
Thus, it gives the standard orthogonal representation $SO(n)$ on ${\mathbb R}^n$, whose
typical 
principal orbit
is the sphere $S^{n-1}$. The Euclidean factor in item (1) 
acts trivially, so that a principal orbit of codimension 2 of the isotropy
representation is
$S^{n-1}\subset S^n\subset{\mathbb R}^{n+1}$.
In the same vein, a typical principal orbit of the isotropy representation
in item (2) is $S^{p-1}\times S^{q-1}\subset
S^{p+q-1}\subset{\mathbb R}^{p+q}$.

In particular, we have $g=1$ or $2$ for the homogeneous spaces in the first two
items of Theorem~\ref{hl}. The isoparametric hypersurfaces with $g=1$ or $2$,
classified by Cartan, are
exactly the ones in the first two items.

Let us study item (3) in Theorem~\ref{hl}, where the principal orbits give rise to all homogeneous
isoparametric hypersurfaces with $g\geq 3$ in the sphere.

Let $G/K$ be a noncompact irreducible symmetric space of
rank $2$ (which is automatically simply connected) with the Cartan decomposition ${\mathcal G}={\mathcal K}\oplus{\mathcal M}$.
Fix a $v\neq 0\in{\mathcal M}$.
We know~\cite[p. 247]{He} there is a $k\in K$ such that
$Ad(k)\cdot v\in{\mathcal A}$,
where ${\mathcal A}$ is the maximal abelian subspace of ${\mathcal M}$. Therefore,
we may assume without loss of generality that $v\in{\mathcal A}$.

\begin{proposition}
With the setup above, an orbit $Ad(K)\cdot v$, where $v\in{\mathcal A},$ is
principal of codimension $2$ if and only if $v$ lies in a Weyl chamber.
\end{proposition}

\begin{proof} The isotropy subgroup $L$ of $Ad(K)$ leaving $v$ fixed has the
Lie algebra
$$
{\mathcal L}:=\{X\in{\mathcal K}:[X,v]=0\}.
$$
We have the root space decomposition
\begin{equation}\label{Weyl}
0=ad_h(Z)=\sum_{\lambda\in\Sigma} \lambda(h) Z_\lambda
\end{equation}
with $h\in{\mathcal A}$ and $Z_\lambda\in{\mathcal G}_\lambda$, and so
$$
{\mathcal G}={\mathcal N}_0\oplus\sum_{\lambda\in\Sigma}{\mathcal G}_\lambda,
$$
where ${\mathcal N}_0$ is the centralizer of ${\mathcal A}$ in ${\mathcal G}$.

If $v$ belongs to a Weyl chamber, then
$\lambda(v)\neq 0$ for all $\lambda\in\Sigma$, so that by~\eqref{Weyl}
$X_\lambda=0$ for all $\lambda\in\Sigma$.
In other words, $X\in{\mathcal L}$ if and only if $X\in {\mathcal M}_0:={\mathcal N}_0\cap {\mathcal K}$,
the centralizer of ${\mathcal A}$ in ${\mathcal K}$.

By~\cite[Lemma 3.6, p. 261]{He}, we know ${\mathcal M}_0$ has the same codimension
in ${\mathcal K}$ as ${\mathcal A}$ in ${\mathcal M}$, i.e.,
$$
\dim(Ad(K)/L)=\dim({\mathcal K}) -\dim({\mathcal L})=\dim({\mathcal M})-\dim({\mathcal A}).
$$
Thus, the isotropy orbit is of codimension $\dim({\mathcal A})=2$.
                                                                
If $v$ lies in a chamber wall, then by~\eqref{Weyl} $X\in{\mathcal L}$
if and only if
$$
X\in{\mathcal M}_0\oplus \sum_{\lambda,\lambda(v)\neq 0} ({\mathcal G}_\lambda\cap{\mathcal K}).
$$
Therefore, the codimension of the orbit of $v$ is larger than 2. 
\end{proof}

\begin{corollary} 
The isotropy representation of an irreducible noncompact
symmetric space of rank $2$ has only two singular orbits and a $1$-parameter family
of diffeomorphic principal orbits of codimension $2$ degenerating to the two singular
orbits.
\end{corollary}

\begin{proof} In the rank 2 case, a Weyl chamber is a sector of the plane of angle measure
$\pi/3$ for $A_2$, $\pi/4$ for $B_2$, and $\pi/6$ for $G_2$. Let us say
$\theta_0<\theta<\theta_0+\pi/l, l=3,4,6,$ defines the chamber. Then the preceding
proposition says that for any unit $v$ assuming angle $\theta$ in the chamber, its isotropic
orbit is homogeneous (and hence isoparametric) of codimension 2 and
is diffeomorphic to $Ad(K)/L$. So, we have a 1-parameter family of
diffeomorphic homogeneous isoparametric hypersurfaces. At the two chamber walls, i.e., when $v$ assumes the angle $\theta_0$ or $\theta_0+\pi/l$, the dimension of
the orbit drops. Meanwhile, since the normalizer of
${\mathcal A}$ serves as the Weyl group by Theorem~\ref{norm} below, we see that the isotropic orbit
of $v$ intersecting the chamber plane at some points $v_1=v,v_2,\cdots,v_{2l}$,
one in each chamber. So the isotropic representation has only two
singular orbits, even though there are $2g$ Weyl chambers.
All other orbits are principal of codimension 2.
\end{proof}

\begin{proposition} With the same setup, $Ad(k)({\mathcal A})$ is the normal plane to the principal orbit $Ad(K)(v)$
at $Ad(k)(v)$ for $v\in{\mathcal A}$.
\end{proposition}

\begin{proof} It suffices to check this at $v$, where the tangent space
of the orbit is
$$
T_v=\{[h,v]:h\in{\mathcal K}\}.
$$
But then for $w\in{\mathcal A}$, we have 
$$
\langle w,[h,v]\rangle=\langle [v,w],h\rangle=0,
$$
since the inner product is proportional to the Killing form,
\end{proof}

\begin{proposition} With the same setup, let $w$
be a unit vector perpendicular to $v$ in ${\mathcal A}$, and extend it to a global
normal field on the principal orbit $Ad(K)\cdot v,|v|=1,$ in the unit sphere
of ${\mathcal M}$. The shape operator $S_w$ of the orbit
at $v$ satisfies that the eigenvalues are
$$
-\lambda(w)/\lambda(v),
$$
where $\lambda$ are reduced positive roots such that $\lambda/2\notin\Sigma$.
The eigenspace associated with the above eigenvalue is
$$
E_\lambda={\mathcal G}_\lambda\oplus{\mathcal G}_{-\lambda}
\oplus{\mathcal G}_{2\lambda}\oplus{\mathcal G}_{-2\lambda}.
$$
In particular, $g$, the number of principal curvatures of the shape operator, is $3$, $4$, or $6$.
If we label the principal curvatures by $\lambda_1>\cdots>\lambda_g$ and their multiplicities by
$m_1,\cdots, m_g$, then $m_i=m_{i+2}$, where the subscripts are modulo $g$. In particular,
the multiplicities are all equal when $g=3$. Moreover, if we choose the angles
$$
\theta_i=(2i-1)\pi/2g,\quad i=1,\cdots, g,
$$
to coordinatize the positive roots, then the principal curvatures are
$$
\lambda_i=\tan(\theta-\theta_i),\quad -\pi/g<\theta<\pi/g,
$$
when $v$ assumes the angle $\theta$ and $w$ the angle $\theta+\pi/2$.
\end{proposition}

\begin{proof} As mentioned in the preceding proposition, a vector $X$ tangent
to the orbit is of the form
$$
X=[k,v]=-\sum_{\lambda\in\Sigma} \lambda(v) X_\lambda
$$
for $k\in{\mathcal K}$. 
Since
$$
n=Ad(K)\cdot w
$$
is a normal vector field to the orbit, the shape operator is
$$
S(X):=-dn(X)=-[X,w]=\sum_{\lambda\in\Sigma} \lambda(w) X_\lambda.
$$
Therefore, we obtain
$$
-\lambda(v)S(X_\lambda)=\lambda(w) X_\lambda.
$$
Since $v$ is regular we have $\lambda(v)\neq 0$ for all $\lambda\in\Sigma$. It
follows that
$$
S(X_\lambda)=-\lambda(w)/\lambda(v) X_\lambda.
$$
The principal curvatures of $S$ are thus $-\lambda(w)/\lambda(v)$, attained
by $\pm\lambda,\pm2\lambda$, so that the eigenspace $E_\lambda$ with the
principal curvature $-\lambda(w)/\lambda(v)$ is  the desired one, where $\lambda$
need only go through the positive roots $\lambda$ for which $\lambda/2\notin\Sigma$,
which form a reduced root system.
The number of positive roots in the $A_2,B_2$, or $G_2$ root system
is $3$, $4$, or $6$, respectively, which is $g$.

We choose the angles $\theta_i=(2i-1)\pi/2g$ to coordinatize the positive roots. 
We see the Weyl group is generated by
\begin{equation}\label{generator}
\theta\mapsto\pi/g-\theta,\quad \theta\mapsto \theta+2\pi/g.
\end{equation}
By Theorem~\ref{norm}, the Weyl group preserves the principal curvatures
and their multiplicities. Hence, $m_i=m_{i+2}$ with index modulo $g$.

Lastly, since
$$
v=(\cos(\theta),\sin(\theta)), \quad w=(-\sin(\theta),\cos(\theta)),\quad \lambda_i=(\cos(\theta_i),\sin(\theta_i)),
$$
we calculate
$$
-\lambda_i(w)/\lambda_i(v)=-\langle w,\lambda_i\rangle/\langle v,\lambda_i\rangle=\tan(\theta-\theta_i).
$$
\end{proof}

Let us now look at
$$
F(\theta):=\sin(g\theta),\quad -\pi/2g<\theta<\pi/2g.
$$
It is left invariant by the two generators of the Weyl group
in~\eqref{generator}. In fact, $F(\theta)$ is the restriction to the unit circle
of the homogeneous polynomial of degree $g$
\begin{equation}\label{A}
F_{\mathcal A}:=\sum_{0}^{[(g-1)/2]}{g\choose 2i+1}(-1)^ix^{g-(2i+1)}y^{2i+1}
\end{equation}
defined by the maximal abelian space ${\mathcal A}$. $F_{\mathcal A}$ is left
invariant by the Weyl group.

\begin{theorem}~\cite[p. 299]{KN} The space of homogeneous
polynomials on ${\mathcal M}$ left invariant by $Ad(K)$ is isomorphic to the space
of homogeneous polynomials on ${\mathcal A}$ left invariant by the
Weyl group.
\end{theorem}

This theorem is called Chevalley Restriction Theorem. In~\cite{KN}, the proof
is given for a compact Lie group, or for a symmetric space of Type II. But the proof there can be modified easily
to arrive at the preceding theorem in view of a characterization of the Weyl group:
\begin{theorem}\label{norm}~\cite[p. 284, p. 289]{He} Let $(G,K,\sigma)$ be an irreducible Riemannian symmetric space
of noncompact type. Let
$$                                                        
M:=\{k\in K:Ad(k)\cdot v=v, \forall v\in{\mathcal A}\},
\quad M':=\{k\in K:Ad(k)\cdot{\mathcal A}\subset{\mathcal A}\}.
$$
Then $M'/M$ is the Weyl group.
\end{theorem}

The last ingredient for constructing the Cartan-M\"{u}nzner polynomial in the homogeneous category is the following theorem of Chevalley.

\vspace{2mm}

\begin{theorem}~\cite{Ch} The space of homogeneous polynomials on a
maximal abelian space ${\mathcal A}$ of dimension $r$ left invariant by the Weyl group
is generated by $r$ algebraically independent polynomials.
\end{theorem}

Since $r=2$ in our case and we have found two generators, namely, $x^2+y^2$
and $F_{\mathcal A}$ on ${\mathcal A}$, the space of homogeneous
polynomials left invariant by $Ad(K)$ on ${\mathcal M}$ of dimension $n$
is thus generated by $(x_1)^2+\cdots+(x_n)^2$ and a homogeneous polynomial $F$ of degree
$g$ whose restriction to the circle is $F_{\mathcal A}$.

$F$, homogeneous of degree $g$, thus leaves each isotropic orbit invariant.
Therefore, we conclude the following.

\begin{theorem}\label{CMM} There is a homogeneous polynomial $F$ of degree $g$, called
the Cartan-M\"{u}nzner polynomial,
for $g=3,4,6$, on ${\mathcal M}$, whose restriction $f$ to the unit sphere
of ${\mathcal M}$ satisfies the property that its range is $[-1,1]$. For
each $c\in(-1,1)$, $f^{-1}(c)$ is a homogeneous (isoparametric) hypersurface
degenerating to two singular submanifolds $f^{-1}(\pm 1)$. The statement is clearly
true when $g=1$ or $2$. All homogeneous hypersurfaces in spheres are constructed this way.
\end{theorem}

We remark that for $g=1$ in the preceding theorem, the polynomial is
$F=x_{n+1}$ over ${\mathbb R}^{n+1}$, while for $g=2$ the polynomial is
$F=(x_1)^2+\cdots+(x_r)^2-(x_{r+1})^2-\cdots-(x_{r+s})^2$ over ${\mathbb R}^{r+s}$.

\vspace{2mm}

\noindent {\bf Example.} We find $F$ in the case $g=3$ when the symmetric space is
$SU(3)/SO(3)$ of Type I and rank 2.

Let ${\mathcal M}$ be the space of 5-dimensional $3$ by $3$ real traceless symmetric matrices. The Cartan decomposition
is
$$
su(3)=so(3)\oplus \sqrt{-1}{\mathcal M},\quad {\mathcal K}=so(3).
$$
${\mathcal M}$ is equipped with the inner product
$$
\langle Y,Y\rangle:=tr(YY)=
\alpha^2+\beta^2+\gamma^2+x^2+y^2+z^2,
$$
which is a multiple of the Killing form of $su(3)$ restricted to ${\mathcal M}$,
where we write
$$
Y:=\begin{pmatrix} \alpha&x/\sqrt{2}&y/\sqrt{2}\\x/\sqrt{2}&\beta&z/\sqrt{2}\\y/\sqrt{2}&z/\sqrt{2}&\gamma\end{pmatrix},\quad \alpha+\beta+\gamma=0.
$$
The isotropic action is the adjoint action
$$
Ad(T):V\in {\mathcal M} \mapsto TVT^{-1}\in{\mathcal M},\quad T\in SO(3).
$$
The diagonal block of ${\mathcal M}$ is the maximal abelian subspace ${\mathcal A}$
of ${\mathcal M}$. The three positive roots are
$$
\alpha_1:=\begin{pmatrix}1&0&0\\0&-1&0\\0&0&0\end{pmatrix}/\sqrt{2},
\quad\alpha_2:=\begin{pmatrix}0&0&0\\0&1&0\\0&0&-1\end{pmatrix}/\sqrt{2},\quad
\alpha_3:=\begin{pmatrix}1&0&0\\0&0&0\\0&0&-1\end{pmatrix}/\sqrt{2},
$$
where $\alpha_1$ and $\alpha_2$ are simple roots. We choose
the unit angle bisector as the standard basis element
$$
e_1:=(2\alpha_1+\alpha_2)/\sqrt{6},
$$
and
$$
e_2:=\alpha_2.
$$
Then $e_1,e_2$ form an orthonormal basis of ${\mathcal A}$, so that an element
in ${\mathcal M}$ relative to $e_1,e_2$ is
$$
X:=ae_1+be_2=\begin{pmatrix}2a/\sqrt{6}&0&0\\0&b/\sqrt{2}-a/\sqrt{6}&0\\0&0&-b/\sqrt{2}-a/\sqrt{6}\end{pmatrix},
$$
and the $Y$ above is
$$
Y:=\begin{pmatrix}2a/\sqrt{6}&x/\sqrt{2}&y/\sqrt{2}\\x/\sqrt{2}&b/\sqrt{2}-a/\sqrt{6}&z/\sqrt{2}\\y/\sqrt{2}&z/\sqrt{2}&-b/\sqrt{2}-a/\sqrt{6}\end{pmatrix}.
$$

Now, it is clear that $\det(Y)$ is $Ad(SO(3))$-invariant. We calculate
$$
\det(X)=\frac{1}{3\sqrt{6}}(a^3-3ab^2)=\frac{1}{3\sqrt{6}}F_{\mathcal A},
$$
where $F_{\mathcal A}$ is given in~\eqref{A}, when we set
$$
a=\cos(\pi/6-\theta),\quad b=\sin(\pi/6-\theta),\quad -\pi/6<\theta<\pi/6.
$$
It follows that
$$
F:=3\sqrt{6}\det(Y)
$$
restricts to $F_{\mathcal A}$ and so
$F$ is the Cartan polynomial given in Theorem~{\rm \ref{CMM}}. A calculation
shows $F$ is exactly the polynomial given in~\eqref{Cartan} by Cartan
in the case when ${\mathbb F}$ is ${\mathbb R}$.

Note that $f$, the restriction of $F$ to the unit sphere, has range $[-1,1]$
and $f^{-1}(\pm 1)$ are the two singular submanifolds, both being the
projective plane. To see this, we set $\theta=\pm\pi/6$. Then, respectively,
$$
X=\begin{pmatrix}2/\sqrt{6}&0&0\\0&-1/\sqrt{6}&0\\
0&0&-1/\sqrt{6}\end{pmatrix},\quad
\begin{pmatrix}1/\sqrt{6}&0&0\\0&1/\sqrt{6}&0\\
0&0&-2/\sqrt{6}\end{pmatrix}.
$$
Let us find the isotropy group $L$ of the isotropy action on $X$, where
$L$ consists of all $T\in SO(3)$ commuting with $X$. We see $L$ is in diagonal block form.
Hence,
$$
L\simeq S(O(1)\times O(2)),
$$
so that the singular orbits are
$$
Ad(SO(3))/L=SO(3)/S(O(1)\times O(2))={\mathbb R}P^2.
$$

A look at the tables for the symmetric spaces of rank $2$
of Types I and II
shows that there are four such spaces with $g=3$, which are
$$
\aligned
&SU(3)/SO(3),\quad SU(3)\times SU(3)/\Delta (SU(3)\times SU(3)),\\
&SU(6)/Sp(3),\quad E_6/F_4,
\endaligned
$$
whose Cartan polynomials of their isotropic orbits are the ones given
in~\eqref{Cartan}.

As in the $SU(3)/SO(3)$ case, the singular orbits of the other three examples
are, respectively, the complex, quaternionic and octonion projective planes.
The principal orbits are tubes around the projective planes.

The following grid table is the collection of all symmetric spaces $G/K$ of Types I and II whose
isotropy representations give homogeneous isoparametric hypersurfaces $M$. There are at most
two multiplicities $(m_{+},m_{-})$ for the $g$ principal curvatures.

\begin{center}
\begin{tabular}{|l|l|l|l|l|}
\hline
$G$&$K$&$\dim M$&$g$&$(m_{+},m_{-})$\\
\hline
$S^1\times SO(n+1)$&$SO(n)$&$n$&1&$(1,1)$\\
\hline
$SO(p+1)\times SO(n+1-p)$&$SO(p)\times SO(n-p)$&$n$&2&$(p,n-p)$\\
\hline
$SU(3)$&$SO(3)$&3&3&$(1,1)$\\
\hline
$SU(3)\times SU(3)$&$SU(3)$&6&3&$(2,2)$\\
\hline
$SU(6)$&$Sp(3)$&12&3&$(4,4)$\\
\hline
$E_6$&$F_4$&24&3&$(8,8)$\\
\hline
$SO(5)\times SO(5)$&$SO(5)$&8&4&$(2,2)$\\
\hline
$SO(10)$&$U(5)$&18&4&$(4,5)$\\
\hline
$SO(m+2),m\geq 3$&$SO(m)\times SO(2)$&$2m-2$&4&$(1,m-2)$\\
\hline
$SU(m+2),m\geq 2$&$S(U(m)\times U(2))$&$4m-2$&4&$(2,2m-3)$\\
\hline
$Sp(m+2),m\geq 2$&$Sp(m)\times Sp(2)$&$8m-2$&4&$(4,4m-5)$\\
\hline
$E_6$&$(Spin(10)\times SO(2))/{\mathbb Z}_4$&30&4&$(6,9)$\\
\hline
$G_2$&$SO(4)$&6&6&$(1,1)$\\
\hline
$G_2\times G_2$&$G_2$&12&6&$(2,2)$\\
\hline
\end{tabular}
\end{center}

\subsection{Ozeki and Takeuchi's work on inhomogeneous examples for $g=4$} Based on M\"{u}nzner's work, Ozeki and Takeuchi~\cite[I]{OT} constructed two classes,
each with infinitely many members, of 
inhomogeneous isoparametric hypersurfaces with $g=4$. This answered Cartan's third question in the negative.

They also classified all isoparametric hypersurfaces with $g=4$ when one of the multiplicities is $2$, which are all homogeneous~\cite[II]{OT}.

An important ingredient in their work is their expansion formula of the
Cartan-M\"{u}nzner polynomial, which was inspired by Cartan's approach to the classification for $g=3$ mentioned in Section 2. The central theme is to study the focal manifolds, the singular set of the 1-parameter family of isoparametric hypersurfaces, to recover properties of the hypersurface. They fixed a point $x$ on either one of the focal manifolds, say, $M_{+}$, and decomposed the ambient Euclidean space by ${\mathbb R}x$,  with coordinate $t$, the tangent space to $M$ at $x$ in $S^n$, where a typical vector is denoted by $y$, and the normal space to $M_{+}$ at $x$ in $S^{n}$, where a typical vector is denoted by $w$ with coordinates $w_i$ with respect to a chosen
orthonormal basis ${\bf n}_0,{\bf n}_1,\cdots,{\bf n}_{m_{+}}$. They expanded the Cartan-M\"{u}nzner polynomial in $t$, with undetermined coefficients in $y$ and $w$, and substituted it into the two equations in~\eqref{CM} to result in
\begin{eqnarray}\label{exp-f}
\aligned
&F(tx+y+w)=t^4+(2|y|^2-6|w|^2)t^2+8(\sum_{a=0}^{m_{+}}p_{a}w_{a})t\\
&+|y|^4-6|y|^2|w|^2+|w|^4-2\sum_{a=0}^{m_{+}}(p_{a})^2
+8\sum_{a=0}^{m_{+}}q_{a}w_{a}
\\
&+2\sum_{a,b=0}^{m_{+}}\langle\nabla p_{a},\nabla p_{b}\rangle w_{a}w_{b},
\endaligned
\end{eqnarray}
where $p_a(y)$ (respectively, $q_a(y)$) is the $a$-th component
of the $2$nd (respectively, $3$rd) fundamental form of $M_{+}$ at $x$.
Furthermore,
$p_a$ and $q_a$ are subject to ten defining equations~\cite[I, pp 529-530]{OT}, of which the first three
assert that, since
$S_{\bf n}$, the $2$nd fundamental matrix of $M_{+}$ in any unit normal
direction ${\bf n}$, has eigenvalues $1,-1,0$ with fixed multiplicities, it
must be that
 \begin{equation}\label{S^3}
 (S_{\bf n})^3=S_{\bf n}
 \end{equation}
 for all $n$. We will return to the ten identities later. 

Ozeki and Takeuchi then introduced {\bf Condition A}, where a point $p$ on a focal manifold, say $M_{+}$, is of Condition A if the shape operators $S_n$  
at $p$ share the same kernel for all $n$. 

Indeed, given a normal basis $n_0,\cdots,n_{m_{+}}$ at $x$ with the associated shape operators $S_0,\cdots, S_{m_{+}}$, let $E_0, E_{1}, E_{-1}$ be the eigenspaces of $S_0$ with eigenvalues $0, 1, -1$, respectively. Relative to the chosen orthonormal bases of $E_0,E_{1},E_{-1}$, we have the matrix representations of the shape operators 
\begin{equation}\label{0a}
S_0=\begin{pmatrix}Id&0&0\\0&-Id&0\\0&0&0\end{pmatrix},\quad
S_a=\begin{pmatrix}0&A_a&B_a\\A_a^{tr}&0&C_a\\B_a^{tr}&C_a^{tr}
&0\end{pmatrix}
\end{equation}
for $1\leq a\leq m_{+},$ where $A_a:E_{-1}\rightarrow E_{1}$,
$B_a:E_0\rightarrow E_{1}$ and $C_a:E_0\rightarrow E_{-1}$.

Condition $A$ means that $B_a=C_a=0$ for all $1\leq a\leq m_{+}$ at $x$. Consider $n:=(n_a+n_b)/\sqrt{2}$ and substitute $S_n$ into~\eqref{S^3} we can extract
\begin{equation}\label{AAtr}
A_aA_b^{tr}+A_bA_a^{tr}=2\delta_{ab}\, Id,\quad A_a^{tr}A_b+A_b^{tr}A_a=2\delta_{ab}.
\end{equation}
This implies that the symmetric matrices
\begin{equation}\label{irre}
T_0:=\begin{pmatrix}I&0\\0&-I\end{pmatrix},\quad T_a:=\begin{pmatrix}0&A_a\\A_a^{tr}&0\end{pmatrix}, \quad 1\leq a\leq m_{+},
\end{equation}
induce a symmetric Clifford $C'_{1+m_{+}}$-module structure on ${\mathbb R}^{2m_{-}}$, so that ${\mathbb R}^{2m_{-}}$ is decomposed into $k$ irreducible $C'_{1+m_{+}}$ modules for some $k$, and thus
$$
2m_{-}= k\,\theta_{m_{+}}=2k\delta_{m_{+}},
$$
where, following the standard notation, $\theta_s$ denotes the dimension of an irreducible $C'_{s+1}$-module while $\delta_s$ denotes that of an irreducible skew-symmetric Clifford $C_{s-1}$-module satisfying the periodicity
$\delta_{s+8}=16\,\delta_s$ with $\delta_1,\cdots,\delta_8$ being 1, 2, 4, 4, 8, 8, 8, 8, respectively. See Definition~\ref{dfn} and what immediately follows it in the next section for a more detailed account.

They focused on the case when $\delta_{1+m_{+}}=1+m_{+}$, i.e., the case when the dimension of an irreducible module of $C_{m_{+}}$ is $1+m_{+}$. It is well known that $m_{+}=$ 1, 3, or 7, which they verified by the above periodicity formula for $\delta_s$. As a consequence $m_{-}$ is a multiple of 2, 4, 8, respectively. Now, complex multiplication on ${\mathbb C}$ give the only irreducible representation of $C_1$, whereas left and right quaternionic or octonion multiplications give the only two distinct irreducible $C_{m_{+}}$-representations in the cases $m_{+}=$ 3, or 7, respectively. Each such irreducible skew-symmetric Clifford representation $A_1,\cdots, A_{m_{+}}$ reconstructs a symmetric Clifford representation $C'_{1+m_{+}}$ via~\eqref{irre}, and vice versa. Putting $k$ such symmetric Clifford representations together, we recover the second fundamental form, given in~\eqref{0a}, of $M_{+}$ of codimension $1+m_{+}=2, 4,$ or 8 in $S^n$.

Having determined the second fundamental form, they introduced Condition B. For each tangent vector $y$ decomposed into $y=y_0+y_1+y_{-1}$ relative to the eigenspaces $E_0,E_1,E_{-1}$ of the shape operator $S_0$, they stipulated that 
$$
q_a(y) =\sum_{b=1}^{m_{+}} r_{ab}(y_0)\, p_b(y), \quad 1\leq a,b\leq m_{+}
$$
for some 1st-degree polynomials $r_{ab}$ linear in $y_0$, to be a candidate for the polynomial in~\eqref{exp-f}. They substituted Condition B into the ten defining equations of an isoparametric hypersurface with four principal curvatures and eventually determined the third fundamental form, and hence the Cartan-M\"{u}zner polynomial. Since the multiplicity pair $(m_{+}, m_{-})$ does not appear on the list of those of homogeneous examples when $m_{+}$ is 3 or 7, they found the first examples of inhomogeneous isoparametric hypersurfaces with four distinct principal curvatures  and multiplicity pair $(3,4k)$ and $(7,8k)$ in the sphere.

Conditions A and B would play a major role in later development.

We note that Takagi~\cite{Ta} classified the case when $g=4$ and one of the multiplicities is 1. They are congruent to the aforementioned examples of Nomizu, and hence are homogeneous. His method is to also expand the Cartan-M\"{u}nzner polynomial at a focal point and analyze the algebraic structures constrained by~\eqref{CM}.

\section{The enlightenment, 1980-1999} 
\subsection{1980-1989} Now that the two questions of Cartan were answered, the next question naturally came down to the possible multiplicity pairs $(m_{+},m_{-})$ of the $g$ principal curvatures. The 1980s started with the paper of Ferus, Karcher, and M\"{u}nzner that constructed an infinite 2-dimensional array of multiplicity pairs each of which is associated with an isoparametric hypersurface, most of them inhomogeneous, which include the examples of Ozeki and Takeuchi.

To motivate their work, let us return to the examples of Nomizu~\eqref{nomizu}. Set
$$
P_0:=\begin{pmatrix}I&0\\0&-I\end{pmatrix},\quad
P_1:=\begin{pmatrix}0&I\\I&0\end{pmatrix}, \quad u:=(x,y)^{tr},
$$
where $I$ is the $k$ by $k$ identity matrix. Then $F$ can be rewritten as
$$
F=|u|^4-2\sum_{i=0}^{1} \langle P_i u,u\rangle^2,\quad P_iP_j+P_jP_i=2\delta_{ij}I.
$$
Ferus, Karcher and M\"{u}nzner's construction is a far-reaching generalization of this.

\begin{definition}\label{dfn} The skew-symmetric (respectively, symmetric) Clifford algebra
$C_n$ (respectively, $C_n'$) over ${\mathbb R}$ is the algebra
generated by the standard basis $e_1,\cdots,e_n$ of ${\mathbb R}^n$ subject to the only constraint
$$
e_ie_j+e_je_i=-2\delta_{ij}I\quad ({\rm respectively},\;\; e_ie_j+e_je_i=2\delta_{ij}I).
$$
\end{definition}
The classification of the Clifford algebras is known~\cite{Hu}:

\vspace{3mm}

\begin{center}
\begin{tabular}{|c|c|c|c|c|c|c|c|c|}
\hline
$n$&1&2&3&4&5&6&7&8\\
\hline
$C_n$&${\mathbb C}$&${\mathbb H}$&${\mathbb H}\oplus {\mathbb H}$&${\mathbb H}(2)$&${\mathbb C}(4)$&${\mathbb R}(8)$&${\mathbb R}(8)\oplus {\mathbb R}(8)$&${\mathbb R}(16)$\\
\hline
$\delta_n$&1&2&4&4&8&8&8&8\\
\hline
$C_n'$&${\mathbb R}\oplus {\mathbb R}$&${\mathbb R}(2)$&${\mathbb C}(2)$&${\mathbb H}(2)$&${\mathbb H}(2)\oplus {\mathbb H}(2)$&${\mathbb H}(4)$&${\mathbb C}(8)$&${\mathbb R}(16)$\\
\hline
$\theta_n$&2&4&8&8&16&16&16&16\\
\hline
\end{tabular}
\end{center}

\vspace{3mm}

\noindent Here, $\delta_n$ is the dimension 
of an irreducible module of $C_{n-1}$, and $\theta_n$ is the dimension of an irreducible module of $C_{n+1}'$.
Moreover, $C_n$ (respectively, $C_n'$) is subject to the periodicity
condition $C_{n+8}=C_n\otimes {\mathbb R}(16)$ (respectively, 
$C_{n+8}'=C_n'\otimes {\mathbb R}(16)$).
The generators $e_1,\cdots,e_n$ acting on
each irreducible module of either $C_n$ or $C_n'$ in the grid table give rise to
$n$ skew-symmetric or symmetric orthogonal matrices
$T_1,\cdots,T_n$ satisfying 
$$
T_iT_j+T_jT_i=\pm 2\delta_{ij}\,Id,
$$ 
yielding a representation of $C_n$ or $C_n'$ on the irreducible module, respectively. Note that we have
$$
\theta_n=2\delta_n.
$$
This is not coincidental. It says that we can construct symmetric representations of $C_{m+1}'$ from skew-symmetric representations
of $C_{m-1}$, and vice versa. Indeed, let us be given $k$
irreducible representations $V_1,\cdots,V_k$ of $C_{m-1}$. Set
$$
V:=V_1\oplus+\cdots\oplus V_k\simeq {\mathbb R}^l,\quad l=k\delta_{m}.
$$
The representations of $e_1,\cdots,e_{m-1}$ on $V_1,\cdots,V_k$ result in
$m-1$ skew-symmetric orthogonal matrices $E_1,\cdots,E_{m-1}$
on $V$. Set
$$
P_0:=\begin{pmatrix}I&0\\0&-I\end{pmatrix},\quad P_1:=\begin{pmatrix}0&I\\I&0\end{pmatrix},
\quad P_{1+i}=\begin{pmatrix}0&E_i\\-E_i&0\end{pmatrix}, \; 1\leq i\leq m-1.
$$
Then
$$
P_iP_j+P_jP_i=2\delta_{ij}\,Id.
$$
$P_0,\cdots,P_{m}$ generate a representation of $C_{m+1}'$ on ${\mathbb R}^{2l}$.

The Cartan-M\"{u}nzner polynomials of the examples of Ferus, Karcher and M\"{u}nzner are
\begin{equation}\label{FKM}
F:=2|u|^4-2\sum_{i=0}^{m}(\langle P_i u, u\rangle)^2,\quad u\in{\mathbb R}^{2l}, \quad l=k\delta_{m}.
\end{equation}
Note that we recover Nomizu's example when $m=1$.

By a straightforward calculation, we conclude the following~\cite{FKM}.

\begin{proposition}
The two multiplicities of the associated isoparametric hypersurface are
$$
(m, k\delta_m-m-1),
$$
where $m,k\in{\mathbb N}$ to make the second entry positive, and the Clifford action operates on the focal manifold of codimension $1+m$ in the ambient sphere. Moreover, for $m\equiv 0\;(\text{mod}\; 4)$, there are $[k/2]+1$ incongruent isoparametric hypersurfaces associated with each multiplicity pair, to be indicated by $[k/2]$ underlines in the following grid table.
\end{proposition}
%We have the following grid table for these multiplicity pairs.

\vspace{3mm}

\begin{center}
\begin{tabular}{|c|c|c|c|c|c|c|c|c|c|c|}
\hline
\diagbox{$k$}{$\delta_m$}&$1$&$2$&$4$&$4$&$8$&$8$&$8$&$8$&$16$&$\cdots$\\
\hline
$1$&--&--&--&--&$(5,2)$&$(6,1)$&--&--&$(9,6)$&$\cdots$\\
\hline
\vspace{1mm}
$2$&--&$(2,1)$&$(3,4)$&$\underline{(4,3)}$&$(5,10)$&$(6,9)$&$(7,8)$&$\underline{(8,7)}$&$(9,22)$&$\cdots$\\
\hline
\vspace{1mm}
$3$&$(1,1)$&$(2,3)$&$(3,8)$&$\underline{(4,7)}$&$(5,18)$&$(6,17)$&$(7,16)$&$\underline{(8,15)}$&$(9,38)$&$\cdots$\\
\hline
\vspace{1mm}
$4$&$(1,2)$&$(2,5)$&$(3,12)$&$\underline{\underline{(4,11)}}$&$(5,26)$&$(6,25)$&$(7,24)$&$\underline{\underline{(8,23)}}$&$(9,54)$&$\cdots$\\
\hline
\vspace{1mm}
$5$&$(1,3)$&$(2,7)$&$(3,16)$&$\underline{\underline{(4,15)}}$&$(5,34)$&$(6,33)$&$(7,32)$&$\underline{\underline{(8,31)}}$&$(9,70)$&$\cdots$\\
\hline
$\vdots$&$\vdots$&$\vdots$&$\vdots$&$\vdots$&$\vdots$&$\vdots$&$\vdots$&$\vdots$&$\vdots$&$\vdots$\\
\hline
\end{tabular}
\end{center}
\vspace{3mm}

Among other things, Ferus, Karcher, and M\"{u}nzner established

\begin{theorem}

\begin{itemize}
\item[(1)] The multiplicity pairs of the homogeneous isoparametric hypersurfaces with four principal curvatures are precisely those listed in the first, second, and fourth columns of the table, together with $(9,6)$ and the two pairs $(2,2)$ and $(4,5)$ not listed on the table.
%The examples with multiplicity pairs
%on the first, second, fourth columns, and $(4,3)$ and $(9,6)$ are exactly the homogeneous examples, barring
%the two with multiplicities $(2,2)$ and $(4,5)$ not on the list.
\item[(2)] The isoparametric hypersurfaces with multiplicity pairs in the
third and seventh columns are exactly the inhomogeneous examples
constructed by Ozeki and Takeuchi.
\end{itemize}
\end{theorem}
So, except for the first, second and fourth columns, we have infinitely many families,
each with infinitely many members, of inhomogeneous isoparametric hypersurfaces
with four principal curvatures. Note that we also have the fact that such a  hypersurface with multiplicity pair $(1,l)$ or $(2,l)$ is congruent to the one
with multiplicity pair $(l,1)$ or $(l,2)$~\cite[6.5]{FKM}. Note also that Cartan classified the cases
when the multiplicities are $\{1,1\}$ and $\{2,2\}$, both being homogeneous~\cite{Ca4}.

Of particular interest is that the authors proved that the focal manifold $M_{-}$, of the inhomogeneous family with multiplicity pair $(m=m_{+},m_{-})=(3, 4k)$ constructed by Ozeki and Takeuchi, is homogeneous while $M_{+}$ is not. On the other hand, the two 
inhomogeneous examples of multiplicity pair $(m=m_{+},m_{-})=(8,7)$ constructed by the authors have the property that both focal manifolds of the indefinite isoparametric hypersurface, i.e., the one with $P_0\cdots P_8\neq \pm Id$, are inhomogeneous, whereas for the definite one, i.e., the one with $P_0\cdots P_8=\pm Id$, $M_{+}$ is homogeneous while $M_{-}$ is not, so that, in particular, the two hypersurfaces are not congruent; moreover, neither of them is congruent to the one with multiplicity pair $(m=m_{+},m_{-})=(7,8)$ constructed by Ozeki and Takeuchi. All of these properties were proved via geometric methods without resorting to the aforementioned classification in the homogeneous category. In fact, they showed through geometric means that most of the examples on the list are inhomogeneous.

Coming next to the arena is the thesis of Abresch~\cite{Ar}, in which he added a projective structure, in the case of $g=4,$ or 6, to the topological structure of M\"{u}nzner so that now Stiefel-Whitney classes come into play with Steenrod squares to obtain the following:

Assume $m_{-}\leq m_{+}$. If $g=4$, then
\begin{description}
\item[($4A$)] $m_{-}+m_{+}+1$ is divisible by $2^k:=\min\{2^\sigma:2^\sigma> m_{-},\; \sigma\in {\mathbb N}\}$, \;\text{or},
\item[($4B_1$)] $m_{-}$ is a power of $2$ and $2m_{-}$ divides $m_{+}+1$,\; \text{or},
\item[($4B_2$)] $m_{-}$ is a power of $2$ and $3m_{-}=2(m_{+}+1)$.
\end{description}
Moreover, if $g=6$, then
$$m_{-}=m_{+}=1,\; \text{or} \; 2.$$

Based on Abresch's work, Dorfmeister and Neher~\cite{DN} succeeded in classifying the case when $g=6$ and $m_{+}=m_{-}=1$. It is the homogeneous space.

Grove and Halperin~\cite{GH} established that when $g=4$, $m_{+}=m_{-}$ only when $(m_{+},m_{-})=(1,1),$ or $(2,2)$. As mentioned above, Cartan indicated without proof~\cite{Ca4} that the hypersurfaces are  homogeneous in both cases, for which an outline was laid out in~\cite{OT} for the $(2,2)$ case and a detailed proof following the outline was done by Tom Cecil (unpublished notes). The crucial step is to establish the existence of a point of Condition A mentioned below~\eqref{S^3}. In fact, in Cartan's formula~\eqref{F(Z)}, one can show that the element $e:=(e_{12}+e_{34})/\sqrt{2}$ assuming $F(e)=-1$ is a point of Condition A, where $e_{ij}\in so(5,{\mathbb R})$, $i<j$, is the matrix whose only nonzero entries are at $(i,j)$ and $(j,i)$ slots with value 1 and -1, respectively.

\subsection{1990-1999} Tang~\cite{Tang} pursued Abresch's setup further and obtained the refined result for $g=4$ that states that no isoparametric hypersurfaces of type $4B_1$ or $4B_2$ exist if $m_{-}\neq 1,2,4,$ or $8$, while Fang~\cite{Fa} followed up to assert that the multiplicity pairs $(2,2)$ and $(4,5)$ are the only possibilities in the case of $4B_2$. For type $4A$, Fang~\cite{Fa0} settled a large portion of the multiplicity problem:

\begin{theorem} Suppose $g=4$ and $m_{+}\leq m_{-}$. Then $m_{+}+m_{-}+1$ is divisible by $\delta(m_{+})$ if $m_{+}\equiv 5, 6, 7 (\text{mod}\; 8)$. In particular, the multiplicity pairs $(m_{+},m_{-}),m_{+}\leq m_{-},$ of 
isoparametric hypersurfaces with four principal curvatures are exactly those in the above grid table of the examples constructed by Ferus, Karcher, and M\"{u}nzner, provided $m_{+}\equiv 5, 6, 7 (\text{mod}\; 8)$.
\end{theorem}

The multiplicity problem was finally settled by the remarkable paper of Stolz~\cite{St}:

\begin{theorem}\label{stolz} Let $g=4$. The multiplicity pairs $(m_{+},m_{-}),m_{+}\leq m_{-},$ of 
isoparametric hypersurfaces with four principal curvatures are exactly those in the above grid table of the examples constructed by Ferus, Karcher, and M\"{u}nzner, besides the pairs $(2,2)$ and $(4,5)$ not in the grid table.
\end{theorem}

He established that if $(m_{+},m_{-}),m_{+}\leq m_{-},$ is neither $(2,2)$ nor $(4,5)$,
then $m_{+}+m_{-}+1$ is a multiple of $2^{\phi(m_{-}-1)}$, where $\phi(n)$ denotes the number
of numerals $s,1\leq s\leq n,$ such that $s\equiv 0,1,2,4\;({\rm mod}\; 8)$. In particular, one can see
easily that such pairs $(m_{+},m_{-})$ are exactly those in
the grid table of Ferus, Karcher, and M\"{u}nzner.

His approach is reminiscent of the theorem of Adams:

\begin{theorem} If there are $k$ independent vector fields on $S^n$, then $n+1$
is a multiple of $2^{\phi(k)}$.
\end{theorem}

The core technique Adams developed for proving the above theorem on vector fields
was to what Stolz reduced his proof.

Fang also showed~\cite{Fa2} that, when $g=6$, the isoparametric hypersurface is diffeomorphic (respectively, homotopic) to the homogeneous example when the equal multiplicity is $1$ (respectively, $2$); the statement in fact holds true in the more general proper Dupin category.

Gary Jensen, Tom Cecil, and I started thinking seriously about the classification of isoparametric hypersurfaces with four principal curvatures around 1998-1999. When I saw Stolz's result, I said to myself:``Ha-ha! What else could such a hypersurface be, except for the ones of Ferus, Karcher, and M\"{u}nzner?'' 

\section{The classification}
\subsection{2000-2009} We spent quite a lot of time working at understanding the underlying geometry of the Ferus-Karcher-M\"{u}nzner examples in the early 2000s. Let us look at~\eqref{FKM} more closely, where we set $m=m_{+}$ for convenience. $M_{+}$ is of codimension $1+m$ in $S^n$ defined by quadrics,
$$
M_{+}:=\{x\in S^n: \langle P_0(x),x\rangle=\cdots=\langle P_m(x),x\rangle=0\}.
$$
Its unit normal sphere at $x$ is 
$$
(UN)_x:=\{P(x): P=\sum_{i=0}^m a_i P_i, \;\;\sum_{i=0}^m a_i^2=1\},
$$
where $P$ constitute a round sphere ${\mathbb S}^m$ in the linear space of symmetric matrices of size $(n+1)\times (n+1)$ equipped with the inner product $\langle A,B\rangle:=-tr(AB)/2$. The map
$$
P\in{\mathbb S}^m\longmapsto P(x)\in UN_x,
$$
is an isometry for each $x\in M_{+}$.

Consider the unit normal bundle $UN$ of $M_{+}$ with the natural projection 
$$
\pi: UN\longrightarrow M_{+}.
$$ 
The Levi-Civita connection on $M_{+}$ naturally splits the tangent bundle of $UN$ into horizontal and vertical bundles ${\mathcal H}$ and ${\mathcal V}$, respectively,
$$
T(UN)={\mathcal V}\oplus {\mathcal H}.
$$
At each $n\in UN$ with base point $x$, the shape operator $S_n$ at $x$ admits three eigenspaces $E_0^n,E_1^n,E_{-1}^n$ with eigenvalues $0,1,-1$. Explicitly, for $n=P(x)$, 
$$
\langle S_n(v),w\rangle=-\langle P(v),w\rangle,\quad v,w\;\;\text{tangent at}\; x,
$$
so that, in particular, it is easily checked by a dimension count that
\begin{equation}\label{dis}
\aligned
&E_0^n=\text{span}\{PQ(x): Q\perp P\; \text{in}\; {\mathbb S}^m\},\\
&E_{1}^n=\{v: P(v)=-v,\;v\perp Q(x)\;\forall Q\in{\mathbb S}^m\},\\
&E_{-1}^n=\{v: P(v)=v,\;v\perp Q(x)\;\forall Q\in{\mathbb S}^m\},
\endaligned
\end{equation}
which are lifted to $UN$ to further split ${\mathcal H}$ into three subbundles ${\mathcal E}_0,{\mathcal E}_{1},{\mathcal E}_{-1}$, respectively, so that 
\begin{equation}\label{tun}
T(UN)={\mathcal V}\oplus {\mathcal E}_0\oplus{\mathcal E}_{1}\oplus{\mathcal E}_{-1}.
\end{equation}
Now, for each $P\in{\mathbb S}^m$, the set
$$
{\mathcal F}_P:=\{P(x):x\in M_{+}\}\subset UN,
$$
defines the section 
\begin{equation}\label{sec}
s_P:M_{+}\longrightarrow UN,\quad x:\longmapsto P(x).
\end{equation}
The tangent space to ${\mathcal F}_P$ at $P(x)$ is, by~\eqref{dis}, 
\begin{equation}\label{delta}
T({\mathcal F}_P)={\mathcal G}\oplus {\mathcal E}_{1}\oplus{\mathcal E}_{-1},
\end{equation}
where ${\mathcal G}$ is the graph of the orthogonal bundle map
$$
P:{\mathcal E}_0\longrightarrow {\mathcal V},\quad PQ(x)\longmapsto Q(x),\;\forall Q\perp P \in {\mathbb S}^m.
$$
In other words, as $P$ varies in ${\mathbb S}^m$, we can define the distribution 
\begin{equation}\label{del}
\Delta:={\mathcal G}\oplus {\mathcal E}_{1}\oplus{\mathcal E}_{-1},
\end{equation}
which is integrable (involutive) with leaves ${\mathcal F}_P$.

Conversely, let $g=4$. Suppose we are given an isoparametric hypersurface $M\subset S^n$ with the focal manifolds $M_{\pm}$. We let $m:=m_{+}$ and let $UN$ be the unit normal bundle of $M_{+}$. As above, we have the decomposition~\eqref{tun}. 
Suppose now there is an orthogonal bundle map
\begin{equation}\label{or}
{\mathcal O}:{\mathcal E}_0\longrightarrow {\mathcal V},
\end{equation}
which gives rise to a distribution $\Delta$ defined similarly as in~\eqref{del} by the graph ${\mathcal G}$ of ${\mathcal O}$. We wish to find a characterization of the integrability of $\Delta$. 

To this end, we let
$X_a,X_p,X_\alpha,X_\mu$ be an orthonormal frame spanning, respectively, ${\mathcal V},{\mathcal E}_0,{\mathcal E}_1,{\mathcal E}_{-1}$ over $T(UN)$, where $a,p,\alpha,\mu$ denote the indexes parametrizing the corresponding spaces with
\begin{equation}\label{range}
\aligned
&1\leq a\leq m=m_{+},\quad m+1\leq p\leq 2m,\\
&2m+1\leq \alpha\leq 2m+m_{-},\quad 2m +m_{-}+1\leq \mu\leq 2m+2m_{-}.
\endaligned
\end{equation}
This convention will be enforced henceforth. We then let $\theta^a,\theta^p,\theta^\alpha,\theta^\mu$ be the dual frame to the orthonormal frame. We set
$$
\omega^i_j:=\langle dX_j,X_i\rangle,
$$
and write
\begin{equation}\label{frame}
\omega^i_j=\sum_k F^i_{jk}\theta^k,
\end{equation}
where without the specification in~\eqref{range} indexes are understood to take values in all possible index ranges. By~\eqref{frame}, 
\begin{equation}\label{FS}
\aligned
&F^\alpha_{pa}=-S_{X_a}(X_\alpha,X_p),\quad F^\mu_{pa}=S_{X_a}(X_\mu,X_p),\\
&F^\mu_{\alpha a}=S_{X_a}(X_\alpha,X_\mu)/2,
\endaligned
\end{equation}
the respective components of the shape operator in the normal direction $X_a$, where by a slight abuse of notation we use $X_p,X_\alpha,X_\mu$ to also denote their pushforwards via $\pi:UN\rightarrow M_{+}$. 

We will see the geometric meaning of $F^\mu_{\alpha p}$ later. 

Isoparametricity of $M$ imposes many constraints on $F^i_{jk}$, which is not our concern here (see~\cite[p.16]{CCJ}). 

Now the orthogonal bundle map ${\mathcal O}$ in~\eqref{or} gives a choice of $X_a$ once $X_p$ are given, namely, we may specify 
\begin{equation}\label{signconvention}
X_{p-m}:=-{\mathcal O}(X_p).
\end{equation}
With this choice it follows that the distribution $\Delta$ in~\eqref{del} is the kernel of $\theta^a+\theta^{a+m}$; differentiating while invoking the structural equations for $d\theta^i$ (see~\cite[(5.1), p. 16]{CCJ}), we obtain

\vspace{2mm}

\begin{proposition}~\cite[p. 137]{Chi} $\Delta$ is integrable if and only if
\begin{equation}\label{3eq}
\aligned
F^\mu_{\alpha\, p}&=F^\mu_{\alpha\, p-m},\\
F^\alpha_{a+m\, b}&=-F^\alpha_{b+m\, a},\\
F^\mu_{a+m\, b}&=-F^\mu_{b+m\, a}.
\endaligned
\end{equation}
\end{proposition}

Let us understand the geometric meaning of this proposition. Each local integral leaf of $\Delta$ now gives rise to a local section
\begin{equation}\label{s}
s: M_{+} \longrightarrow UN, \quad x\longmapsto n(x),
\end{equation}
similar to the one in~\eqref{sec}. Now, a special feature of $g=4$ is that any normal $n(x)\in\pi^{-1}(x)$ also lives in $M_{+}$ (see the discussions in Section 4.2); for clarity of notation, we denote $n(x)$ by $x^\#\in M_{+}$. Furthermore, the normal space $N_{x^\#}$ to $M_{+}$ at $x^\#$ is
\begin{equation}\label{nord}
N_{x^\#}={\mathbb R}x\oplus E_0,
\end{equation}
where $E_0, E_1,E_{-1}$ are the eigenspaces of the shape operator $S_{n(x)}$ at $x$ with eigenvalues $0,1,$ and $-1$, respectively. We may now find the eigenspaces $E_0^{\#},E_1^{\#},E_{-1}^{\#}$ of the shape operator $S_x $ at $x^{\#}\in M_{+}$ to be
$$
\aligned
&E_0^{\#}=n(x)^\perp \subset \;N_x=\;\text{the normal space at}\; x,\\
&E_{\pm 1}^{\#}=E_{\pm 1}.
\endaligned
$$
It becomes evident now, in view of~\eqref{FS} and~\eqref{nord}, that 

\vspace{2mm}

{\em $F^\mu_{\alpha\, p}$ represents the component $S(X_\alpha,X_\mu)/2$ of the second fundamental form of $M_{+}$ at $x^{\#}$ in the normal direction $X_p\in N_{x^{\#}}$}, 

\vspace{2mm}

\noindent for which
we have the identity~\cite[(5.6), p.16]{CCJ}
\begin{equation}\label{FF}
\sum_{a=1}^m (F^\mu_{\alpha\, a} F^\nu_{\beta\, a}+F^\mu_{\beta\, a} F^\nu_{\alpha\, a})=\sum_{p=m+1}^{2m} (F^\mu_{\alpha\, p} F^\nu_{\beta\, p}+F^\mu_{\beta\, p} F^\nu_{\alpha\, p}).
\end{equation}

Note that the local map
$$
f_s: x\in M_{+}\longrightarrow M_{+},\quad x\longmapsto x^{\#},
$$
that the section $s$ in~\eqref{s} induces is a local isometry on $M_{+}$, where $(f_s)_*$ maps $v_{\pm}\in E_{\pm 1}$ to ${\mp} v_{\pm}\in E_{\pm}^{\#}$, and maps $v_0\in E_0$ to $w_0:={\mathcal O}(v_0)\in E_0^{\#}$.

\vspace{2mm}

\begin{proposition}~\cite[p. 138]{Chi} Assume $\Delta$ is integrable. The local isometries $f_s$ of $M_{+}$ extend to ambient isometries of $S^n$ for all $s$ if and only if 
\begin{equation}\label{4theq}
\omega^a_b-\omega^{a+m}_{b+m}=\sum_p L^p_{b\,a}(\theta^{p-m}+\theta^p)
\end{equation}
for some smooth functions $L^p_{b\,a}$.
\end{proposition}
The idea is to show that the local isometries $f_s$, preserving the first fundamental form, also preserve the second fundamental form and the normal connection form if and only if the four sets of identities hold. As a result, the analyticity of isoparametricity now implies that the local isometries are in fact global ones. 

We can now fix an $x_0\in M_{+}$ and consider the unit normal $m$-sphere $\pi^{-1}(x_0)$. Through each $n\in \pi^{-1}(x_0)$ there passes a unique section $s_n$ of $\Delta$ whose associated local isometry $f_{s_n}$ is now a global one induced by an isometry $P_n$ of $S^n$ that extends $f_{s_n}$. We thus have an $S^m$-worth of isometries $P$ of the ambient sphere inducing isometries of $M_{+}$. Analyzing how the $S^m$-worth of the ambient isometries $P$ interact with the expansion formula~\eqref{exp-f} of Ozeki and Takeuchi confirms the following.

\begin{proposition}\label{pivo}\cite[142-154]{Chi}
Assuming~\eqref{3eq} and~\eqref{4theq}, the $S^m$-worth of ambient isometries $P$ form a round sphere in the linear space of symmetric matrices of size $(n+1)\times(n+1)$ equipped with the standard inner product, so that the isoparametric hypersurface $M$ is one constructed by Ferus, Karcher, and M\"{u}nzner.
\end{proposition}

The proof in~\cite[33-51]{CCJ} of this proposition was through extremely long calculations by differential forms with remarkable cancellations! In contrast, the different proof in~\cite{Chi} is considerably shorter and more conceptual. On the other hand, the former is entirely local and the intuition behind it is also rather clear, in that since isoparametric hypersurfaces are defined by an overdetermined differential system, exterior-differentiating enough times should result in sufficiently many constraints for the conclusion.

Note that the three equations in~\eqref{3eq} are algebraic in $F^i_{jk}$ while the fourth one in~\eqref{4theq} is a system of PDEs in them. It would be desirable if~\eqref{4theq} could be suppressed. This is indeed possible if $M_{+}$ is ``sufficiently curved''. More precisely, we introduce the following definition.

\vspace{2mm}

\begin{definition}~\cite[p. 19]{CCJ} The unit normal bundle $UN$ of $M_{+}$ satisfies the spanning property at some $n$ over base point $x$, if there is an $X$ in $E_{1}$ such that
$$
S(X, \cdot): E_{-1} \longrightarrow n^{\perp}
$$ 
is surjective, and there is a $Y$ in $E_{-1}$ such that
$$
S(\cdot,Y):E_{1}\longrightarrow n^{\perp}
$$
is surjective, at $x$. Here, $E_0,E_1,E_{-1}$  are the eigenspaces of the shape operator $S_n$ at $x$ with eigenvalues $0,1,$ and $-1$, respectively, and $n^{\perp}$ is the orthogonal complement of $n$ in the normal space at $x$.
\end{definition}
Equivalently, the spanning property is equivalent to the local conditions that the Euclidean vector-valued bilinear form
\begin{equation}\label{bxy}
B(X,Y):=(\sum_{\alpha\mu} F^\mu_{\alpha\, 1}\,x_\alpha y_\mu,\sum_{\alpha\mu} F^\mu_{\alpha\, 2}\,x_\alpha y_\mu,\cdots,\sum_{\alpha\mu} F^\mu_{\alpha\, m}\,x_\alpha y_\mu)
\end{equation}
for $\alpha,\mu$ in the specified range in~\eqref{range} satisfies that $B(X,\cdot):{\mathbb R}^{m_{-}}\mapsto {\mathbb R}^m$ is surjective for some $X$, and $B(\cdot, Y):{\mathbb R}^{m_{-}}\rightarrow {\mathbb R}^m$ is surjective for some $Y$.

\begin{proposition}\label{tal} \cite[Proposition 19, p. 28]{CCJ} Suppose the unit normal bundle $UN$ of $M_{+}$ satisfies the spanning property at some $n$. Then around $n$ the first equation in~\eqref{3eq}, i.e., $F^{\mu}_{\alpha\, a}=F^\mu_{\alpha\, a+m},$ implies the remaining equations in~\eqref{3eq} and~\eqref{4theq}.
\end{proposition}
The proof utilizes the spanning property and various identities~\cite[pp. 16-17]{CCJ} of covariant derivatives of $F^i_{jk}$. 

In view of Propositions~\ref{pivo} and~\ref{tal}, it suffices to find conditions to warrant 
\begin{equation}\label{F^_}
F^\mu_{\alpha\, a}=F^\mu_{\alpha\, a+m},\quad \forall \alpha,\mu, \quad\text{with the spanning property},
\end{equation}
 for the isoparametric hypersurface with four principal curvatures to be one constructed by Ferus, Karcher, and M\"{u}nzner.

Note that Takagi's classification said at the end of Section 4 readily follows now. In this case, $m=1$, $a=1$ and $p=2=a+m$ in~\eqref{range}. From~\eqref{FF}, one reads off 
$$
F^\mu_{\alpha\, a}=\pm F^\mu_{\alpha\, a+m}
$$ 
for all $\alpha$ and $\mu$, where we can assume the sign is positive. 
Meanwhile, we need to verify that the spanning property holds. Suppose it is not true and $B(X,\cdot)$ is the zero map for every $X$ in~\eqref{bxy}. Then it must be that $B=0$ with $m=1$ gives that $F^\mu_{\alpha\, 1}=0$ for all $\mu,\alpha$. Now, with $a=1$,~\eqref{0a} extracts out of~\eqref{S^3} the identity
\begin{equation}\label{ab}
A_1A_1^{tr}+2B_1B_1^{tr}=Id,\quad A_1=\begin{pmatrix}F^\mu_{\alpha\, a}\end{pmatrix}=0,
\end{equation}
so that 
$$
B_1^{tr}=\begin{pmatrix}F^{\alpha=2m+1}_{p\, a},\cdots,F^{\alpha=2m+m_{-}}_{p\, a}\end{pmatrix} 
$$
satisfies $B_1B_1^{tr}=Id/2$. This is impossible, whence follows Takagi's classification.

It is at this point that algebraic geometry comes into play. We refer the reader to~\cite{Chi4} for a rather detailed account of the commutative algebra to be employed in the following. Since the algebro-geometric method in~\cite{CCJ} is superseded by the simpler and more effective method that prevails in~\cite{Ch1},~\cite{Ch2},~\cite{Chi3}, to be discussed later, I will only indicate briefly the inductive steps that are engaged by looking at the case when $m=2$ classified by Ozeki and Takeuchi in~\cite[II]{OT}.
Let us first look at the spanning property.

\begin{lemma}\label{irr}~\cite[II, p. 45]{OT} If $m_{-}\geq m+2$, then each of the $m$ polynomials 
$$
p_a(x,y):=\sum_{\alpha\,\mu}\, F^\mu_{\alpha\, a}\, x_\alpha y_\mu,\quad 1\leq a\leq m,
$$
is irreducible. 
\end{lemma}

\begin{proof} Note that 
$$
4p_a(x,y)=\langle \begin{pmatrix} x\\y\end{pmatrix},\begin{pmatrix}0&A_a\\A_a^{tr}&0\end{pmatrix}\begin{pmatrix}x\\y\end{pmatrix}\rangle,
$$
where $A_a$ is given in~\eqref{0a}. Hence, 
$$
\text{rank}(p_a)=2\;\text{rank}(A_a),
$$
where $\text{rank}(p_a)$ is that of the symmetric matrix
$$
U:=\begin{pmatrix}0&A_a\\A_a^{tr}&0\end{pmatrix}.
$$
Let 
$$
V:=S_a-U
$$
(for notational ease, we use $U$ to also denote its augmentation by zeros to match the size of $S_a$). We have $\text{rank}(S_a)\leq \text{rank}(U)+\text{rank}(V)$. Since $S_a$ is similar to $S_0$, we know $\text{rank}(S_a)=2m_{-}$, whereas it is clear that $\text{rank}(U)=2\;\text{rank}(A_a)$ and $\text{rank}(V)\leq 2m$. Putting these together, we derive
\begin{equation}\label{2m}
2(m_{-}-m)\leq 2\;\text{rank}(A_a)=\text{rank}(p_a).
\end{equation}
If $p_a$ is reducible, then $p_a=fg$ for two linear polynomials $f$ and $g$ of the form
$$
f=\sum_\alpha a_\alpha x_\alpha +\sum_\mu a_\mu y_\mu,\quad g=\sum_\alpha b_\alpha x_\alpha +\sum_\mu b_\mu y_\mu.
$$
Let 
$$
a:=\begin{pmatrix}a_\alpha&a_\mu\end{pmatrix}^{tr},\quad b:=\begin{pmatrix}b_\alpha&b_\mu\end{pmatrix}^{tr}
$$
of size $2m_{-}\times 1$. Then $4p_a=(ab^{tr}+ba^{tr})/2$ has rank $\leq 2$. Thus~\eqref{2m} implies $m_{-}\leq m+1$, a contradiction.
\end{proof}

\begin{lemma}\label{inde} If $m_{-}\geq m+1$, then the polynomials $p_a$, $1\leq a\leq m$, are linearly independent.
\end{lemma}

\begin{proof}  Suppose they are linearly dependent. There are nonzero constants, $c_1,\cdots,c_m,$ not all zero, such that
$$
0=c_1\, A_1+\cdots+c_m\,A_m =\sqrt{c_1^2+\cdots+c_m^2} A_n
$$
for some unit normal $n$. So, we may assume without loss of generality that $A_1=0$. But then~\eqref{ab} results in $B_1B_1^{tr}=Id/2$, which means that the row vectors of $B_1$ of size $m_{-}\times m$ are linearly independent, so that $m_{-}\leq m$. This is absurd.
\end{proof}

Note that both lemmas are equally good if we complexify the polynomials $p_a$, as they have real coefficients. 

In the case of $m=2$ in the classification of Ozeki and Takeuchi~\cite[II]{OT}, since the two polynomials $p_1$ and $p_2$ are irreducible and linearly independent, one cannot be a constant multiple of the other. We conclude that if $f\,p_1+g\,p_2=0$ for two polynomials $f$ and $g$ in $x$ and $y$, then $p_2$ divides $f$. We can put this succinctly in the language of regular sequences in commutative algebra.

\begin{definition}\label{regularsequence}~\cite[Definition 1, p. 84]{Chi4}\label{def} 
A regular sequence in the complex polynomial ring $P[l]$ in $l$ variables
is a sequence $p_0,\cdots,p_k$ in $P[l]$ such that firstly the variety
defined by $p_0=\cdots=p_k=0$ in ${\mathbb C}^l$ is not empty. Moreover,
$p_i$ is 
a non-zerodivisor in the quotient ring $P[l]/(p_0,\cdots,p_{i-1})$ for
$1\leq i\leq k$;
in other words, any relation
$$
p_0f_0+\cdots+p_{i-1}f_{i-1}+p_if_i=0
$$
will result in $f_i$ being in the form
$$
f_i=p_0h_0^i+\cdots+p_{i-1}h_{i-1}^i
$$
for some $h_0^i,\cdots,h_{i-1}^i\in P[l]$ for $1\leq i\leq k$.
\end{definition}
A regular sequence imposes strong algebraic independence amongst its elements~\cite[2.3, p. 91]{Chi4}.

\begin{proposition}{\rm (Special case)} \label{sition}~\cite[p. 62]{CCJ} Assume $m=2$ and $m_{-}\geq 4$. We use $p_1^{\mathbb C}$ and $p_2^{\mathbb C}$ to denote the complexification of $p_1$ and $p_2$. Let
$V_2^{\mathbb C}$ be the variety carved out by $p_1^{\mathbb C}=p_2^{\mathbb C}=0$. Then $\dim(V_2^{\mathbb C})=2m_{-}-2.$ In particular, the spanning property is true.
\end{proposition}

\begin{proof} $p_1^{\mathbb C}$ and $p_2^{\mathbb C}$ are irreducible by Lemma~\ref{irr}. It is well known that the irreducible $p_1^{\mathbb C}$ cuts out an irreducible variety $V_1$ and the irreducible $p_2^{\mathbb C}$ cuts out the variety $V_2$ of pure codimension 1 in $V_1$ since they are linearly independent~\cite[Theorem 5, p.58]{Sh}. In particular, the real counterpart $V_2$ of $V_2^{\mathbb C}$ satisfies
$$
\dim(V_2)\leq 2m_{-}-2.
$$
We claim the equality holds. To this end, consider the map
$$
V_2\stackrel{\iota}\longrightarrow{\mathbb R}^{m_{-}}\times{\mathbb R}^{m_{-}}\stackrel{\pi_1}\longrightarrow{\mathbb R}^{m_{-}},
$$
where $\iota$ is the standard embedding and $\pi_1$ is the projection onto the first summand. Also, for each $x\in {\mathbb R}^{m_{-}}$, consider the map
$$
S_2^x: y\longrightarrow (p_1(x,y),p_2(x,y)).
$$
Note that 
\begin{equation}\label{kernel}
\dim(\text{kernel}(S_2^x))\geq m_{-}-2>0,
\end{equation} 
and the kernel of $S_2^x$ is exactly the set 
$$
\{y\in{\mathbb R}^{m_{-}}: (x,y)\in (\pi_1\circ\iota)^{-1}(x)\}.
$$
Hence, the map $\pi_1\circ\iota$ is surjective. Now, the set 
$$
Z:=\{x\in{\mathbb R}^{m_{-}}: \text{kernel}(S_2^x)\;\text{assumes the minimum dimension}\; t\}
$$
is Zariski open. So, by Sard's theorem, there is an irreducible component $W$ of $V_2$ whose image via $\pi_1\circ\iota$ contains a regular value $x\in Z$, so that $(\pi_1\circ\iota)^{-1}(q)\simeq{\mathbb R}^t$ for $q$ in a neighborhood of $x$; in other words,
$(\pi_1\circ\iota)^{-1}(Z)$ around $x$ in $W$ is a product. We conclude, by~\eqref{kernel},
$$
\dim(W)=t+m_{-}\geq m_{-}+m_{-}-2=2m_{-}-2.
$$
Consequently, the claim 
\begin{equation}\label{claim}
\dim(W)=2m_{-}-2
\end{equation}
is proven. In particular, $t=m_{-}-2$.

A by-product is that around this regular value $x$, we enjoy the property 
$$
\dim((\pi_1\circ\iota)^{-1}(x))=t=m_{-}-2,
$$
 or, that $S_2^x$ is surjective, so that the spanning property is true for some $x$. Likewise, the spanning property is true for some $y$.
 \end{proof}

Let us next investigate the validity of~\eqref{F^_}. In view of~\eqref{0a} and the notation set around it, we denote the components of the second fundamental form by
\begin{equation}\label{p0}
\tilde{p}_i(v):=\langle S_i(v),v\rangle,\quad 0\leq i\leq m,
\end{equation}
and we define
$$
{\mathcal D}=\{z\in E_{+}\oplus E_{-}:  |z|=1,\; \tilde{p}_i(z)=0,\; i=0,\cdots,m\}.
$$

\begin{lemma}\cite[Proposition 25, p. 51]{CCJ} ${\mathcal D}=(E_{+}\oplus E_{-})\cap M_{+}.$
\end{lemma}

\begin{proof} This follows when we set $t=w_0=\cdots=w_m=0$ in~\eqref{exp-f}.
\end{proof}

Notation as around~\eqref{0a} and the index range convention~\eqref{range} prevailing, similar to the discussions below~\eqref{s}, let us denote $n_0$ by $x^\#\in M_{+}$. Three lines above~\eqref{FF} we gave the geometric meaning of $F^\mu_{\alpha \, p}$ at $x^\#$ vs. $F^\mu_{\alpha\, a}$ at $x$. An immediate consequence of the preceding lemma is the following crucial observation.

\begin{corollary}\label{on}\cite[Proposition 28, p. 53]{CCJ} We have that the zero locus $V$ of
$$
p_a:=\sum_{\alpha\,\mu} F^\mu_{\alpha\, a}x_\alpha y_\mu,\quad 1\leq a\leq m,
$$ 
and the zero locus of 
$$
\overline{p}_a:=\sum_{\alpha\,\mu} F^\mu_{\alpha\, p}\, x_\alpha y_\mu,\quad m\leq p\leq 2m,
$$
are identical in ${\mathbb R}P^{m-1}\times {\mathbb R}P^{m-1}$.
\end{corollary}

\begin{proof} The eigenspaces of the shape operator $S_x$ at $x^\#$ with eigenvalue $\pm 1$ are identical, respectively, with the eigenspaces of $S_{n_0}$ at $x$ with eigenvalues $\pm 1$. The preceding corollary then implies that 
$$
{\mathcal D}=\{z=(x,y)\in E_{+}\oplus E_{-}: |z|=1,\, \sum_{\alpha\,\mu} F^\mu_{\alpha\, i}x_\alpha y_\mu=0, 0\leq i\leq m\}.
$$
Furthermore, observe that in~\eqref{p0}
$$
\tilde{p}_0=|x|^2-|y|^2=0
$$
on ${\mathcal D}$ while $|x|^2+|y|^2=1$ so that $|x|=|y|=1/2$. The result follows by projectivization.
\end{proof}
Note that the same conclusion need not necessarily hold on ${\mathbb C}P^{m-1}\times {\mathbb C}P^{m-1}$ when we complexify the associated polynomials, since the complexification of the real irreducible components need not exhaust all the complex irreducible components. Nevertheless, the situation is clear in the complex irreducible case, as follows.

\begin{corollary}\label{tion}\cite[p. 63]{CCJ}\label{prime} Notation as in the preceding corollary, if the zero locus $V^{\mathbb C}$ of $p_a^{\mathcal C},1\leq a\leq m$, is an irreducible variety in ${\mathbb C}P^{m-1}\times{\mathbb C}P^{m-1}$ and, moreover, the complex dimension of $V^{\mathbb C}$ equals the real dimension of $V$, then $F^\mu_{\alpha\, a}=F^\mu_{\alpha\, a+m}$ in~\eqref{F^_} after an orthogonal basis change.
\end{corollary}

\begin{proof} Since the zero locus $V^{\mathbb C}$ of $p_a^{\mathbb C},1\leq a\leq m,$ is irreducible, its ideal $I:=(p_1^{\mathbb C},\cdots,p_m^{\mathbb C})$ is prime. On the other hand, the preceding corollary says that $\overline{p}_1,\cdots,\overline{p}_m$ vanish on $V$ having the same real dimension as the complex dimension of $V^{\mathbb C}$, we conclude that $\overline{p}_1^{\mathbb C},\cdots,\overline{p}_m^{\mathbb C}$ also vanish on $V^{\mathbb C}$ so that by Hilbert's Nullstellensatz, 
$$
\overline{p}_a^{\mathbb C}=\sum_{b=1}^m\, r_{ab}\, p_b^{\mathbb C},\quad 1\leq a\leq m,
$$
for some polynomials $r_{ab}$. We deduce that
$$
\overline{p}_a=\sum_{b=1}^m\, c_{ab}\, p_b,\quad 1\leq a\leq m,
$$
for some constants $c_{ab}$, since $p_a$ and $\overline{p}_a$ are all of bidegree $(1,1)$. In other words,
$$
F^\mu_{\alpha\, a+m}=\sum_{b=1}^m \,c_{ab}\,F^\mu_{\alpha\, b},\quad 1\leq a\leq m,
$$
whence~\eqref{FF} gives that $\begin{pmatrix} c_{ab}\end{pmatrix}$ is indeed an orthogonal matrix. The conclusion follows by a suitable orthogonal basis change.
\end{proof}

Let us now return to the classification of Ozeki and Takeuchi in the case $m=2$. By Lemma~\ref{irr}, Lemma~\ref{inde},~\eqref{F^_},~\eqref{claim}, Proposition~\ref{sition}, Corollary~\ref{on}, and Corollary~\ref{tion}, we can conclude that the isoparametric hypersurface of multiplicity pair  $(m=2,m_{-})$, $m_{-}\geq 4,$ is one constructed by Ferus, Karcher, and M\"{u}nzner, and hence is homogeneous, so long as we can verify that the zero locus of $p_1^{\mathbb C}$ and $p_2^{\mathbb C}$ is irreducible. This is indeed the case. It relies on the criterion of Serre for irreducible varieties.

\begin{theorem}\label{Serre}~\cite[Theorem 1, p. 85; Theorem 2, p. 89]{Chi4} Let $I_k:=(f_1,\cdots,f_k)$ be the ideal generated by a regular sequence
$$
f_1,\cdots,f_k, \quad k\leq l, 
$$
in the complex polynomial ring $P[l]$ in $l$ variables $z_1,\cdots,z_l$, whose zero locus is $V_k$. Let $J_k$ be the subvariety of $V_k$ consisting of all points of $V_k$ where the Jacobian matrix
\begin{equation}\label{jacobian}
\partial(f_1,\cdots,f_k)/\partial(z_1,\cdots,z_l)
\end{equation} 
is not of full rank $k$. Suppose the codimension of $J_k$ is $\geq 1$ in $V_k$. Then $V_k$ is reduced, i.e., $I_k$ is a radical ideal.

Moreover, if $V_k$ is connected and the codimension of $J_k$ is $\geq 2$ in $V_k$, then $V_k$ is irreducible, i.e., $I_k$ is a prime ideal. In particular, the connectedness condition is automatically satisfied when $f_1,\cdots,f_k$ are homogeneous polynomials.
\end{theorem}

Notation as above, we know the homogeneous $p_1^{\mathbb C}$ and $p_2^{\mathbb C}$ form a regular sequence by Lemma~\ref{irr} and the discussion above Definition~\ref{regularsequence}. It thus suffices to check that the codimension 2 estimate holds true in Theorem~\ref{Serre}, where $k=2$ and $l=2m_{-}$, for the classification of Ozeki and Takeuchi to go through. An elementary linear algebra argument establishes the following.

\begin{lemma}\label{le12.1}~\cite[Lemma 49, p. 64]{CCJ} We set $m:=m_{+}$ as usual.
Notation is as in~\eqref{0a}. There is an
orthonormal basis in $E_{1}$ and an orthonormal basis in $E_{-1}$ such that
relative to these bases we have, in~\eqref{0a},
\begin{description}
\item[(1)] $B_{1}=C_{1}=\begin{pmatrix}0&0\\0&\sigma\end{pmatrix}$ with $\sigma={\rm diag}(\sigma_1,\cdots,\sigma_r)$, $\sigma_s>0,  \forall s$,
\item[(2)] $A_{1}=\begin{pmatrix}I&0\\0&\Delta\end{pmatrix}$, where $\Delta$
is an $r\times r$ matrix in block form
$$
\Delta = \begin{pmatrix}\Delta_{1}&0&0&0&\dots\\0&\Delta_{2}&0&0&\dots
\\0&0&\Delta_{3}&0&\dots\\\vdots&\vdots&\vdots&\vdots&\vdots\end{pmatrix}
$$
with $\Delta_{1}=0$ and $\Delta_{i},i\geq 2,$ nonzero skew-symmetric
matrices in block form
$$
\Delta_i = \begin{pmatrix}0&f_{i}&0&0&\dots\\-f_{i}&0&0&0&\dots\\0&0&0&f_{i}&\dots\\
0&0&-f_{i}&0&\dots\\\vdots&\vdots&\vdots&\vdots&\vdots\end{pmatrix}, \;\text{and},
$$
\item[(3)] $\Delta_i^2=-(1-2\sigma_i^2)\, Id.$
\end{description}
\end{lemma}

Clearly, the lemma holds for any fixed $a\geq 1$ in the expression of $S_a$.

\begin{corollary}\label{co12.1}
$\dim(Ker(A_{a}))=\dim(\Delta_{1})\leq r={\rm rank}(B_{a})\leq m$ for all $a\geq 1$.
\end{corollary}

With $m=2$, we first estimate the dimension of the
subvariety $Z_{2}$ of  ${\mathbf C}^{m_{-}}\times {\mathbf C}^{m_{-}}$ at
each point
of which the Jacobian matrix~\eqref{jacobian} of $p_{1}^{\mathbf C},p_{2}^{\mathbf C}$
is of rank $<2$.
At $(x,y)\in Z_{2}$, the differentials
$dp_{1}^{\mathbf C},dp_{2}^{\mathbf C}$ are linearly
dependent, i.e., there are $c_{1},c_{2}\in {\mathbf C}$, depending on
$(x,y)$, such that
$$
0=\sum_{a=1}^2 \,c_{a}dp_{a}^{\mathbf C}=\sum_{\alpha} (\sum_{a,\mu}
c_{a}F^{\mu}_{\alpha a}y_{\mu})dx_{\alpha}
+\sum_{\mu}(\sum_{a,\alpha}
c_{a}F^{\mu}_{\alpha a}x_{\alpha})dy_{\mu},
$$
which requires that the coefficients of $dx_{\alpha}$ be zero and the
coefficients of $dy_{\mu}$ be zero.
Thus
$$
Z_{2}=\{(x,y)\in
{\mathbf C}^{m_{-}}\times {\mathbf C}^{m_{-}}
:\exists (c_{1},c_{2}),
\sum_{a}c_{a}A_{a}^{tr}x=\sum_{a}c_{a}A_{a}y=0\}.
$$
Accordingly, for a fixed $(c_{1},c_{2})$ let us define
$$
Z_{(c_{1},c_{2})}:=\{(x,y)\in
{\mathbb C}^{m_{-}}\times {\mathbb C}^{m_{-}}
:\sum c_{a}A_{a}^{tr}x=\sum_{a}c_{a}A_{a}y=0\}.
$$
Consider the incidence space $Y_{2}$ in
${\mathbf C}P^1\times {\mathbf C}^{m_{-}}\times {\mathbf C}^{m_{-}}$
given by
$$
\{([c_{1}:c_{2}],x,y): (x,y)\in Z_{(c_{1},c_{2})}\}.
$$
The standard projection of $Y_{2}$ to
${\mathbf C}^{m_{-}}\times {\mathbf C}^{m_{-}}$ maps $Y_{2}$ onto $Z_{2}$.
Let 
$$
\pi:Y_2\longrightarrow {\mathbb C}P^1
$$
be the standard projection of $Y_{2}$ to ${\mathbb C}P^1$.
Then with respect to $\pi$ we have
\begin{equation}\label{eq12.145}
\dim(Z_{2})\leq \dim(Y_{2})\leq \dim({\rm base})+\dim({\rm fiber}),
\end{equation}
where $\dim({\rm fiber})$ is the
maximal dimension of all fibers.
It is easier to estimate the dimension
of the fibers $\pi^{-1}\{[c_1:c_2]\} = Z_{(c_{1},c_{2})}$. In fact, it comes down to estimating the dimension of
$$
T_{(c_{1},c_{2})}:=\{y\in {\mathbf C}^{m_{-}}:\sum_{a} c_{a}A_{a}y=0\}
$$
for a fixed $(c_{1},c_{2})$, because this estimate will also be a valid upper bound for the dimension of %$U_{(c_1,\dots,c_2)} := 
$\{x \in {\mathbb C}^{m_{-}} : \sum_a c_a A_a^{tr} x =0 \}$, thus giving us the estimate
$$
\dim(Z_{(c_{1},c_{2})})\leq 2\dim(T_{(c_{1},c_{2})}).
$$

%To estimate the dimension of $Z_{(c_{1},\dots,c_{n})}$.

\vspace{2mm}

\noindent {\bf Case (1)}. $c_{1},c_{2}$ are either all real or all purely
imaginary. This is the easier case. It is essentially Corollary~\ref{co12.1}. We have
$$
\dim(T_{(c_{1},c_{n})})\leq r\leq m
$$
and
\begin{equation*}
\dim (Z_{(c_1,c_2)}) \leq 2 \dim T_{(c_1,c_2)} \leq 2m
\end{equation*}

\noindent {\bf Case (2)}. $c_{1},c_{2}$ are not all real and not all purely imaginary.
Write
$$
c_{k}=\alpha_{k}+\sqrt{-1}\beta_{k}.
$$
%where not all $\alpha_{k}$ and not all $\beta_{k}$ are zero.
We may assume without loss of generality that
$$
c_{1}S_{e_{1}}+c_{2}S_{e_{2}}=(\alpha_{1}
+\sqrt{-1}\beta_{1})S_{e_{1}}+\sqrt{-1}\beta_{2}S_{e_{2}}.
$$
By restricting to the $A$-block in $S$ again we see that 
%$(\sum_{a}c_{a}A_{a})y=0$ is reduced to
$$
\beta_{2}A_{2}y
=\sqrt{-1}(\alpha_{1}+\sqrt{-1}\beta_{1})A_{1}y;
$$
we may assume both coefficients are nonzero, or else
we would be back to Case (1). Hence we are now handling
\begin{equation}\label{eq12.15}
(A_{2}-zA_{1})y=0
\end{equation}
for some nonzero $z\in {\mathbf C}$. By Lemma~\ref{le12.1}, we may assume
$$
A_{1}=\begin{pmatrix}I&0\\0&\Delta\end{pmatrix}.
$$
Write
$$
A_{2}=\begin{pmatrix}\Theta&\Lambda\\\Omega&\Gamma\end{pmatrix}
$$
of the same block sizes as $A_{1}$. Inspecting the identity
$$
A_1A_2^{tr}+A_2A_1^{tr}+B_1B_2^{tr}+B_2B_1^{tr}=0,
$$
extracted out of~\eqref{S^3}
%By the second equation of~\eqref{eq:5:6},
%which is
%$$
%A_{2}{\T A_{1}}+A_{1}{\T A_{2}}+2(B_{2}{\T B_{1}}
%+B_{1}{\T B_{2}})=0,
%$$
we obtain
\begin{equation}\label{eq12.16}
\Theta+\Theta^{tr}=0.
\end{equation}
With this and $r\leq m=2$, one can eventually come up with the estimate~\cite[p. 71]{CCJ}
%since by~\eqref{eq12.2} and~\eqref{eq12.3}, there is an $r \times r$ diagonal matrix $\sigma$ such that
%$$
%B_{1}={\T B_{1}}=\begin{pmatrix}0&0\\0&\sigma\end{pmatrix}.
%$$
$$
\dim(T_{(c_{1},c_{2})})\leq (m_{-}+r)/2\leq (m_{-}+m-1)/2,
$$
whose details are not our concern here. Note that the upper bound $(m_{2}+m_{1}-1)/2$, instead of the weaker
$(m_{2}+m_{1})/2$,
holds because we know $m_{2}+m_{1}$ is an odd number if $2\leq m<m_{-}$
by M\"{u}nzner~\cite[II]{Mu} (or, by the more general Abresch~\cite{Ar} mentioned in Subsection 5.1), and as a result
\begin{equation}\label{extra}
\dim(Z_{(c_{1},c_2)})\leq 2\dim(T_{(c_{1},c_2)})\leq
m_{-}+m-1=m_{-}+1.
\end{equation}
On the other hand, the base of $\pi$ consists of finitely many points since generic $[c_1:c_2]$ would make $Z_(c_1,c_2)$ null. Putting these together,~\eqref{eq12.145} yields
\begin{equation}\label{est}
\dim(Z_2)\leq m_{-}+1.
\end{equation}
To achieve an estimate for $\dim(J_2)$, where $J_2$ and $V_2$ are defined in~\eqref{jacobian} with $m=2$, by a well known fact~\cite[Corollary 5, p. 57]{Sh}, the zero locus $V_2$ of $p_1^{\mathbb C}$ and $p_2^{\mathbb C}$ has the property
$$
\dim(V_2)\geq 2m_{-}-2,
$$
so that by the fact that $J_2\subset Z_2$ we achieve the {\em a priori} estimate
$$
\dim(J_2)\leq \dim(V_2)-2
$$
if $m_{-}\geq 5$, in which case the ideal $(p_1^{\mathbb C},p_2^{\mathbb C})$ is prime as a result of Serre's criterion. The isoparametric hypersurface with the multiplicity pair $(2,m_{-})$, $m_{-}\geq 5$, is thus homogeneous. Although this suffices for the conclusion in the case of $m=2$ and $m_{-}\geq 4$ as $m_{-}$ are all odd~\cite[II, p. 49]{OT}, we can perform a general cutting procedure to reach the conclusion for $m_{-}\geq 4$. Indeed, consider the surjective map
$$
f:{\mathbb C}^{m_{-}}\times {\mathbb C}^{m_{-}}\longrightarrow {\mathbb C}^2, \quad (x,y) \longmapsto (p_1^{\mathbb C}(x,y),p_2^{\mathbb C}(x,y)).
$$
$Z_2$ is the set where $df$ is of rank $<2$ and $J_2=f^{-1}(0)\cap Z_2$. 

Let $W_{k}$, $0\leq k\leq 1,$ be the subvarieties of $Z_2$ where $df$ is of rank $\leq 1-k$. We have $W_0\supset W_1$. Let $X_j:=W_j\setminus W_{j+1}$. Then $Z_2$ is stratified into $X_0,X_1$ with $X_j$ Zariski open in $W_j$, where $df$ is of rank $1-j$ on $X_j$. Around each $(x,y)$ in $X_0$, we may assume $\nabla(p_2^{\mathbb C})$ is a multiple of $\nabla(p_1^{\mathbb C})$ considered as column vectors. Since 
$$
\begin{pmatrix}p_1^{\mathbb C}&p_2^{\mathbb C}\end{pmatrix}=\begin{pmatrix}x&y\end{pmatrix}\begin{pmatrix}\nabla(p_1^{\mathbb C})&\nabla(p_2^{\mathbb C})\end{pmatrix},
$$
we see that near generic $(x,y)$ in $X_0$, the image of $f$ is of dimension 1. Therefore, a generic line cut through the origin performed in the target space of $f$ cuts down the dimension by 1 from $Z_2$ in the dimension estimate of $J_2=f^{-1}(0)$ so that 
$\dim(J_2)\leq \dim(Z_2)-1\leq m_{-}$ by~\eqref{est} and thus $\dim(J_2)\leq \dim(V_2)-2$ when $m_{-}\geq 4$, whence follows the classification of Ozeki and Takeuchi when $m=2$ and $m_{-}\geq 4$..

In general, an inductive procedure~\cite[pp. 60-61]{CCJ} takes care of all $m$:

\begin{proposition}\label{red}\cite[Proposition 46, p. 61]{CCJ} Notation as in Theorem~\ref{Serre} above, set $m=m_{+}$ as usual and assume $m_{-}\geq m+2$. For $n\leq m$, let 
$$
Z_n=\{(x,y)\in
{\mathbf C}^{m_{-}}\times {\mathbf C}^{m_{-}}
:\exists (c_{1},\cdots,c_n),
\sum_{a=1}^n\, c_{a}A_{a}^{tr}x=\sum_{a=1}^n\, c_{a}A_{a}y=0\},
$$
and 
$$
f_n:(x,y)\in{\mathbf C}^{m_{-}}\times {\mathbf C}^{m_{-}}\longrightarrow  {\mathbb R}^n, \quad (x,y)\longmapsto (p_1(x,y),\cdots,p_n(x,y))
$$
with $J_n=f^{-1}(0)\cap Z_n$. If $m_{-}\geq 2m$, then 
\begin{equation}\label{codim2}
\dim(J_n)\leq \dim(V_n)-2
\end{equation}
for all $n\leq m$. If $m=2m-1$, then $dim(J_n)\leq \dim(V_n)-2$ for all $n\leq m-1$ while $dim(J_m)\leq \dim(V_m)-1$.
\end{proposition}
As a consequence of the preceding proposition, an argument similar to the one outlined above for the case $m=2$ results in the classification theorem:

\begin{theorem}\label{clathm}\cite[Theorem 47, p. 61]{CCJ} If $m_{-}\geq 2m_{+}-1$, then an isoparametric hypersurface with four principal curvatures is one constructed by Ferus, Karcher, and M\"{u}nzner, where the Clifford action operates on $M_{+}$.
\end{theorem}

Note that when $m_{-}=2m-1$ in Proposition~\ref{red}, the ideal $I_m=(p_1^{\mathbb C},\cdots,p_m^{\mathbb C})$ is only radical by Theorem~\ref{Serre}, so that {\em a priori} the arguments in Corollary~\ref{tion} do not appear to work. However, as in~\eqref{claim}, the real variety $V_{m}$ cut out by $p_1,\cdots,p_m$ is of real dimension $2m_{-}-m$. Let $W$ be an irreducible component of $V_{m}^{\mathbb C}$ containing an irreducible component of $V_m$ of real dimension $2m_{-}-m$; the complex dimension of $W$ is $2m_{-}-m$. Then Nullstellensatz holds locally, around generic point $z:=(x,y)$ of $W$, for the ideal $I_m$ so that, in the notation of Corollary~\ref{tion},
\begin{equation}\label{ms}
%F^\mu_{\alpha\, a+m}\,y_\mu=\sum_{b=1}^m \,c_{ab}\,F^\mu_{\alpha\, b}\, y_\mu,\quad 
F^\mu_{\alpha\, a+m}\,x_\alpha y_\mu=\sum_{b=1}^m \,c_{ab}\,F^\mu_{\alpha\, b}\, x_\alpha y_\mu,  \quad\text{summed on}\;\alpha,\mu,
\end{equation}
for $1\leq a\leq m$ and some rational functions $c_{ab}(x,y)$ smooth around $z$. 
Since the projection sending $W\subset {\mathbb R}^{m_{-}}\times{\mathbb R}^{m_{-}}$ to $x\in{\mathbb R}^{m_{-}}$ is surjective around $z$ by the local structure of $W$ given in Proposition~\ref{sition}, we see
\begin{equation}\label{same}
F^\mu_{\alpha\, a+m}\,x_\alpha y_\mu=\sum_{b=1}^m \,c_{ab}(x)\,F^\mu_{\alpha\, b}\, x_\alpha y_\mu
\end{equation}
by Taylor-expanding along ${\mathbb R}^{m_{-}}$ around $x$, where $c_{ab}(x)$ is the collection in the analytic $c_{ab}(x,y)$ of the terms depending only on $x$, and so partial differentiating with respect to $y_\mu$ gives
$$
F^\mu_{\alpha\, a+m}\,x_\alpha =\sum_{b=1}^m \,c_{ab}(x)\,F^\mu_{\alpha\, b}\, x_\alpha,     \quad\text{summed on}\; \alpha.
$$
Taking second order partial derivatives with respect to $x$, denoted by $c_{ab}''$, we obtain
$$
c_{ab}'' F^{\mu}_{\alpha\, b} x_\alpha = 0
$$
for all $\mu$, which implies $c_{ab}''=0$ by the surjectivity of the map 
\begin{equation}\label{sm}
S_m^x: y\longmapsto (p_1(x,y),\cdots,p_m(x,y)),\quad p_a=\sum_{\alpha\mu}\, F^\mu_{\alpha a}\, x_\alpha y_\mu,
\end{equation}
as detailed in Proposition~\ref{sition}. Hence, $c_{ab}(x)$ are of degree at most 1 in $x$ and hence must be a constant by comparing types in~\eqref{same}. We arrive at the same conclusion as in the case when $I_m$ is a prime ideal.
In particular, this takes care of the case $(m,m_{-})=(2,3)$ in Ozeki and Takeuchi's classification since their explicit formula~\cite[II, p. 49]{OT} for $p_1$ and $p_2$ results in that $V_2^{\mathbb C}$ is reduced when $m_{-}=3$~\cite[Remark 53, p. 73]{CCJ}.

As good as it gets, Theorem~\ref{clathm} exactly reaches the borderline to account for all multiplicity pairs except for the four exceptional cases $(m=m_{+},m_{-})=(3,4), (4,5),(6,9),$ and $(7,8)$, by a look at Stolz's multiplicity result Theorem~\ref{stolz}.

I spent the two years 2008-2009 on and off thinking about an effective method to break the barrier to forge beyond. The first thing that came in mind was that we had not utilized the second fundamental form of $M_{+}$ to full generality, where the $E_0$-component had been entirely ignored when, inspired by Ozeki and Takeuchi, studying the restricted form $p_a=\sum_{\alpha\mu}\, F^\mu_{\alpha\, a}x_\alpha y_\mu$ defined only on $E_1\oplus E_{-1}$. Moreover, a glance at~\eqref{0a} showed that every component $\tilde{p}_a, 0\leq a\leq m,$ of the second fundamental form of $M_{+}$ is irreducible since they are equivalent to
$$
\tilde{p}_0=|x|^2-|y|^2
$$ 
after a coordinate change, so long as $m_{-}\geq 2$, obtained almost for free. 
%Furthermore, the fact that the codimension 2 estimate~\eqref{codim2} holds for all $n\leq m-1$, so that $I_k$ are all prime ideals for $k\leq m-1,$  guarantees that $p_1^{\mathbb C},\cdots,p_{m}^{\mathbb C}$ form a regular sequence~\cite[Proposition 39]{CCJ}, %which in turn warrants that $p_1^{\mathbb C},\cdots,p_m^{\mathbb C}$ cut ``transversally'' (even on the ideal level) so that the cuts by $p_1^{\mathbb C}, \cdots, p_n^{\mathbb C}$, $n\leq m,$ all produce the crucial dimension $2m_{-}-n$ of the resulting zero %locus. 
It seemed reasonable to replace $p_1,\cdots, p_m$ by $\tilde{p}_0,\cdots,\tilde{p}_m$ and find criteria to warrant that $\tilde{p}_0,\cdots,\tilde{p}_m$ form a regular sequence by working through the successive codimension 2 estimates in~\eqref{codim2}. The idea paid off.

\subsection{2010-2019} From now on we stick to the convention that $m:=m_{+}\leq m_{-}$. 
For ease of notation, we denote the components of the second fundamental form of $M_{+}$ by $p_0,\cdots,p_m$ without the tilde, and, moreover, we drop the superscript ${\mathbb C}$ in $p^{\mathbb C}$ of a real polynomial $p$ whenever we indicate that $p$ lives in $P[l]$, the polynomial ring in $l$ complex variables, where $l:=2m_{-} +m$ is the dimension of $M_{+}$.

\vspace{2mm}

\noindent {\bf The case of multiplicity pair (3,4)}

\vspace{2mm}

In the expansion formula~\eqref{exp-f}, the components of the second and third fundamental forms of $M_{+}$ 
are intertwined in ten convoluted equations. The first three say that the
shape operator $S_n$ satisfies $(S_n)^3=S_n$ for any normal direction $n$,
which is agreeable with the fact that the eigenvalues of $S_n$ are $0,1,-1$ with
fixed multiplicities. Following the notation of Ozeki and Takeuchi~\cite[I]{OT}, set
$$
<p_a,q_b>:=\langle\nabla p_a,\nabla q_b\rangle,\quad 0\leq a, b\leq m.
$$

The fourth and fifth combined and the sixth are
\begin{equation}\label{paqb}
<p_a,q_b>+<p_b,q_a>=0,
\end{equation}
\begin{equation}
<<p_a,p_b>,q_c>+<<p_c,p_a>,q_b>+<<p_b,p_c>,q_a>=0,
\end{equation}
for distinct $a, b,c$. The seventh is
\begin{equation}\label{syzygy}
p_0q_0+\cdots+p_{m}q_{m}=0.
\end{equation}

\vspace{1mm}

Set $G:=\sum_{a=0}^{m}(p_a)^2.$ The last three are

\begin{equation}\label{last3}
16\sum_{a=0}^{m} (q_a)^2=16G\,|y|^2-<G,G>,
\end{equation}
\begin{equation}\label{last2}
\aligned
&8<q_a,q_a>=8(<p_a,p_a>|y|^2-(p_a)^2)+<<p_a,p_a>,G>\\
&-24G-2\sum_{b=0}^{m}<p_a,p_b>^2,
\endaligned
\end{equation}
\begin{equation}\label{last1}
\aligned
&8<q_a,q_b>=8(<p_a,p_b>|y|^2-p_a\,p_b)+<<p_a,p_b>,G>\\
&-2\sum_{c=0}^{m}<p_a,p_c><p_b,p_c>,\quad a,b\;\text{distinct}.
\endaligned
\end{equation}
~\eqref{syzygy} caught my eye while others seemed dauntingly entangling; it is the well known syzygy equation in commutative algebra. A property the syzygy equation enjoys is that when $p_0,\cdots,p_m$ form a regular sequence in $P[l]$, $l=2m_{-}+m$, we have 
\begin{equation}\label{qp}
q_a=\sum_{b} r_{ab}\,p_b,\quad r_{ab}=-r_{ba}
\end{equation}
for some first degree polynomials $r_{ab}$~\cite[Proposition 4]{Chi4}. This is exactly Condition B of Ozeki and Takeuchi. With this observation, I worked out in~\cite{Ch1} the {\em a priori} codimension 2 estimate in $P[l]$ for the components $p_0,\cdots,p_m$ of the second fundamental form of $M_{+}$. Indeed, following the convention in~\eqref{range}, let us write
\begin{equation}\label{eeq5}
\aligned
r_{ab}&=\sum_{\alpha} T_{ab}^{\alpha}\,u_\alpha+\sum_{\mu}T_{ab}^{\mu}\,v_\mu
+\sum_{p}T_{ab}^p\,z_p,\\
p_0&=\sum_{\alpha}(u_\alpha)^2-\sum_{\mu}(v_\mu)^2,\\
p_a&=4\sum_{\alpha\mu}F^\mu_{\alpha\, a}\,u_\alpha v_\mu
-2\sum_{\alpha p}F^\alpha_{p\, a}\,u_\alpha z_p
+2\sum_{\mu p}F^\mu_{p\,a}\,v_{\mu}z_p.\\
\endaligned
\end{equation}
%where 
%$$
%S^a_{\alpha\mu}:=\langle S(X_\alpha,Y_\mu),n_a\rangle,
%$$
%etc., with $X_\alpha,Y_\mu,$ and $Z_p$ the orthonormal bases for the coordinates $u_\alpha,v_\mu,$
%and $z_p$ of $E_1, E_{-1},$ and $E_0$, respectively. 
Note that $y$ in~\eqref{last3} through~\eqref{last1} denotes a tangential vector whose components are $u_\alpha, v_\mu,$ and $z_p$. Substituting~\eqref{eeq5} into~\eqref{qp} and comparing polynomial types yields that the coefficient of $(u_\alpha)^3$, denoted $q_a^{\alpha\alpha\alpha}$, in the polynomial expression of $q_a$ is
$$
q_a^{\alpha\alpha\alpha}=T^\alpha_{a\,0},
$$
while the right hand side of~\eqref{last3} ensures that there are no $(u^\alpha)^6$-terms. We conclude
\begin{equation}\label{t0a}
T_{a\,0}^\alpha=0=-T_{0\,a}^\alpha;\quad \text{likewise,}\; T^\mu_{a\,0}=-T^\mu_{0\,a}=0.
\end{equation}
In other words, 

\vspace{1mm}

{\em $r_{0a}$ consists of only $z_p$-terms in the expansion of $q_0$.}

\vspace{1mm}

\noindent Meanwhile, it is known that 
%by~\eqref{eq4},~\eqref{eq5},~\eqref{eq6}, and~\eqref{eq7}, keeping in mind
$q_0$ is homogeneous of degree 1 in $u_\alpha,v_\mu$ and $z_p$
by~\cite[I, p. 537]{OT}. We expand the right hand side of~\eqref{qp} to ascertain that the coefficient of the
$u_\alpha v_\mu z_p$-term
of $q_0$, denoted $q_0^{\alpha\mu p}$, is
\begin{equation}\label{eeq8}
q_0^{\alpha\mu p}=4\sum_{b\geq 1}T^p_{0\,b}\,F^\mu_{\alpha\, b}.
\end{equation}              

Here comes something particularly nice. It turns out that $q_0$ at $x$ encodes information of the second fundamental form at $x^{\#}:=n_0\in M_{+}$. More precisely, it is shown in~\cite{Ch1} that 
$$
q_0^{\alpha\mu p}=4 F^\mu_{\alpha\, p}.
$$
(See also~\cite[p. 97]{Chi4} for a different proof using the expansion formula~\eqref{exp-f}.)
%where the right hand side is given three lines above~\eqref{FF}. 
It follows that
$$
F^\mu_{\alpha\, p} =\sum_b\, f^p_b\,F^\mu_{\alpha\, b},\quad f^p_b:=T^p_{0\,b}.
$$
 We need only verify that $\begin{pmatrix}f^p_b\end{pmatrix}$ is an orthogonal matrix for the first equation in~\eqref{3eq} to hold. It is remarkable that this is indeed true as a consequence of a piece of commutative algebra~\cite[p. 153]{Ku} (see also~\cite[p. 91]{Chi4}):

\begin{proposition}\label{homog} Let $p_1,\cdots,p_k$ be a regular sequence in $P[l]$.
Let $F(t_1,\cdots,t_k)$ be a homogeneous polynomial of degree $d$
in $k$ variables with coefficients in $P[l]$.
Suppose $F(p_1,\cdots,p_k)=0$. Then 
all the coefficients of $F$ belong to $I=(p_1,\cdots,p_k)$.
\end{proposition}
The key idea of showing the orthogonality of $(f^p_b)$ is to rewrite~\eqref{last2} as a polynomial homogeneous
in all $p_ap_b$ whose coefficients are homogeneous polynomials of degree 2,
so that these coefficients are linear combinations of all $p_a$ by the preceding proposition.
Specifically, the coefficient of $(p_0)^2$ is
$$
16\sum_{a=1}^{m_1}(r_{0a})^2-16(\sum_\alpha(u_\alpha)^2+\sum_\mu(v_\mu)^2
+\sum_p(z_p)^2)
+4<p_0,p_0>,
$$
which is a linear combination of $p_0,p_1,\cdots,p_{m}$. Knowing that $r_{0a}$ are
functions of $z_p$ alone by~\eqref{t0a}, we invoke~\eqref{eeq5} and compare variable
types to conclude 
\begin{equation}\label{orthogonal}
\sum_{a=1}^{m}(r_{0a})^2=\sum_{p=m+1}^{2m}(z_p)^2,
\end{equation}
which thus asserts the orthogonality of the matrix
$\begin{pmatrix}f^p_b\end{pmatrix}$.
We may now assume 
\begin{equation}\label{am}
T^{a+m}_{0\, b}=f^{a+m}_b=\delta^a_b,\quad \text{so that}\; F^\mu_{\alpha\, a+m}=F^\mu_{\alpha\, a},
\end{equation}
and as a result of~\eqref{t0a} and~\eqref{am} derive
$$
r_{0b}=\sum_{a}\delta^{a}_{b}z_{a+m}=z_{b+m}.
$$
With the Einstein summation convention, we calculate
\begin{equation}\nonumber
\aligned
&q_0=r_{0b}\,p_b\\
&=2(\delta^{a}_{b}z_{a+m})
(S^b_{\alpha\mu}\,x_\alpha y_\mu+S^b_{\alpha\; c+m}\,x_\alpha z_{c+m}
+S^b_{\mu\; c+m}\,y_{\mu}z_{c+m}).
\endaligned
\end{equation}
Hence, we obtain
$$
\sum_{a b c\alpha}(\delta^{a}_{b}z_{a+m})\,(S^b_{\alpha\; c+m}u_\alpha z_{c+m})=0,
$$
or equivalently, 
$$
\sum_{ac}S^a_{\alpha\, c+m}\,z_{c+m}z_{a+m}=0.
$$
In other words, we have
\begin{equation}\label{anti-sym2}
F^\alpha_{c+m\, a}=-F^\alpha_{ a+m\, c}.
\end{equation}
Similarly, 
\begin{equation}\label{anti-sym3}
F^\mu_{c+m\, a}=-F^\mu_{a+m\, c}.
\end{equation}
In particular,~\eqref{3eq} are valid now. What is satisfying is that $m=m_{+}< m_{-}$ is the only condition to warrant that~\eqref{4theq} is true, so that Proposition~\ref{pivo} gives that the isoparametric hypersurface is one constructed by Ferus, Karcher, and M\"{u}nzner, as follows.

\begin{proposition}\label{isopa}~\cite[Proposition 4]{Ch1}
 Assume  $m=m_{+}<m_{-}$. Suppose Condition B of Ozeki and Takeuchi holds, which is the case when
the components of the second fundamental form $p_0,p_1,\cdots,p_m$
of $M_{+}$ form a regular sequence in $P[l]$, $l=2m_{-}+m$.
Then the isoparametric hypersurface is of the type constructed by Ferus, Karcher, and M\"{u}nzner with the Clifford action operating on $M_{+}$.
\end{proposition}
The proof uses another easier spanning property~\cite[Proposition 7, p. 18]{CCJ}, which states that if $m<m_{-}$ then $(F^\mu_{\alpha\, 1},\cdots,F^\mu_{\alpha m})$, $\forall \alpha, \mu$, span ${\mathbb R}^m$, to verify that~\eqref{3eq} implies~\eqref{4theq} if $m<m_{-}$, by utilizing certain identities in~\cite[Proposition 19]{CCJ}.

The preceding proposition is pivotal for the classification of the remaining cases, to be discussed later.

In view of the preceding proposition, we need only find conditions to guarantee that $p_0,\cdots,p_m$ form a regular sequence in $P[l]$, which is where Serre's criterion Theorem~\ref{Serre} of codimension 2 estimate comes in again. We record the inductive scheme to generate a regular sequence.

\begin{lemma}\label{generatereg}~\cite[Proposition 39]{CCJ} Let $p_0,\cdots,p_m\in P[l]$ be linearly independent
homogeneous polynomials of equal degree $\geq 1$.
For each $0\leq k\leq m-1$, let $V_k$ be the variety defined by $p_0=\cdots=p_k=0$,
and let $J_k$ be the subvariety of\, $V_k$, where the Jacobian
$$
\partial(p_0,\cdots,p_k)/\partial(z_1,\cdots,z_l)
$$
is not of full rank $k+1$. If the codimension of $J_k$ in $V_k$ is $\geq 2$
for all $0\leq k\leq m-1$, then $p_0,\cdots,p_m$ form a regular sequence.
\end{lemma}

%By this proposition, the classification of isoparametric hypersurfaces
%with four principal curvatures now boils down to exploring
%Lemma~\ref{generatereg} to warrant that the components
%$p_0,\cdots,p_{m_1}$ of the second fundamental form of $M_{+}$ constitute a regular sequence.
%To this end, let us look at the $p_0,\cdots,p_k,k\leq m_{1}-1$. Following Lemma~\ref{generatereg}
%we must estimate the codimension of $J_k$ in $V_k$ by understanding the rank of
%the Jacobian matrix of $p_0,\cdots,p_k$.

Let us parametrize ${\mathbb C}^{2m_{-}+m}$ by
points $(u,v,w)$, where $u,v\in{\mathbb C}^{m_{-}}$ and $w\in{\mathbb C}^m$.
%with coordinates $u_\alpha,v_\mu$, and $w_p$ parametrizing, respectively, $E_1,E_{-1},E_0$ as before. 
For $k\leq m$, let %where $1\leq\alpha,\mu\leq m_2$, and $1\leq p\leq m_1$. For $0\leq k\leq m_1$, let
$$
V_k:=\{(u,v,w)\in{\mathbb C}^{2m_{-}+m}:p_0(u,v,w)=\cdots=p_k(u,v,w)=0\}
$$
be the variety carved out by $p_0,\cdots,p_k$ in $P[l]$. We first estimate the dimension of the subvariety $X_k$ of ${\mathbb C}^{2m_{-}+m}$ defined by
$$
X_k:=\{(u,v,w)\in{\mathbb C}^{2m_{-}+m}
:\text{rank of Jacobian of}\; p_0,\cdots,p_k<k+1\}.
$$
This amounts to saying %that $dp_0,\cdots,dp_k$ are linearly dependent, or, 
that there are
constants $c_0,\cdots,c_k$ such that
\begin{equation}\label{linear}
c_0dp_0+\cdots+c_kdp_k=0.
\end{equation}
Since $p_a=\langle S_a(x),x\rangle$, we see $dp_a=2\langle S_a(x),dx\rangle$
for $x=(u,v,w)^{tr}$;
therefore, by~\eqref{linear},
$$
X_k=\{(u,v,w):(c_0S_0+\cdots+c_kS_k)\cdot (u,v,w)^{tr}=0\}.
$$
for $[c_0:\cdots:c_k]\in{\mathbb C}P^k$,
where $\langle S_{a}(X),Y\rangle=\langle S(X,Y),n_a\rangle$ is the shape operator of the focal
manifold
$M_{+}$ in the normal direction $n_a$.
Since $J_k=X_k\cap V_k$, by Lemma~\ref{generatereg} we wish to establish 
\begin{equation}\nonumber
\dim(X_k\cap V_k)\leq\dim(V_k)-2
\end{equation}
for $k\leq m-1$ to verify that $p_0,p_1,\cdots,p_{m}$ form a
regular sequence. 
%since\begin{equation}\label{prim}
%J_k=X_k\cap V_k.
%\end{equation}

Note that for a fixed
$\lambda=[c_0:\cdots:c_k]\in {\mathbb C}P^k$,
if we set
$$
{\mathscr S}_\lambda:=\{(u,v,w): (c_0S_0+\cdots+c_kS_k)\cdot (u,v,w)^{tr}=0\},
$$
then we have
\begin{equation}\label{unio}
X_k=\cup_{\lambda\in{\mathbb CP}^k} {\mathscr S}_\lambda.
\end{equation}
Thus, it boils down to estimating the dimension of ${\mathscr S}_\lambda$.

We break it into two cases.
If $c_0,\cdots,c_{k}$ are either all real or all purely imaginary, then
$$
\dim({\mathscr S}_\lambda)=m,
$$
since $c_0S_{n_0}+\cdots+c_kS_{n_k}=cS_n$ for some unit normal vector $n$ and
some nonzero real or purely imaginary constant $c$, and we know that the null space
of $S_n$ is of dimension $m$ for all normal $n$.

On the other hand, if $c_0,\cdots,c_k$ are not all
real and not all purely imaginary, then after a normal basis change,
we may assume that
\begin{equation}\label{slamb}
{\mathscr S}_\lambda=\{(u,v,w):(S_{1^*}-\iota_\lambda\, S_{0^*})
\cdot(u,v,w)^{tr}=0\}
\end{equation}
for some complex number $\iota_\lambda$ relative to a new orthonormal normal basis $n^*_0,n^*_1,\cdots,n^*_{k}$
in the linear span of $n_0,n_1,\cdots,n_k$ by the Gram-Schmidt process.
% explicitly, $n_0^*$ and $n_1^*$
%are obtained by decomposing $n:=c_0n_0+\cdots+c_kn_k$ into its real and imaginary
%parts $n=\alpha+\sqrt{-1}\beta$ and define $n_0^*$ and $n_1^*$ by performing the Gram-Schmidt
%process.
In matrix terms, the equation in~\eqref{slamb} assumes the form
\begin{equation}\label{Matrix}
\begin{pmatrix}0&A&B\\A^{tr}&0&C\\B^{tr}&C^{tr}&0\end{pmatrix}\begin{pmatrix}x\\y\\z\end{pmatrix}
=\iota_\lambda\begin{pmatrix}I&0&0\\0&-I&0\\0&0&0\end{pmatrix}\begin{pmatrix}x\\y\\z\end{pmatrix},
\end{equation}
where $x,y$, and $z$ are (complex) eigenvectors of $S_{0^*}$ with
eigenvalues $1,-1$, and $0$, respectively.

The decomposition Lemma~\ref{le12.1}
ensures that we can normalize the matrix on the left hand side
of~\eqref{Matrix} to
decompose $x,y,z$ into $x=(x_1,x_2),y=(y_1,y_2),z=(z_1,z_2)$
with $x_2,y_2,z_2\in{\mathbb C}^{r_\lambda}$, where $r_\lambda$ is the rank of $B$, or
intrinsically, 

\vspace{1mm}

{\em $m-r_\lambda$ is the dimension of the intersection of the kernels of $S_{0^*}$ and $S_{1^*}$.}

\vspace{1mm}

 With respect to this decomposition 
either $x_1=y_1=0$, or both are nonzero
with 
\begin{equation}\label{iota}
\iota_\lambda=\pm\sqrt{-1}.
\end{equation} 
In both cases we have
$x_2=-y_2$ and can be solved in $z_2$ 
so that $z$ may be chosen to be a free variable. Hence, either $x_1=y_1=0$, in which case
$$
\dim({\mathscr S}_\lambda)=m,
$$
or both $x_1$ and $y_1$ are nonzero,
in which case $y_1=\pm\sqrt{-1}x_1$ and so
\begin{equation}\label{eSt}
\dim({\mathscr S}_\lambda)=m+m_{-}-r_\lambda,
\end{equation}
where $x_1$ contributes the dimension count $m_{-}-r_\lambda$ while $z$ does $m$. Now, 
%since by~\eqref{prim} and
by~\eqref{unio} we see
\begin{equation}\label{jk}
J_k=X_k\cap V_k=\cup_{\lambda\in{\mathbb C}P^k} ({\mathscr S}_\lambda\cap V_k),
\end{equation}
where $V_k$ is also defined by %$p_0=\cdots=p_k=0$ and also by 
$p_{0^*}=\cdots=p_{k^*}=0$.
Let us cut ${\mathscr S}_\lambda$ by
$$
0=p_{0^*}=\sum_{\alpha}(x_\alpha)^2-\sum_{\mu}(y_\mu)^2
$$
to achieve an {\em a priori} estimate of $\dim(J_k)$.

\vspace{2mm}

\noindent Case 1: $x_1$ and $y_1$ are both nonzero. This is the case of
nongeneric $\lambda\in{\mathbb C}P^k$. We substitute $y_1=\pm \sqrt{-1}x_1$ and $x_2$
and $y_2$ in terms of $z_2$ into $p_{0^*}=0$ to
deduce that
$$
0=p_{0^*}=(x_1)^2+\cdots+(x_{m_{-}-r_\lambda})^2+z\; \text{terms};
$$
hence, $p_{0^*}=0$ cuts ${\mathscr S}_\lambda$ to reduce the dimension by 1.
In other words, now by~\eqref{eSt},
\begin{equation}\label{Sub}
\dim(V_{k}\cap{\mathscr S}_\lambda)\leq (m+m_{-}-r_\lambda)-1\leq m+m_{-}-1.
\end{equation}
%noting that $V_k$ is also cut out by $p_{0^*},p_{1^*},\cdots,p_{k^*}$. 
Meanwhile,
only a subvariety of $\lambda$ of
dimension $k-1$ in ${\mathbb C}P^k$ assumes $\iota_\lambda=\pm\sqrt{-1}$;
in fact, this subvariety
is the smooth hyperquadric 
\begin{equation}\label{quadric}
{\mathcal Q}_{k-1}:=\{\lambda=[c_0:\cdots:c_{k}]: c_0^2+\cdots+c_k^2=0\}
\end{equation}
in ${\mathbb C}P^k$. This is because
if we write $(c_0,\cdots,c_k)=\alpha+\sqrt{-1}\beta$ where
$\alpha$ and $\beta$ are
real vectors, then $\iota_\lambda=\pm\sqrt{-1}$
is equivalent to the conditions that
$\langle\alpha,\beta\rangle=0$ and $|\alpha|^2=|\beta|^2$. In other words,
the nongeneric
$\lambda\in{\mathbb C}P^k$ constitute the smooth hyperquadric.
Therefore, by~\eqref{jk}, an irreducible component ${\mathcal W}$ of 
$J_k$ over nongeneric $\lambda$
will satisfy

\begin{equation}\label{gen}
\dim({\mathcal W})\leq\dim(V_{k}\cap{\mathscr S}_\lambda)+k-1\leq
m+m_{-}+k-2.
\end{equation}
(Total dimension $\leq$ base dimension $+$ fiber dimension.)

\vspace{2mm}
 
\noindent Case 2: $x_1=y_1=0$. This is the case of generic $\lambda$, where
$\dim({\mathscr S}_\lambda)=m$, so that an irreducible component ${\mathcal V}$ of 
$J_k$ over generic $\lambda$ will satisfy
$$
\dim({\mathcal V})\leq m+k\leq m+m_{-}+k-2,
$$
as we may assume $m_{-}\geq 2$, noting that
the case $m=m_{-}=1$ follows from Takagi's classification for $m=1$ as mentioned above.

Putting these two cases together, we conclude that
\begin{equation}\label{Crucial}
\dim(J_k)=\dim(X_k\cap V_k)\leq m+m_{-}+k-2.
\end{equation}
On the other hand, since $V_k$ is cut out by $k+1$ equations
$p_0=\cdots=p_k=0$, we have
\begin{equation}\label{Lower-bound}
\dim(V_k)\geq m+2m_{-}-k-1.
\end{equation}
Therefore,
\begin{equation}\label{cod2}
\dim(J_k)\leq \dim(V_k)-2
\end{equation}
when $k\leq m-1$, taking $m_{-}\geq 2m-1$ into account.

In summary, with the assumption $m_{-}\geq 2m-1,$ we have established~\eqref{cod2} for $k\leq m-1$,
so that the ideal $(p_0,p_1,\cdots,p_k)$ is prime when $k\leq m-1$ by Serre's criterion Theorem~\ref{Serre}.
Lemma~\ref{generatereg} 
then implies that $p_0,p_1,\cdots,p_{m}$ form a regular sequence.
It follows by Proposition~\ref{isopa} 
that the isoparametric hypersurface is of the type constructed by Ferus, Karcher, and M\"{u}nzner. 

This approach, done  in~\cite{Ch1}, gives a considerably simpler proof of Theorem~\ref{clathm} above.
%\begin{theorem}\label{normalcool} Assume $m_{-}\geq 2m-1$. Then the isoparametric hypersurface
%with four principal curvatures is of the type constructed by Ferus, Karcher, and M\"{u}nzner.
%\end{theorem}

The extra bonus to this approach is that in~\cite{Ch1} I could also classify the exceptional case when the multiplicity pair is $(3,4)$, in which case it was known that there had been two existing examples, one is the inhomogeneous one of Ozeki and Takeuchi, where the Clifford action operates on $M_{+}$, and the other is the homogeneous one, where the Clifford action operates on $M_{-}$. I also knew that for the homogeneous example $p_0,p_1,p_2,p_3$ of $M_{+}$ did not form a regular sequence anymore; otherwise, Proposition~\ref{isopa} would give that the hypersurface was the one of Ozeki and Takeuchi.

It turned out that Condition A returned in an unexpected way to settle the case of multiplicity pair $(3,4)$. 
Note that the quantity $r_\lambda$ is entirely discarded in~\eqref{Sub}.
Condition A enables us to come up with a finer estimate on the right hand side of~\eqref{Sub} by utilizing $r_\lambda$.
%in which the quantity $r_\lambda$ is entirely discarded. This is where Condition A comes in.

Indeed, if we stratify the hyperquadric ${\mathcal Q}_{k-1}$, $k\leq m-1,$ given in~\eqref{quadric}, of nongeneric
$\lambda$ for which the dimension estimate may create complications,
into subvarieties 

\vspace{2mm}

{\em ${\mathcal L}_j$, of some dimension $d_j\leq k-1$, over which $r_\lambda=j$,} 

\vspace{2mm}

\noindent then by~\eqref{Sub} an irreducible component $\Theta_j$ of 
$V_k\cap (\cup_{\lambda\in {\mathcal L}_j} {\mathscr S}_\lambda)$
will satisfy                            
\begin{equation}\label{W}
\dim(\Theta_j)\leq\dim(V_{k}\cap{\mathscr S}_\lambda)+d_j\leq
m+m_{-}+d_j-1-j.
\end{equation}
%%\noindent Case 2. $x_1=y_1=0$. This is the case of generic $\lambda$, where
%%$\dim({\mathscr S}_\lambda)=m_1$, so that an irreducible component ${\mathcal V}$
%%of $W_k\cap (\cup_{\lambda\in{\mathcal G}}{\mathscr S}_\lambda)$,
%%where ${\mathcal G}$ is the Zariski open set of ${\mathbb C}P^k$ of generic
%%$\lambda$, will satisfy
%%\begin{equation}\label{v}
%%\dim({\mathcal V})\leq m_1+k.
%%\end{equation}
%%On the other hand, since $W_k$ is cut out by $k+1$ equations, we have
%%\begin{equation}\label{lower-bound}
%%\dim(W_k)\geq m_1+2m_2-k-1.
%%\end{equation}
%%\begin{lemma}\label{lm1} When $(m_1,m_2)=(4,5)$ $($respectively, $(m_1,m_2)=(6,9)$$)$
%%and $j\geq 2$, there holds in
%%equation~\eqref{W} the estimate
%%\begin{equation}\label{wk}
%%\dim({\mathcal W}_j)\leq \dim(W_k)-2
%%\end{equation}
%%for $k\leq m_1-1=3$ $($respectively, $k\leq 5$$)$.
%%\end{lemma}
%%\begin{proof} For~\eqref{wk} to be true, we must have both 
%%\begin{eqnarray}\nonumber
%%\aligned
%%m_1+m_2+k-2-j&\leq m_1+2m_2-k-3,\\
%%m_1+k&\leq m_1+2m_2-k-3
%%\endaligned
%%\end{eqnarray}
%%by~\eqref{W},~\eqref{v} and~\eqref{lower-bound}. The second inequality is $2m_2\geq 2k+3$, which is always true, while the
We run through the same arguments as those following~\eqref{Sub} to deduce that
the codimension 2 estimate~\eqref{cod2} holds true over ${\mathcal L}_j$ when
\begin{equation}\label{refinedest}
m_{-}\geq 2k+1-j-c_j,\quad \rm{where}
\end{equation}
$$c_j:=k-1-d_j=\rm{codimension\;of}\;{\mathcal L}_j\; \rm{in}\;{\mathcal Q}_{k-1}.$$ 

Note that the inequality $\dim({\mathcal V})\leq m+k$ below~\eqref{gen} for generic $\lambda$ in ${\mathbb C}P^{k-1}$ automatically results in  the codimension 2 estimate
$\dim({\mathcal V})\leq\dim(V_k)-2$, since $m_{-}>m$ in the remaining four exceptional cases. Thus, it suffices to consider only those $\lambda\in {\mathcal Q}_{k-1}$
for $k\leq m-1$ from now on.

Let us look at the case when $(m,m_{-})=(3,4)$, where now $0\leq k\leq m-1=2$. First, observe that~\eqref{refinedest} is automatically
satisfied when $j\geq 1$; the same is also true for all $j$ when $k=1$. So, we assume $k=2$ and $j=0$ henceforth.

With $k=2,j=0$, let $\lambda_0\in {\mathcal L}_0$ be generically chosen; we have $r_{\lambda_0}=j=0$.
Suppose that $M_{+}$ is free of points of Condition A everywhere. 
Let us span $\lambda_0$ by the orthonormal $n_0^*$ and $n_1^*$ completed to an orthonormal basis $n_0^*,n_1^*,n_2^*,n_3^*$. Since $r_{\lambda_0}=0$, the matrices $B=C=0$ and $A=I$ in~\eqref{Matrix} for $S_{1^*}$. For notational clarity,
let us denote the associated $B$ and $C$ blocks of the shape operator matrices $S_{a^*}$ by $B_{a^*}$ and $C_{a^*}$ for the normal basis elements
$n_1^*,
\cdots,n_{3}^*$.
It follows that $p_{0^*}=0$ and $p_{1^*}=0$ cut ${\mathscr S}_{\lambda_0}$ in the variety
\begin{equation}\nonumber
\{(x,\pm \sqrt{-1}x,z):\sum_\alpha (x_\alpha)^2=0\}.
\end{equation}
We may assume $(B_{2^*},C_{2^*})$ is nonzero as $M_{+}$ has no points of Condition A. 
Since $z$ is a free variable, $p_{2^*}=0$ will have nontrivial $xz$-terms
\begin{eqnarray}\nonumber
\aligned
0=p_{2^*}&=\sum_{\alpha p}S_{\alpha p}x_\alpha z_p+\sum_{\mu p}T_{\mu p}y_\mu z_p+x_\alpha y_\mu\;\text{terms}\\
&=\sum_{\alpha p}(S_{\alpha p}\pm\sqrt{-1}T_{\alpha p})x_\alpha z_p+x_\alpha y_\mu\;\text{terms},
\endaligned
\end{eqnarray}
taking $y=\pm\sqrt{-1}x$ into account, where $S_{\alpha p}:=\langle S(X^*_\alpha,Z^*_p),n^*_2\rangle$
and $T_{\mu p}:=\langle S(Y^*_\mu,Z^*_p),n^*_2\rangle$ are (real) entries of $B_{2^*}$ and $C_{2^*}$,
respectively, and $X^*_{\alpha}$,
$Y^*_{\mu}$, and $Z^*_{p}$ are orthonormal
eigenvectors for the eigenspaces of $S_{0^*}$ with eigenvalues $1,-1,$ and $0$,
respectively;
hence, the dimension
of ${\mathscr S}_{\lambda_0}$ will be cut down by 2 by $p_{0^*},p_{1^*},p_{2^*}=0$. In conclusion, modifying~\eqref{Sub} we have
\begin{equation}\nonumber
\dim(V_{k}\cap{\mathscr S}_\lambda)\leq m+m_{-}-2,
\end{equation} 
for all $\lambda\in{\mathcal L}_0$. As a consequence,
the right hand side of~\eqref{refinedest}, which is no bigger than 5 for $j=0$, is now no bigger than 4
with the additional cut $p_{2^*}=0$
so that the codimension 2 estimate goes through for ${\mathcal L}_0$ as well. It follows by Proposition~\ref{isopa} that the isoparametric hypersurface is in fact the one 
of Ozeki and Takeuchi, which thus has points of Condition A. This is a piece of absurdity to the assumption that $M_{+}$ has no points of 
Condition A. Therefore, $M_{+}$ does admit points of Condition A.
The result of Dorfmeister and Neher ~\cite{DN1} (see also~\cite{Chichi}) implies that the isoparametric hypersurface is then necessarily of the type of Ferus, Karcher, and M\"{u}nzner. In particular, it is either the inhomogeneous one by Ozeki and Takeuchi, or the homogeneous one.

The classification in the case of multiplicity pair $(3,4)$ is now achieved, as was done in~\cite{Ch1}.

%\begin{theorem} Let $(m_1,m_2)=(3,4)$. Then the isoparametric hypersurface is either the homogeneous one, or is the inhomogeneous one constructed by Ozeki and Takeuchi.
%\end{theorem}

\vspace{2mm}

\noindent {\bf The case of multiplicity pairs $(4,5)$ and $(6,9)$}

\vspace{2mm}

The role that Condition A plays in the classification in the case of multiplicity pair $(3,4)$ is that, assuming nonexistence of points of Condition A on $M_{+}$ makes the codimension 2 estimate go through. On the other hand, nonexistence of points of Condition A on $M_{+}$ in the case of multiplicity pairs $(4,5)$ and $(6,9)$ always holds true by the discussions following~\eqref{irre}. Guided by the $(3,4)$ case, I suspected that nonexistence of points of Condition A on $M_{+}$ could lead to something fruitful in the case of multiplicity pairs $(4,5)$ and $(6,9)$. It surely did.

The homogeneous example of multiplicity pair $(4,5)$ is a principal orbit of the action $U(5)$ on $so(5,{\mathbb C})$ given by 
$$
g\cdot Z=\overline{g}Zg,
$$ 
while the homogeneous example of multiplicity pair $(6,9)$, being the focal manifold $M_{-}$ of an isoparametric hypersurface of the type of Ferus, Karcher, and M\"{u}nzner on which the Clifford action operates, can be realized as the Clifford-Stiefel manifold
\begin{eqnarray}\nonumber
\aligned
M_{-}=\{&(\zeta,\eta)\in S^{31}\subset{\mathbb R}^{16}\times{\mathbb R}^{16}:\\
&|\zeta|=|\eta|=1/\sqrt{2},\zeta\perp\eta,\check{J}_i(\zeta)\perp\eta,i=1,\cdots,8\},
\endaligned
\end{eqnarray}
where $\check{J}_1,\cdots, \check{J}_8$ are the unique (up to equivalence) irreducible representation of the
(anti-symmetric)
Clifford algebra $C_8$ on ${\mathbb R}^{16}$ constructed by the octonion algebras as follows. 
%The representation $\check{J}_1,\cdots,\check{J}_8$ can be constructed out of the octonion algebra as follows.
Let $e_1,e_2,\cdots,e_8$ be the standard basis of the octonion algebra
${\mathbb O}$ with $e_1$ the multiplicative unit. Let $J_1,J_2,\cdots,J_7$
be the matrix representations of the octonion multiplications by
$e_2,e_2,\cdots,e_8$ on the right over ${\mathbb O}$. Then
\begin{equation}
\check{J}_i=\begin{pmatrix}J_i&0\\0&-J_i\end{pmatrix},1\leq i\leq 7,\qquad
\check{J}_8=\begin{pmatrix}0&I\\-I&0\end{pmatrix}.
\end{equation}

In~\cite[Sections 2.3-2.4]{Ch2}, I calculated their second fundamental forms at specifically chosen points of $M_{+}$, which does not lose generality because of homogeneity, and observed that in the expression~\eqref{0a} we have that $B_a=C_a$ are of rank 1 for $1\leq a\leq 4$ (respectively, $1\leq a\leq6$) in the $(4,5)$ (respectively, $(6,9)$) case, satisfying the property that in each $B_a$ the only nonzero entry is in the last row at column $a$ with value $1/\sqrt{2}$. This peculiarity is not coincidental and is in fact a consequence of codimension 2 estimates.

\begin{lemma}\label{LM} Assume $(m,m_{-})=(4,5)$ or $(6,9)$. Then either the isoparametric hypersurface is the inhomogeneous one of multiplicity pair $(6,9)$ of Ferus, Karcher, and M\"{u}nzner, or
$r_\lambda =1$ for all $\lambda$ in
${\mathcal Q}_{m-1}$.
\end{lemma}

\begin{proof} It is straightforward to see that when $(m,m_{-})=(4,5)$ (respectively, $(m,m_{-})=(6,9)$)
and $j\geq 2$, 
the codimension 2 estimate~\eqref{refinedest} goes through for all $k\leq m-1=3$ (respectively, $k\leq m-1=5$); likewise, the codimension 2 estimate goes through when $k\leq m-2$ for all $j$. Thus, we may assume $k=m-1$ and $j\leq 1$ without loss of generality to consider $\lambda\in{\mathcal Q}_{k-1}={\mathcal Q}_{m-2}\subset {\mathcal Q}_{m-1}$.

Suppose 
\begin{equation}\label{sup}
\sup_{\lambda\in {\mathcal Q}_{m-1}} \, r_\lambda\,\geq\, 2.
\end{equation}
We may so arrange such that generic $\lambda\in {\mathcal Q}_{k-1}={\mathcal Q}_{m-2}$ assumes $r_\lambda\geq 2$.

%The multiplicity pair $(4,5)$
%cannot allow any points of Condition A on
%$M_{+}$.  Hence,
%one of the four pairs of
%matrices $(B_1,C_1),(B_2,C_2),(B_3,C_3)$ and $(B_4,C_4)$ of the shape
%operators $S_{n_1},S_{n_2},S_{n_3}$ and $S_{n_4}$, similar to the one given
%in~\eqref{Matrix}, must be nonzero; we may assume
%one of $(B_1,C_1),(B_2,C_2)$ and $(B_3,C_3)$ is nonzero in the neighborhood of a given
%point, over which generic $\lambda\in{\mathcal Q}_{k-1}$ has $r_\lambda\geq 2$.

%The remark in the first paragraph of this proof reduces the proof to considering those nongeneric $\lambda\in Q_{k-1}$ for which 
%$r_{\lambda}\leq 1$.

\vspace{2mm}

\noindent {\bf Case 1}. On ${\mathcal L}_1$ where $r_\lambda=1$, we have that the codimension 2 estimate still goes through.
This is because~\eqref{W} is now replaced by 

\begin{equation}\label{WW}
\dim(\Theta_j)\leq m_1+m_2+k-3-j=m+m_{-}+k-4
\end{equation}
with $j=1$, due to the fact that such $\lambda$, being nongeneric in ${\mathcal Q}_{k-1}$ as $r_\lambda\geq 2$ for generic $\lambda$,
constitute a subvariety of ${\mathcal Q}_{k-1}$ of dimension at most $k-2$. It follows by~\eqref{Lower-bound} that 
$$
\dim(\Theta_j)\leq \dim(V_k)-2,\quad j=1.
$$

\vspace{1mm}

\noindent {\bf Case 2}. On ${\mathcal L}_0$ where $r_\lambda=0$,~\eqref{WW} now reads
$$
\dim(\Theta_j)\leq m+m_{-}+k-3
$$
with $j=0$. We need to cut back one
more dimension to make the equality in~\eqref{WW} valid.
Since $r_\lambda=0$, we see $B^*_1=C^*_1=0$ and $A^*=I$ in~\eqref{Matrix} for
$S_{n^*_1}$, where $\lambda$ is the 2-plane spanned by $n_0^*$ and $n_1^*$ completed into the basis $n_0^*,n_1^*,n_2^*,\cdots,n_m^*$.
It follows that $p^*_0=0$ and $p^*_1=0$ cut ${\mathscr S}_\lambda$ in the
variety
\begin{equation}\label{free}
\{(x,\pm \sqrt{-1}x,z):\sum_\alpha (x_\alpha)^2=0\}.
\end{equation}
We may assume $(B^*_2,C^*_2)$ is nonzero because of nonexistence of points of Condition A on $M_{+}$.
Since $z$ is a free variable in~\eqref{free},
$p^*_2=0$ will have nontrivial $xz$-terms, 
\begin{eqnarray}\label{p2}
\aligned
0=p^*_2&=\sum_{\alpha p}S_{\alpha p}x_\alpha z_p
+\sum_{\mu p}T_{\mu p}y_\mu z_p+\sum_{\alpha\mu}
U_{\alpha\mu}x_\alpha y_\mu\\
&=\sum_{\alpha p}(S_{\alpha p}\pm\sqrt{-1}T_{\alpha p})x_\alpha z_p,
\endaligned
\end{eqnarray}
taking $y=\pm\sqrt{-1}x$ into account and
remarking that as a result the $x_\alpha y_\mu$ terms are gone because 
$U_{\alpha\mu}=-U_{\mu\alpha}$, where
$$
U_{\alpha\mu}:=\langle S(X_\alpha^*,X_\mu^*),n_2^*\rangle,\;
S_{\alpha p}:=\langle S(X^*_\alpha,Z^*_p),n^*_2\rangle,\; T_{\mu p}:=\langle S(Y^*_\mu,Z^*_p),n^*_2\rangle
$$ 
are entries of $A^*_2$, $B^*_2$ and $C^*_2$,
respectively, and $X^*_{\alpha},1\leq\alpha\leq m_2$,
$Y^*_{\mu},1\leq\mu\leq m_2$ and $Z^*_{p},$ are orthonormal
eigenvectors for the eigenspaces of $S_{n^*_0}$ with eigenvalues $1,-1,$ and $0$,
respectively, as usual. The skew-symmetry of the matrix $U$ comes from the identity
$$
A_jA_1^{tr}+A_1A_j^{tr}+2(B_jB_1^{tr}+B_1B_j^{tr})=0,\quad j\neq 1,
$$ 
obtained by inserting~\eqref{0a} into~\eqref{S^3}, where $(A_1,B_1)$ are as given in Lemma~\ref{le12.1} above.

Hence the dimension                              
of ${\mathscr S}_\lambda$ will be cut down by 2 by $p^*_0,p^*_1,p^*_2=0$, so that again
\begin{equation}\label{EQ2}
\dim(V_{k}\cap{\mathscr S}_\lambda)\leq m+m_{-}-2.
\end{equation} 
%noting that $p^*_0,p^*_1,p^*_2=0$ also cut out $W_k$.
In conclusion, we deduce
\begin{equation}\label{UB}
\dim(\Theta_j)\leq\dim(V_{k}\cap{\mathscr S}_\lambda)+k-2\leq
m+m_{-}+k-4,
\end{equation}
so that the codimension 2 estimate~\eqref{cod2} goes through. The codimension 2 estimate holds true verbatim for the case of multiplicity pair $(6,9)$ with obvious dimension modifications.

Now, the validity of~\eqref{cod2} implies that the 
isoparametric hypersurface is the inhomogeneous one of Ferus, Karcher, and M\"{u}nzner due to Proposition~\ref{isopa}, provided~\eqref{sup} holds.

We therefore conclude that $r_\lambda\leq 1$ for all $\lambda\in {\mathcal Q}_{m-1}$ if the hypersurface is not the inhomogeneous one of Ferus, Karcher, and M\"{u}nzner, in which case we claim that generic $r_\lambda=1$.
%We consider the $(4,5)$ case; the other case is verbatim.
Suppose the contrary. Then $r_\lambda=0$ for all $\lambda$ in ${\mathcal Q}_{m-1}$.
It would follow that 
%$B_1$ of $S_{n_1}$ is identically zero by considering
%$\lambda=[1:\sqrt{-1}:0:0:0]$, because then $B_0^*$ and $B_1^*$ associated with
%$S_{n_0^*}$ and $S_{n_1^*}$ are zero. Likewise, 
$B_a=0$ for all $a\geq 1$, which 
would imply that the isoparametric hypersurface is of Condition A.
This is impossible.

Lastly, for a $\lambda$ with $r_\lambda=0$ we have $A$ in~\eqref{Matrix} is the
identity matrix by~\eqref{A}, so that its rank is full (=5 or 9). It follows that generic $\lambda$
in ${\mathcal Q}_{m-1}$ will have the same full rank property. However, for a $\lambda$
with $r_\lambda=1$, the structure of $\Delta$ in Lemma~\ref{le12.1} implies that
$\Delta=0$ so that
such $A$, which are also generic, will be of rank 4 or 8. This is a contradiction. In other words, $r_\lambda=1$ for all $\lambda\in {\mathcal Q}_{m-1}$.
\end{proof}

\begin{corollary}\label{det}  Suppose $(m,m_{-})=(4,5)$ or $(6,9)$. Then either 
the isoparametric hypersurface is the inhomogeneous one of multiplicity pair $(6,9)$ of Ferus, Karcher and M\"{u}nzner, or
$r_\lambda=1$ for all $\lambda\in {\mathcal Q}_{m-1}$. In the latter case, the second fundamental form of $M_{+}$ is identical with  that of $M_{+}$ of the homogeneous isoparametric hypersurface of the respective multiplicity pair.
\end{corollary}
 %and will only look at the $(4,5)$ case; the other case is verbatim with obvious changes on index ranges. 
Here goes the idea. Assume the hypersurface is not the inhomogeneous one of multiplicity pair $(6,9)$. $r_\lambda=1$ by the preceding proposition.

For a $\lambda$ in ${\mathcal Q}_{m-1}$ spanned by $n_0^*$ and
$n_1^*$, we extend them to a smooth local
orthonormal frame $n_0^*,n_1^*,\cdots,n_{m}^*$ such that $S_{n_0^*}$
is the square matrix on the right hand side
of~\eqref{Matrix} while $S_{n_1^*}$ is the square one on the left hand side,
where $A$ and $B$ are the respective $A_1$ and $B_1$ in Lemma~\ref{le12.1}, for which the skew-symmetric 1-by-1 matrix
$\Delta=0$ in item (2) and so by item(3) the 1-by-1 matrix $\sigma=1/\sqrt{2}$.

Nonexistence of points of Condition A on $M_{+}$ gives us
%a $\lambda_0$ at which $S_{n_2^*}$ in matrix form is such that the
a nonzero matrix $B_2^*$ associated with $S_{n_2^*}$.
%; this property will continue to be true in a neighborhood of $\lambda_0$. 
Modifying ~\eqref{free}, $p^*_0=0$ and $p^*_1=0$ now cut
${\mathscr S}_\lambda$ in the variety
\begin{equation}\label{}
\aligned
&\{(x_1,\cdots,x_4,\frac{t}{\sqrt{2}\iota_\lambda},
\iota_\lambda x_1,\,\cdots,\,\iota_\lambda x_4,\,-\frac{t}{\sqrt{2}\iota_\lambda},\,z_1,\cdots,\,z_3,\,t):\\
&\sum_{j=1}^4 (x_j)^2=0\}
\endaligned
\end{equation}
where $\iota_\lambda=\pm\sqrt{-1}$ as given in~\eqref{iota}, and
\begin{eqnarray}\nonumber
\aligned
x&=(x_1,\,x_2,\,x_3,\,x_4,\,x_5=t/\sqrt{2}\iota_\lambda),\\
y&=(y_1,\,\cdots,\,y_5)=(\iota_\lambda\, x_1,\,\iota_\lambda\, x_2,\,\iota_\lambda\, x_3,\,\iota_\lambda \,x_4,\,-t/\sqrt{2}\iota_\lambda),\\
z&=(z_1,\,z_2,\,z_3,\,z_4=t).
\endaligned
\end{eqnarray}
Meanwhile,~\eqref{p2} becomes

\begin{equation}\label{cut}
\aligned
&0=\sum_{\alpha=1, p=1}^{4,3}(S_{\alpha p}\pm\sqrt{-1}
T_{\alpha p})\,x_\alpha z_p\\
&+ \text{terms not involving}\; x_1,\cdots, x_4, z_1,\cdots,z_3.
%\sum_{p=1}^4(S_{5p}-T_{5p})\,tz_p/(\sqrt{2}\mu_\lambda)\\
%&+\sum_{\alpha\leq 4}\,(S_{\alpha 4}\pm\sqrt{-1}
%T_{\alpha 4}-U_{\alpha 5}/\sqrt{2}\mu_\lambda) \,tx_\alpha
%+\sum_{\mu\leq 4}\,U_{5\mu}\,tx_\mu/\sqrt{2}+U_{55}\,t^2/2.
\endaligned
\end{equation}
%The assumption that $B_2^*$ (or $C_2^*$)
%assumes an
%extra nonzero row other than the last one implies that one more dimension cut can be achieved
Since $x_1,\cdots,x_4,z_1,\cdots,z_3$ are independent variables, the nontriviality of the displayed term on the right hand side of the preceding equation implies that the dimension cut can be reduced 
by 1 so that we have, by~\eqref{EQ2} and~\eqref{UB},
$$
\dim(\Theta_j)\leq m+m_{-}+k-4,\quad j=1,
$$
for $k\leq 3$, so that the codimension 2 estimate~\eqref{cod2} goes through in the neighborhood
of $\lambda$, which is absurd as the hypersurface would then be the inhomogeneous one of the multiplicity pair $(6,9)$ in view of Proposition~\ref{isopa}.
We therefore conclude that 
$(B_j,C_j),j\geq 2,$ around $\lambda$ are of the form
\begin{equation}\label{nlty}
B_j=\begin{pmatrix}0&d_j\\b_j&c_j\end{pmatrix},\quad
C_j=\begin{pmatrix}0&g_j\\e_j&f_j\end{pmatrix},\quad \forall j\geq 2,
\end{equation}
for some real numbers $c_j$ and $f_j$ with 0 of size $4\times 3$ corresponding to the displayed term of~\eqref{cut}.

It turns out, as done in~\cite[Remark 3]{Ch2}, that this suffices to conclude 
that $$d_j=g_j=c_j=f_j=0,\quad\forall j\geq 2.$$ So now we have
$$
A_1=\begin{pmatrix}I&0\\0&0\end{pmatrix},
\quad A_j=\begin{pmatrix}\alpha_j&0\\0&0\end{pmatrix},j=2,3,4,
\quad B_j=C_j=\begin{pmatrix}
0&0\\b_j&0\end{pmatrix},
$$
all of the same block sizes, satisfying
$$
\alpha_j\alpha_k+\alpha_k\alpha_j=-2\delta_{jk}I,\qquad \langle b_j,b_k\rangle=\delta_{jk}/2,
$$
by the identity
\begin{equation}\label{AB}
A_iA_j^{tr}+A_jA_i^{tr}+B_iB_j^{tr}+B_jB_i^{tr}=2\delta_{ij}\, Id.
\end{equation}

As a consequence, first of all, we can perform an orthonormal basis change on
$n_2^*,n_3^*,n_4^*$ so that the resulting new $b_j$ is $1/\sqrt{2}$ at the
$j$th
slot and is zero elsewhere. Meanwhile, we can perform an orthonormal basis
change of the $E_{1}$ and $E_{-1}$ spaces
so that $I$ and $\alpha_j,2\leq j\leq 4,$ are exactly the matrix
representations
of the right multiplication of $1,i,j,k$ on 
${\mathbb H}$ without affecting the row vectors $b_j,2\leq j\leq 4$.
This is precisely the second fundamental form of the homogeneous example.

Now that the second fundamental form of $M_{+}$ is identical with that of the homogeneous example of the same multiplicity pair when the isoparametric hypersurface is not the inhomogeneous one of Ferus, Karcher and M\"{u}nzner,
one can explore the defining equations~\eqref{paqb} through~\eqref{last1} to determine definitively, as was done in~\cite{Ch2}, that the third fundamental form is also nothing other than that of the homogeneous example, most succinctly expressed in the respective quaternionic and octonion framework. The hypersurface is thus the homogeneous one and whence follows the classification. 

\vspace{2mm}

\noindent {\bf The case of multiplicity pair $(7,8)$.} In retrospect, we explored two avenues in the preceding three exceptional cases of multiplicity pairs $(3,4),(4,5)$, and $(6,9)$ to achieve the codimension 2 estimates. The first route was via the right hand side of the {\em a priori} estimate~\eqref{refinedest}, 
$$%\begin{equation}\label{apriori}
m_{-}\geq 2k+1-j-c_j,\quad\forall k\leq m-1,
$$%\end{equation}
where large rank $j$ and nonzero codimension $c_j$ are employed to reduce the situation to $j=0$ or 1, to be followed by the second route to handle the estimate in~\eqref{EQ2}, where, by introducing more cuts via $p_a=0,$ $a\geq 2$, we were able to cut down the upper bound of $\dim(V_k\cap {\mathscr S}_\lambda)$ from the one in~\eqref{Sub} to that in~\eqref{EQ2}, 
$$
\dim(V_{k}\cap{\mathscr S}_\lambda)\leq m+m_{-}-2,
$$
so that the right hand side of~\eqref{refinedest} was down by 1 to achieve
the improved estimate
$$%\begin{equation}\label{cut}
m_{-}\geq 2k-j-c_j, \quad \forall k\leq m-1.
$$%\end{equation} 

The worst case scenario in this procedure is that $p_a=0$, for all $a\geq 2$, always contain ${\mathscr S}_\lambda$ so that no dimension cut can be achieved of $V_k\cap{\mathscr S}_\lambda$, which was avoided in the above three cases with the help of nonexistence of points of Condition A and that the isoparametric hypersurface is not the inhomogeneous one of Ozeki and Takeuchi or of Ferus, Karcher, and M\"{u}nzner.

To understand this worst case, I introduced in~\cite{Chi3} the notion of $r$-nullity:
\begin{definition} We say a normal basis $n_0,n_1,n_2,\cdots,n_m$ is normalized if $S_0,\cdots,S_m$ are as in~\eqref{0a} with $S_1$ normalized as in Lemma~{\rm \ref{le12.1}} for which $B_1$ is of rank $r$.
\end{definition}
\begin{definition}
Given a normalized normal basis $n_0,\cdots, n_m$, let ${\mathbb C}^{m_{-}}\simeq {\mathbb C}E_{+}$,
${\mathbb C}^{m_{-}}\simeq {\mathbb C}E_{-}$ and
${\mathbb C}^{m_{+}}\simeq {\mathbb C}E_{0}$ be parametrized by $x,y$ and $z$, respectively,
where $E_{+},E_{-}$ and $E_0$ are the eigenspaces of $S_0$ with eigenvalues
$1,-1$, and $0$, respectively. Let $x:=(x_1,x_2),y:=(y_1,y_2)$ and $z:=(z_1,z_2)$
with $x_2,y_2,z_2\in{\mathbb C}^r$. Let $p_0, \cdots,p_m$ be the components of the second fundamental form at the base point of $n_0$.

We say a normal basis element $n_l$, $l\geq 2$, is $r$-null if $p_l$
is identically zero when we restrict it to the linear constraints
\begin{equation}\label{cons}
y_1=\iota \, x_1,\quad y_2=-x_2,\quad z_2=\sigma^{-1}(\Delta+\iota\, Id)x_2,
\quad \iota=\pm \sqrt{-1}.
\end{equation}

We say the normal basis is $r$-null if $n_l$ are $r$-null for all $l\geq 2$. 
\end{definition}
\noindent Note that the conditions in~\eqref{cons} define ${\mathscr S}_\lambda$ given in~\eqref{slamb} when $r_\lambda=r$. 

The algebro-geometric definition has a differential-geometric characterization.
\begin{lemma}~\cite[Lemma 3.1]{Chi3} Let $n_0,\cdots,n_m$ be a normalized normal basis.
A normal basis element $n_l,l\geq 2,$ is $r$-null if and only if 
the upper left $(m_{-}-r)$-by-$(m_{+}-r)$ block of $B_l$ and $C_l$
of $S_{l}$ are zero. 
\end{lemma}

It is now clear by the preceding lemma that Condition A is equivalent to that all normalized normal bases are 0-null at the relevant point of $M_{+}$.
Moreover, what we showed in Corollary~\ref{det} for the cases of multiplicity pairs $(4,5)$ or $(6,9)$ is that, all normalized normal bases of $M_{+}$ are 1-null if the isoparametric hypersurface is not the inhomogeneous one of Ferus, Karcher, and M\"{u}nzner, from which we determined below~\eqref{nlty} that the second fundamental form of $M_{+}$ coincides with that of the homogeneous example with the respective multiplicity pair.

What is remarkable is that something similar holds true in the case of multiplicity pair $(7,8)$ as well, only much more complicated this time.

\vspace{2mm}

\begin{proposition}~\cite[Sections 3-6]{Chi3} Assume $(m,m_{-})=(7,8)$ and the isoparametric hypersurface is not the inhomogeneous one constructed by Ozeki and Takeuchi. Then away from points of Condition A, $M_{+}$ is generically $4$-null, i.e., generically chosen normalized normal bases are $4$-null. Moreover, we may assume $A_a$ and $B_a$ are of the form
\begin{equation}\label{mtx}
A_a=\begin{pmatrix}z_a&0\\0&w_a\end{pmatrix},\; B_a=\begin{pmatrix}0&0\\0&c_a\end{pmatrix},\; C_a=\begin{pmatrix}0&0\\0&f_a\end{pmatrix},\;\ 1\leq a\leq 3,
\end{equation}
\begin{equation}\label{mtxx}
A_a=\begin{pmatrix} 0&\beta_a\\\gamma_a&\delta_a\end{pmatrix},\; B_a=\begin{pmatrix}0&d_a\\b_a&c_a\end{pmatrix},\; C_a=\begin{pmatrix} 0&g_a\\b_a&f_a\end{pmatrix},\;4\leq a\leq 7,
\end{equation}
where $A_a$ are of size $8\times 8$, $B_a$ and $C_a$ are of size $8\times 7$, and the lower right blocks of all matrices are of size $4\times 4$.
\end{proposition}

The idea is that if we define
$$
R:=\sup_{\lambda\in {\mathcal Q}_{m-1}} \, r_\lambda,\quad m=7,
$$
then a generically chosen normalized normal basis is $R$-null~\cite[Corollary 3.2]{Chi3}. Furthermore, if $R\geq 5$, then 
the codimension 2 estimate goes through away from points of Condition A~\cite[Sections 3-5]{Chi3}, so that by Proposition~\ref{isopa} the isoparametric hypersurface is the one of Ozeki and Takeuchi. Thus, $R\leq 4$ if it is not the one of Ozeki and Takeuchi, to be assumed from now on. Then at a generically chosen $n_0$ with base point $x_0$, the normal vector $n_0^{\#}:=x_0$ at the ``mirror point'' $x_0^{\#}=n_0$ of $x_0$ in $M_{+}$ is also generic spanning the same $\lambda\in{\mathcal Q}_{m-1}$, $m=7,$  of nullity $R$. We have the explicit dictionary to translate the data of the shape operators at $x$ and $x^{\#}$. Explicitly,
\begin{equation}\label{gOOd}
A_a:=\begin{pmatrix}S^a_{\alpha\mu}\end{pmatrix},\quad B_a:=\begin{pmatrix}S^a_{\alpha p}\end{pmatrix},\quad C_a:=\begin{pmatrix}S^a_{\mu p}\end{pmatrix},\; 1\leq a\leq 7.
\end{equation}
Let the counterpart matrices at $x_0^{\#}$ and their blocks
be denoted by the same notation with an additional $\#$.
Then, 
\begin{equation}\label{gooD}
A_p^\#:=\begin{pmatrix}S^p_{\alpha\mu}\end{pmatrix},\quad B_p^{\#}=\begin{pmatrix}S^a_{\alpha p}\end{pmatrix},
\quad C_p^{\#}=-\begin{pmatrix}S^a_{\mu p}\end{pmatrix}, \quad 8\leq p\leq 14,
\end{equation} 
following the index convention~\eqref{range}. As a consequence, we obtain many zeros as indicated in~\eqref{mtx} and~\eqref{mtxx}, and in particular, we obtain $R=4$ by a resulting Clifford representation~\cite[Proposition 6.1]{Chi3}.

Though a first glance at~\eqref{mtx} and~\eqref{mtxx} suggests that it would still be a long way home to determine the second fundamental form of $M_{+}$, in sharp contrast with the case of multiplicity pairs $(4,5)$ and $(6,9)$, where $1$-nullity rather quickly pins it down. However, a crucial observation is that
the matrices $\begin{pmatrix}\sqrt{2}c_a&w_a\end{pmatrix},\, 1\leq a\leq 3,$ form a 
Clifford multiplication of type $[3,4,8]$,
\begin{equation}\label{eq10000}
F:{\mathbb R}^3\times {\mathbb R}^4\rightarrow{\mathbb R}^8,\quad F(u_a,v_\alpha)=\;\text{the}\;\alpha\text{th row of}\;\begin{pmatrix}\sqrt{2}c_a&w_a\end{pmatrix},
\end{equation}
for orthonormal bases $u_a$ and $v_\alpha$, satisfying $|F(x,y)|=|x||y|$ derived from~\eqref{AB}. 

A second crucial observation is that, with the conversion~\eqref{mtx} through~\eqref{gooD}, the first columns of\, $b_4,\cdots,b_7$ are, respectively, the first, second, third, and fourth columns of $c_1^\#$. Similarly, the second {\rm (}vs. third{\rm )} columns of $b_4,\cdots,b_7$ are the respective columns of $c_2^\#$ {\rm (}vs. $c_3^\#${\rm )}. For instance, at $x_0^{\#}$, the normalized $c_1^{\#}$ is
$$
c_1^\#=\begin{pmatrix}\sigma_1&0&0&0\\0&\sigma_1&0&0\\0&0&\sigma_2&0\\0&0&0&\sigma_2\end{pmatrix},
$$
by Lemma~\ref{le12.1} due to 4-nullity at $x_0^{\#}$. Hence, the first columns of $b_4,\cdots,b_7$ are, respectively, 
$$
(\sigma_1,0,0,0)^{tr},\;\; (0,\sigma_1,0,0)^{tr},\;\;(0,0,\sigma_2,0)^{tr},\;\;(0,0,0,\sigma_2)^{tr},
$$
from which a third crucial observation can be drawn. Indeed, it was shown in~\cite[Lemma 7.1]{Chi3} that a generic linear combination
$$
b(x):=x_1b_4+\cdots+x_4b_7
$$ 
is of rank no more than 2, so that we may in fact assume that all of $b_4,\cdots,b_7$ have zero third column, by the fact that the Koszul complex
$$
0\longrightarrow R\stackrel{x\wedge}\longrightarrow \Lambda^1 R^4\stackrel{x\wedge}\longrightarrow\Lambda^2R^4\stackrel{x\wedge}\longrightarrow\Lambda^3R^4\stackrel{x\wedge}\longrightarrow \Lambda^4 R^4\rightarrow 0,
$$
where $R:={\mathbb R}[x_1,x_2,x_3,x_4]$ is the polynomial ring in four variables and  $x\wedge$ means taking the wedge product against $x$, is a free resolution~\cite[Chapter 17]{Ei}. The assumption that $b(x)$ is generically of rank 2 means that the wedge product of the second column $v_2$ and third column $v_3$ of $b(x)$ lives in the kernel of
$$
\longrightarrow\Lambda^2R^4\stackrel{x\wedge}\longrightarrow\Lambda^3R^4,\quad v_2\wedge v_3\mapsto x\wedge(v_2\wedge v_3)=0,
$$
so that either $v_2\wedge v_3=0$, in which case they differ by a constant multiple, or, $v_2\wedge v_3=x\wedge w$ for some $w\in R^4$, so that we may assume the first two columns of $b(x)$ are both $x$ up to a constant multiple. As a consequence, the conversion says that we may assume $c_3^{\#}=0$. % fits~\eqref{eq10000}.

In~\cite[Section 5]{CW}, $F$ in~\eqref{eq10000} with the constraint $c_3^{\#}=0$ were classified and given by
%We summarize before we proceed further. When both $x$ and $x^\#$ are generic in $M_{+}$ with the chosen $4$-nullity bases as specified in Remark~$\ref{rkk}$, we have~\eqref{mtx}
%where, interchanging $x$ and $x^\#$ by symmetry, we may assume, by~\eqref{8000} and~\eqref{NQE*}, that
$$ 
c_1^{\#}=\epsilon\, I,\quad c_3^{\#}=0
$$
for some $\epsilon>0$. 
%%Item {\rm (}$2${\rm )} of Corollary~$\ref{inde}$ then implies that all $d_a^{\#}\neq 0,4\leq a\leq 7,$ because now $B_3^{\#}=0$ by~\eqref{mtx} and the first four rows of $B_1^{\#}$ and $B_2^{\#}$ are zero. 
%%Also, since $C_3^{tr}C_3=B_3^{tr}B_3$, we obtain $C_3=0$ so that a similar situation holds for $C_a,1\leq a\leq 7$ as well. 
%Moreover, the third columns of the four $4$-by-$3$ matrices $b_4,\cdots,b_7$ at $x$ are zero in accordance with $c_3^\#=0$. By~\eqref{NQE*} applied at $x^{\#}$, the matrix
%%~\eqref{comb} in Appendix II gives the explicit form of a linear combination of
 $c_2^\#$ is of the form
 \begin{equation}\label{NQE}
 c_2^\#=a\,Id+b\begin{pmatrix}I&0\\0&\pm I\end{pmatrix},\quad I=\begin{pmatrix}0&-1\\1&0\end{pmatrix},\;\; b\neq 0,
 \end{equation}
for some $a$ and $b$. By conversion,
%with $c_1^\#=\epsilon\, Id$ and $c_3^\#=0$, the three matrices can be converted, by {\bf NOTE} in Lemma~{\rm \ref{cccor}}, to the data
\small\begin{equation}\label{EQQQ}
b_4=\begin{pmatrix}\epsilon&a&0\\0&b&0\\0&0&0\\0&0&0\end{pmatrix}, \,b_5=\begin{pmatrix}0&-b&0\\\epsilon&a&0\\0&0&0\\0&0&0\end{pmatrix},\, b_6=\begin{pmatrix}0&0&0\\0&0&0\\\epsilon&a&0\\0&\pm b&0\end{pmatrix},\, b_7=\begin{pmatrix}0&0&0\\0&0&0\\0&\mp b&0\\\epsilon&a&0\end{pmatrix}
\end{equation}\normalsize
at $x_0$, whose linear combinations are of generic rank $2$. 

In particular, a glance at 
$B_a,1\leq a\leq 7,$ in~\eqref{mtx} shows that their third columns are all zero, or equivalently, that there is a common kernel of generic dimension 1 for all the shape operators $S_n$ for all $n$. 
This gives us a clear geometric picture:

\vspace{1mm}

{\em  When the isoparametric hypersurface with multiplicities $(m,m_{-})=(7,8)$ is not the one constructed by Ozeki and Takeuchi, consider the quadric ${\mathcal Q}_6$ of oriented $2$-planes in the normal space at a generic point $x_0\in M_{+}$. We know a generic element $(n_0,n_1)$ in ${\mathcal Q}_6$ is $4$-null, or equivalently, the intersection $V$ of the kernels of $S_{n_0}$ and $S_{n_1}$ is $3$-dimensional. By the preceding lemma, there is a nonzero unit vector $v\in V$ common to all kernels of the shape operators at $x_0$. We choose an orthonormal basis $e_1,e_2,e_3=v$ spanning $V$.
When viewed at the mirror point $x_0^\#=n_0\in M_{+}$, $e_1,e_2,e_3$ are converted to three normal basis vectors of which the three matrices $c_1^\#,c_2^\#,c_3^\#$ in~\eqref{mtx} are of the form $c_1^\#=\epsilon\,Id, c_3^\#=0$, and $c_2^\#$ given in~\eqref{NQE}.
By a symmetric reasoning, all this holds true as well at $x_0$ when both $x_0$ and $x_0^\#$ are generic.}

\vspace{1mm}

We are ready to see that the isoparametric hypersurface is one of the two constructed by Ferus, Karcher, and M\"{u}nzner with the Clifford action operating on $M_{-}$.  To this end, for a normal basis $n_0,\cdots,n_m$ of $M_{+}$ with base point $x_0$, let us set
\begin{equation}\label{*}
x_0^*:=(x_0+n_0)/\sqrt{2},\quad n_0^*:=(x_0-n_0)/\sqrt{2}.
\end{equation}
$x_0^*$ is a point on $M_{-}$ and $n_0^*$ is normal to $M_{-}$ at $x_0^*$. The normal 
space to $M_{-}$ at $x_0^*$ is
${\mathbb R}n_0^*\oplus E_{1}.$ 
Furthermore, the $(+1)$-eigenspace $E_{1}^*$
of the shape operator $S_{n_0^*}$ is spanned by $n_1,\cdots,n_{m_{+}}$, the $(-1)$-eigenspace $E_{-1}^*$ of $S_{n_0^*}$ is $E_0$, and
the $0$-eigenspace $E_0^*$ of $S_{n_0^*}$ is $E_{-1}$.

Referring to~\eqref{0a}, let the counterpart matrices at $x_0^*$ and their blocks
be denoted by the same notation with an additional *.
Then, for $\alpha=1,\cdots,m_{-},$

\begin{equation}\label{Good}
A_\alpha^*=-\sqrt{2}\begin{pmatrix}S^a_{\alpha p}\end{pmatrix},\quad
B_\alpha^*=-1/\sqrt{2}\begin{pmatrix}S^a_{\alpha\mu}\end{pmatrix},\quad
C_\alpha^*=-1/\sqrt{2}\begin{pmatrix}S^p_{\alpha\mu}\end{pmatrix}.
\end{equation}

There follows from~\eqref{NQE} important features for the block matrices in~\eqref{mtx} and~\eqref{mtxx}.

\vspace{2mm}

\begin{lemma}\label{ammel}~\cite[Lemma 7.6, Corollary 7.1]{Chi3} After a frame change over $E_1$ and $E_{-1}$ we have the following.
\begin{description}
\item[(1)] The spectral data $(\sigma,\Delta)$ given in Lemma~{\rm \ref{le12.1}} equal $(Id/\sqrt{2}, 0)$.
\item[(2)] 
$$
d_a=\begin{pmatrix}d_1^a\\0\end{pmatrix},\quad g_a=\begin{pmatrix}0\\g_1^a\end{pmatrix},\quad a=1,\cdots, 7,
$$
where 0 and $d_1^a=g_1^a$ are of size $2\times 4$.
\item[(3)] $c_a=f_a$, $\beta_a=(\gamma_a)^{tr}$, and $\delta_a$ is skew-symmetric for all $1\leq a\leq 7$.
\end{description}
\end{lemma}
%\begin{corollary} With the data in Lemma~\ref{ammel}, we have .
%\end{corollary}

%The corollary follows from $\Delta=0$ and the identity 
%$$
%A_iC_jB_j^{tr}+B_iC_j^{tr}A_j^{tr}+A_jC_iB_j^{tr}\quad\text{is skew-symmetric}
%$$
%derived from that the shape operator $S_n$ satisfies $(S_n)^3=S_n$, where we set $i\geq 4$ and $j=1$ in~\eqref{mtxx}; we obtain $\beta_i=\gamma^{tr}_i$ for $i\geq 4$ (it is trivially true for $1\leq i\leq 3$.) 

%\noindent {\bf NOTE:}: {\em When viewed at $x^\#$ the first columns of\, $b_4,\cdots,b_7$ are, respectively, the first, second, third, and fourth columns of $c_1^\#$, i.e.,
%$$
%c_1^\#=\begin{pmatrix}\sigma_1&0&0&0\\0&\sigma_1&0&0\\0&0&\sigma_2&0\\0&0&0&\sigma_2\end{pmatrix},
%$$
%in~\eqref{sharp}. Similarly, the second {\rm (}vs. third{\rm )} columns of $b_4,\cdots,b_7$ are the respective columns of $c_2^\#$ {\rm (}vs. $c_3^\#${\rm )}. }

%(Likewise, there are counterpart matrices when we replace $\alpha$ by $\mu$ at the points $(x_0^*)^\#\in M_{-}$.)
Proposition~\ref{pivo}, in which the equations~\eqref{3eq} and~\eqref{4theq} are modified with appropriate index changes by the recipe~\eqref{Good}, characterizes the Ferus-Karcher-M\"{u}nzner examples whose Clifford action operates on $M_{-}$. It
reads, in view of the notation of~\eqref{frame} and~\eqref{FS},
\begin{equation}\label{cliff}
\aligned
&A_\alpha^*=A_\mu^*,\\
&(a,\mu)\;\text{entry of }\; B_\alpha^*=-(a,\alpha)\;\text{entry of}\; B_\mu^*,\\
&(p,\mu)\;\text{entry of }\; C_\alpha^*=-(p,\alpha)\;\text{entry of}\; C_\mu^*,\\
&\omega^i_j-\omega^{i'}_{j'}=\sum_k L^i_{jk}\, (\theta^k+\theta^{k'}),
\endaligned
\end{equation}
for some smooth functions $L^i_{jk}$, where $i,j,k$ are in the $\alpha$ index range and $i',j',k'$ are in the $\mu$ index range with the respective index values (i.e., $i$ indicates $\alpha=i$ and $i'$ indicates $\mu=i+m_{-}$, etc.), recalling~\eqref{range} through~\eqref{FS}.

In fact, employing~\eqref{Good}, the first three equations in~\eqref{cliff} for $M_{-}$ take the form
\begin{equation}\label{EQ}
\aligned
&B_a=C_a,\; \forall a,\\
&A_a\;\text{is skew-symmetric},\;\forall a,\\
&A_a^\#\;\text{is skew-symmetric},\;\forall a,
\endaligned
\end{equation}
over $M_{+}$, which is exactly Lemma~\ref{ammel} after a slight frame change by swapping rows. Note, at $x_0^*$, we can now change the sign of the last four $\alpha$-rows of $A_i$ without affecting the skew-symmetry of $\delta_i$ and the property $d_i=g_i,c_i=f_i$, so that now
$\beta_i=\gamma_i^{tr},\quad 1\leq i\leq 7,$ at $x$ are converted to satisfy the second and third skew-symmetric conditions in~\eqref{cliff} at $x_0^*$.

It remains to establish the fourth equation in~\eqref{cliff}. A slight modification of~\cite[Lemma 2, p. 11]{Ch1}, the last item holds true if either $\alpha=i$ or $\alpha=j$ indexes a basis vector in the image of the linear map
\begin{equation}\label{H}
H:{\mathcal E}_{+}^*\oplus {\mathcal E}_{-}^*\rightarrow {\mathcal E}_0^*,\quad (e_a,e_p)\mapsto \sum_{\alpha}S^a_{\alpha p}e_\alpha,
\end{equation}
which is easily seen to be the direct sum of all $e_{\alpha=l}$ for $l\neq 3, 4$ (i.e., the 3rd and 4th rows of $B_a$ are zero for all $1\leq a\leq 7$). Thus, we need only show that the last item of~\eqref{cliff} is valid for $i=3,j=4$ in the $\alpha$-range. Referring to the discussions around~\eqref{del}, now adopted for $M_{-}$, since the distribution $\Delta$ is the kernel of~$\theta^a+\theta^{a+m}, a=1, \cdots,m,$ we see the right hand side of the fourth equation of~\eqref{cliff} is automatically zero over $\Delta$. Thus to establish the identity, it suffices to show that $\omega^3_4-\omega^{3'}_{4'}$ on the left hand side annihilates $\Delta$, i.e., for $v:=e_{l'}-e_l\in{\mathcal F}$ we must verify $\omega^3_4(v)=\omega^{3'}_{4'}(v)$, or,
\begin{equation}\label{4v}
\omega^3_4(e_l)=-\omega^{3'}_{4'}(e_{l'}), \quad l=1,\cdots, 7,
\end{equation}
 since $F^i_{jk}=0$ whenever exactly two indexes fall in the same $\alpha,\mu,a$, or $p$ range~\cite[(2.9), p. 9]{CCJ}.

%Since the calculation is pointwise, we first look at the geometry before we proceed. 
Now, for $x\in M_{+}$ and $n$ in the unit normal sphere to $M_{+}$ at $x$, the map
\begin{equation}\label{diffeo}
f:(x,n)\mapsto (x^*,n^*)=((x+n)/\sqrt{2},(x-n)/\sqrt{2})
\end{equation}
sets up a diffeomorphism between the normal bundles of $M_{+}$ and $M_{-}$. Fix a point $(x_0,n_0)$ in the unit normal bundle of $M_{+}$, consider the two sets  
$$
S_{+}:=\{(x,n): x+n=x_0+n_0\}, \quad S_{-}:=\{(x,n):x-n=x_0-n_0 \}.
$$
$S_{\pm}$ are two 8-dimensional spheres, which can be seen by taking derivative of $x\pm n=c$ for a constant $c$, whose typical tangent space to $S_{\pm}$ is the eigenspace ${\mathcal E}_{\pm}$ at $(x,n)$, respectively. 

The diffeomorphism $f$ maps $S_{+}$ to a sphere whose tangent space at $(x_0^*,n_0^*)$ is the vertical ${\mathcal V}^*$ of the unit normal bundle of $M_{-}$ because
$$
f: (x,\,n)\in S_{+}\longmapsto (c/\sqrt{2},\,c/\sqrt{2}-\sqrt{2}n),
$$
so that it is the fiber of the unit normal bundle of $M_{-}$ over $x_0^*$; likewise, $f$ maps $S_{-}$ to a sphere whose tangent space at $(x_0^*,n_0^*)$ is the horizontal ${\mathcal E}_0^*$ because 
$$
f:(x,\,n)\in S_{-}\rightarrow (-c/\sqrt{2}+\sqrt{2}x,\,c/\sqrt{2}).
$$
Thus to calculate the quantities in~\eqref{4v}, it suffices to observe that~\eqref{H} gives us the information
$$
\dim(\bigcap_{a=1}^7\text{kernel}(B_a^{tr}))=2,
$$
which is a consequence of item (2) of Lemma~\ref{ammel}. This translates to $S_{+}$ to say that the tangent space to $S_{+}$ at $(x,n)$ is identified with $E_+$ of the second fundamental form $S_n$, in which there naturally sits a 2-dimensional plane that is the
intersection of all kernels of the $B_m^{tr}$-block of $S_m$, with $m$ perpendicular to $n$ at $x$, which forms a 2-plane bundle ${\mathcal P}_{+}$ over $S_{+}$. By the same token there is a 2-plane bundle
${\mathcal P}_{-}$ over $S_{-}$ which comes from the intersection of all kernels of the $C_m^{tr}$-block of $S_m$ with $m$ perpendicular to $n$ at $x$. %Now, the above fact that after swapping rows
Now, since $d_a=g_a,1\leq a\leq 7$, means that ${\mathcal P}_{+}$ and ${\mathcal P}_{-}$ are parametrized identically in the coordinates, once we set up the coordinate system of the ambient Euclidean space by the eigenspace decomposition
$$
{\mathbb R}x\oplus{\mathbb R}n\oplus E_0\oplus E_{+}\oplus E_{-}
$$
of the shape operator $S_n$ at $x$ for $(x,n)\in S_{+}$, where the third and fourth rows of $B_a$ are zero for all $1\leq a\leq 7$. %we may assume, after swapping the third and fourth rows with the first and second, that . 
%${\mathcal P}_{+}$ and ${\mathcal P}_{-}$ are parametrized identically in the coordinates
%In other words, in the coordinates we can parametrize $S_{+}$ and $S_{-}$ via an isometry $\iota$ in which ${\mathcal P}_{+}$ is brought to ${\mathcal P}_{-}$. 
As a consequence, via the diffeomorphism $f$ in~\eqref{diffeo}, a local basis $(e_3,e_4)$ spanning ${\mathcal P}_{+}$ is converted to one on the image sphere whose tangent space at $(x_0^*,n_0^*)$ is ${\mathcal V}$, and local basis $(e_{3'},e_{4'})$ spanning ${\mathcal P}_{-}$ is converted to one on the image sphere whose tangent space at $(x_0^*,n_0^*)$ is ${\mathcal E}^*_0$. Thus on the image sphere we derive %through the isometry $\iota$ we see that
$$
\omega^3_4(e_l)=-\langle de_3(e_l),e_4\rangle=\langle de_{3'}(e_{l'}),e_{4'}\rangle=\omega^{3'}_{4'}(e_{l'}),
$$
which gives~\eqref{4v}, remarking that the negative sign in the first equality is a result of the sign convention~\eqref{signconvention}. %remarking that there the extra sign is a result of the sign convention in our identification map $Q$ in~\eqref{Qsign}, whose choice is in agreement with that of an isoparametriic hypersurface constructed by Ferus, Karcher, and M\"{u}nzner. 

The four equations in~\eqref{cliff} are satisfied. Thus the isoparametric hypersurface is one of the two constructed by Ferus, Karcher, and M\"{u}nzner, if it is not the one constructed by Ozeki and Takeuchi.  

Lastly, for $g=6$, Miyaoka classified~\cite{Mi2, Mi4} (see also~\cite{Mi1}) the case when the multiplicity pair is $(2,2)$. It is homogeneous. The key to her proof is to show that Condition A holds at all points on either focal manifold. 

The classification of isoparametric hypersurfaces in the sphere has thus been completed.

\section{A few questions} There remain a few fundamental questions in the spherical case that I find especially interesting to be listed here to conclude the article. 

Is there a geometric way to prove that the number $g$ of principal curvatures is $1,2,3,4,$ or 6 ? M\"{u}nzner's proof~\cite[II]{Mu} is topological. A recent paper of Fang~\cite{Fa1} gave another topological proof of M\"{u}nzner's result. 

On the other hand, is there a geometric proof for the multiplicity pairs $(m_{+},m_{-})$, $m_{+}\leq m_{-}$, when $g=4$? Though Stolz's approach is topological and works for more general compact proper Dupin hypersurfaces, our classification in~\cite{CCJ} exhausts the multiplicity pairs so long as $m_{-}\geq 2m_{+}-1$. So, the question to ask is whether there is a geometric way to show that the multiplicity pair must be $(2,2), (4,5), (3,4), (6,9),$ or $(7,8)$ when $m_{+}\leq m_{-}\leq 2m_{+}-2$? Similarly, is there a geometric proof that $m_{+}=m_{-}=1,$ or $2$ when $g=6$?

Immervoll~\cite{Im} gave a different proof of the result in~\cite{CCJ} by isoparametric triple systems that Dorfmeister and Neher developed~\cite{DN2}. Is there a classification of the exceptional cases by the approach of isoparametric triple systems?

For $g=6$, derive an expansion formula for the Cartan-M\"{u}nzner polynomial similar to the one done by Cartan for $g=3$ and by Ozeki and Takeuchi for $g=4$. In addition to Miyaoka's geometric proof of homogeneity of such isoparametric hypersurfaces, is there a proof by utilizing the expansion formula similar to the one given by Cartan for the case $g=3$ that also enjoys equal multiplicity? As Miyaoka's proof pointed out, Condition A should play a decisive role.

\vspace{6mm}

\noindent {\em Department of Mathematics and Statistics, Washington University, St. Louis, MO $63130$; chi@wustl.edu}


\begin{thebibliography}{100}
%\bibitem{Ad} N.\ Addington, {\em Note on Spin}, pages.uoregon.edu/adding/notes/accidental.pdf
\bibitem{Ar} U.\ Abresch, {\em Isoparametric hypersurfaces with four or six principal
curvatures,} {Math.\ Ann.} {\bf 264}(1983), 283-302.
%\bibitem{Be} M.\ Berger, {Les vari\'{e}t\'{e}s riemanniennes homog\`{e}nes normales simplement connexes \`{a} courbure strictement positive}, Ann.\ Scuola Norm.\ Sup.\ Pisa {\bf 15}(1961), 179-246.
\bibitem{BCO} J.\ Berndt, S.\ Console and C.\ Olmos, {Submanifolds and Holonomy},
Chamoman and Hall/CRC, Boca Rotan, London, New York, Washington D.C., 2003
\bibitem{Ca1} \'{E}.\ Cartan, {\em Familles de surfaces isoparam\'{e}triques dans
les espaces \`{a} courbure constante}, {Ann.\ di Mat.} {\bf 17}(1938), 177-191.
\bibitem{Ca2} \'{E}.\ Cartan, {\em Sur des familles remarquables d'hypersurfaces
isoparam\'{e}triques dans les espaces sph\'{e}riques}, {Math.\ Z.} {\bf 45}(1939), 335-367.
\bibitem{Ca3} \'{E}.\ Cartan, {\em Sur quelques familles remarquables d'hypersurfaces},
{C.\ R.\ Congr\`{e}s Math.\ Li\`{e}ge} (1939), 30-41.
\bibitem{Ca4} \'{E}.\ Cartan, {\em Sur des familles d'hypersurfaces isoparam\'{e}triques
des espaces sph\'{e}riques \`{a} $5$ et \`{a} $9$ dimensions}, {Revista Univ.\ Tucum\'{a}n}{\bf 1}
(1940), 5-22.
%\bibitem{C} T.\ Cecil, {\em Taut and Dupin Submanifolds}, {Tight and Taut Submanifolds}, MSRI publications {\bf 32}(1997), 135-180, Cambridge University Press, New York.
\bibitem{CCJ} {T. Cecil, Q.-S. Chi and G. Jensen}, {\em Isoparametric hypersurfaces with four
principal curvatures,}\, {Ann. Math.}\,{\bf 166}(2007), 1-76.
\bibitem{CR} T. Cecil and P. Ryan, {Geometry of Hypersurfaces}, Springer-Verlag, New York, 2015.
\bibitem{Ch} C.\ Chevalley, {\em Invariants of finite groups generated by reflections},
Amer.\ J.\ Math. {\bf 77}(1955), 778-782.
\bibitem{Chi} Q.-S.\ Chi, {\em Isoparametric hypersurfaces with four principal curvatures
revisited}, Nagoya Math. J. {\bf 193}(2009), 129-154. 
\bibitem{Ch1} Q.-S.\ Chi, {\em Isoparametric hypersurfaces with four principal curvatures,
II}, Nagoya Math.\ J. {\bf 204}(2011), 1-18.
\bibitem{Chichi} Q.-S.\ Chi, {\em A new look at Condition A}, Osaka J.\ Math. {\bf 49}(2012), 133-166.
\bibitem{Ch2} Q.-S.\ Chi, {\em Isoparametric hypersurfaces with four principal curvatures, III}, {J.\ Diff.\ Geom.} {\bf 94}(2013), 469-504.
\bibitem{Chi4} Q.-S.\ Chi, {\em Ideal theory and classification of isoparametric hypersurfaces}, Contemporary Mathematics {\bf 646}(2015), 81-104, A. M. S.
\bibitem{Chi3} Q.-S.\ Chi, {\em Isoparametric hypersurfaces with four principal curvatures, IV}, to appear in J.\ Diff.\ Geom.
\bibitem{Chi6} Q.-S.\ Chi, {\em Classification of isoparametric hypersurfaces}, Proceedings of the Sixth International Congress of Chinese Mathematicians, ALM {\bf 36}(2017), 437-451. 
\bibitem{CW} Q.-S. Chi and Haiyang Wang, {\em Orthogonal multiplications of type $[ 3, 4, p], p\leq 12$}, Beitr\"{a}ger zur Algebra und Geometrie/Contributions to Algebra and Geometry {\bf 59}(2018), 167-197.
\bibitem{CO} S.\ Console and C. Olmos, {\em Clifford systems, algebraically constant second fundamental form and isoparametric hypersurfaces}, Manuscripta Math. {\bf 97}(1998), 335-342.
%\bibitem{De} O.\ Dearricott, {\em $n$-Sasakian manifolds}, T\^{o}hoku Math.\ J. {\bf 60}(2008), 329-347.
 % , MR3761408, Zbl 06844623.
\bibitem{D} J. Dadok, {\em Polar coordinates induced by actions of compact Lie groups}, Trans. Amer. Math. Soc. {\bf 288}(1985), 125-137.
\bibitem{DN2} J.\ Dorfmeister and E.\ Neher, {\em An algebraic approach to isoparametric hypersurfaces in spheres} I and II, {T\^{o}hoku Math.\ J.} {\bf 35}(1983), 187-224 and {\bf 35}(1983), 225-247.
\bibitem{DN1} J.\ Dorfmeister and E.\ Neher, {\it Isoparametric triple systems of algebra type}, Osaka J.\ Math. {\bf 20}(1983), 145-175.
\bibitem{DN} J.\ Dorfmeister and E.\ Neher, {\it Isoparametric hypersurfaces, case $g = 6, m = 1$}, Communications in Algebra {\bf 13} (1985), 2299-2368.
\bibitem{Ei} D.\ Eisenbud, {Commutative Algebra with a View Toward Algebraic Geometry}, Springer-Verlag, New York, 1995.
\bibitem{EH} J.-H. Eschenberg and E. Heintze, {\em Polar representations and symmetric spaces}, J. reine angew. Math. {\bf 507}(1999), 93-106.
\bibitem{EH1} J.-H. Eschenberg and E. Heintze, {\em On the classification of polar representations}, Math. Z. {\bf 232}(1999), 391-398.
\bibitem{Fa0} F. Fang, {\em Multiplicities of principal curvatures of isoparametric hypersurfaces}, Preprint, Max Planck Institut f�r Mathematik (1996).
\bibitem{Fa} F. Fang, {\em On the topology of isoparametric hypersurfaces with four distinct principal curvatures}, Proc.\ A. M. S. {\bf 127}(1999), 259-264. 
\bibitem{Fa2} F. Fang, {\em Topology of Dupin hypersurfaces with six principal curvatures}, Math. Z.{\bf 231}(1999), 533-555.
\bibitem{Fa1} F. Fang, {\em Dual submanifolds in rational homology spheres}, Sci. China Math. {\bf  60}(2017), 1549-1560.
\bibitem{FKM} {D. Ferus, H. Karcher and H.-F. M\"{u}nzner}, {\em Cliffordalgebren
und neue isoparametrische hyperfl\"{a}schen,}\, {Math. Z.}\,{\bf 177}(1981), 479-502.
%\bibitem{GT} J.\ Ge and Z.\ Tang, {\em Isoparametric functions and exotic spheres}, {to appear in J.\ Reine Angew.\ Math.}
%\bibitem{GT1} J.\ Ge and Z.\ Tang, {\em Geometry of isoparameric hypersurfaces in Riemannian manifolds}, preprint, arXiv:1066.2577.v1.
\bibitem{GH} K. Grove, S. Halperin, {\em Dupin hypersurfaces, group actions and the double mapping cylinder}, J. Diff.\ Geom. {\bf 26}(1987), 429-459.
\bibitem{He} S.\ Helgason, {Differential Geometry, Lie Groups, and Symmetric Spaces}, Academic Press, Inc., 1978.
\bibitem{HL} W.\ Y.\ Hsiang and H.\ B.\ Lawson, {\em Minimal submanifolds with low cohomogeneity}, J.\ Diff. Geom. {\bf 5}(1971), 1-38.
\bibitem{Hu} {D. Husemoller}, {Fiber Bundles}, 3rd, ed., Graduate Text in Mathematics 20, Springer-Verlag, New York, Berlin, Heidelberg. 
\bibitem{Im} S.\ Immervoll, {\em The classification of isoparametric hypersurfaces with four distinct principal curvatures}, {Ann.\ Math.}, {\bf 168}(2008), 1011-1024.
\bibitem{Ka} H.\  Karcher, {\em A geometric classification of positively curved symmetric spaces and the isoparametric construction of the Cayley plane}, Ast\'{e}risque {\bf 163-164}(1988), 111-135.
\bibitem{KK} N.\ Knarr and L.\ Kramer, {\em Projective planes and isoparametric hypersurfaces}, Geometriae Dedicata {\bf 58}(1995), 193-202.
\bibitem{KN} S.\ Kobayashi and K.\ Nomizu, {Foundations of Differential Geometry},
Vol. I, 1963, Vol. II, 1969, John Wiley \& Sons, New York, London.
\bibitem{Ku} E.\ Kunz, Introduction to Commutative Algebra and Algebraic Geometry, Birkh\"{a}user, Boston, 1985.
\bibitem{La} E.\ Laura, {\em Sopra la propagazione di onde in n mezzo indefinito}, In: {Scritti matamatici offerti ad Enrico D'Ovidio}, 253-278, Bocca, Torino, 1918.
\bibitem{Le} T.\ Levi-Civita, {\em Famiglie di superficie isoparametriche nell'ordinario spazio euclide}, Rend. Aca. Naz. Lincei {\bf XXVI}(1937), 355-362.
\bibitem{Lo} O.\ Loos, {Symmetric Spaces}, Vols. I, II, W. A. Benjamin, Inc., New York, Amsterdam, 1969.
%\bibitem{MO1} H.\ Ma and Y.\ Ohnita, {\em On Lagrangian submanifolds in complex hyperquadrics and isoparametric hypersurfaces in spheres}, {Math.\ Z.} {\bf 261}(2009), 749-785.
%\bibitem{MO2} H.\ Ma and Y.\ Ohnita, {\em Hamiltonian stability of the Gauss images of homogeneous isoparametric hypersurfaces. I}, {J.\ Diff.\ Geom.} {\bf 97}(2014), 275-348.
\bibitem{Mi1} {R. Miyaoka}, {\em The Dorfmeister-Neher theorem on isoparametric hypersurfaces},\,
{Osaka J. Math.}\,{\bf 46}(2009), 695-715.
%\bibitem{Mi3} {R. Miyaoka}, {\em Transnormal functions on a Riemannian manifold}, {Diff.\ Geom.\ Appls.} {\bf 31}(2013), 130-139.
\bibitem{Mi2} {R. Miyaoka}, {\em Isoparametric hypersurfaces with $(g,m)=(6,2)$}, {Ann.\ Math.} {\bf 177}(2013), 53-110.
\bibitem{Mi4} {R. Miyaoka}, {Errata of ``Isoparametric hypersurfaces with (g, m) = (6, 2)''}, {Ann.\ Math.} {\bf 183} (2016) 1057-1071.
\bibitem{Mu} H.\ F.\ M\"{u}nzner, {\em Isoparametrische Hyperfl\"{a}chen in Sph\"{a}ren}, I, {Math.\ Ann.} {\bf 251}(1980), II, {\bf 256}(1981).
\bibitem{No1} K.\ Nomizu, {\em Some results in E. Cartan's theory of isoparametric families of hypersurfaces}, Bull.\ Amer.\ Math.\ Soc. {\bf 79}(1973), 1184-1188. 
\bibitem{No2} K.\ Nomizu, {\em \'{E}lie Cartan's work on isoparametric families of hypersurfaces}, Proc.\ Sympos.\  Pure Math. {\bf 27}(1975), 191-200.
\bibitem{On} B.\ O'Neill, {Elementary Differential Geometry}, Academic Press, New York and London, 1966.
\bibitem{OT} H.\ Ozeki and M.\ Takeuchi, {\em On some types of isoparametric hypersurfaces in spheres}, I, T\^{o}hoku Math.\ J. {\bf 27}(1975), 515-559, and II, {\bf 28}(1976), 7-55. 
\bibitem{Ry} P.\ Ryan, {\em Homogeneity and some curvature conditions for hypersurfaces}, T\^{o}hoku Math.\ J. {\bf 21}(1969), 363-388.
\bibitem{Se1} B.\ Segre, {\em Una propriet\'{a} caratteristica di tre sistemi $\infty^1$ di duperficie}, Atti Aca. Sc. Torino {\bf LIX} (1924), 666-671.
\bibitem{Se2} B.\ Segre, {\em Famiglie di ipersuperfie isoparemetriche negli spazi euclidei ad un qualunque numero di dimensioni}, Rend. Aca. Naz. Lincei {\bf XXVII}(1938), 203-207.
\bibitem{Sh} I.\ R.\ Shafarevich, {Basic Algebraic Geometry}, Springer-Verlag, Berlin, Heidelberg, New York, 1977.
\bibitem{So} C.\ Somigliana, {\em Sulle relazionefra il principio di huygens e l'ottica geometrica}, Atti Aca. Sc. Torino {\bf LIV} (1918-1919), 357-364. 
\bibitem{St} {S. Stolz}, {\em Multiplicities of Dupin hypersurfaces}\, {Inven. Math.}\,{\bf 138}(1999), 253-279.
\bibitem{Ta} R.\ Takagi, {\em A class of hypersurfaces with constant principal curvatures in a sphere}, J.\ Diff.\ Geom. {\bf 11}(1976), 225-233.
\bibitem{TT} R.\ Takagi and T.\ Takahashi, {\em On the principal curvatures of homogeneous hypersurfaces in a sphere}, Differential Geometry in honor of K. Yano, 469-481, Kinokuniya, Tokyo, 1972.
\bibitem{Tang} Z. Tang, {\em Isoparametric hypersurfaces with four distinct principal curvatures}, Chinese Science Bulletin {\bf 36}(1991), 1237-1240.
%\bibitem{TXY} Z.\ Tang, Y.\ Xie and W.\ Yan, {\em Gromov-Lawson-Scheon-Yau theory and isoparametric foliations}, {Comm.\ Anal.\ Geom.} {\bf 20}(2012), 989-1018.
%\bibitem{TY} Z.\ Tang and W.\ Yan, {\em Isoparametric foliation and Yau conjecture on the first eigenvalue}, {J.\ Diff.\ Geom.} {\bf 94}(2013), 521-540.
\bibitem{Th}  G.\  Thorbergsson,  {\em A  survey  on  isoparametric hypersurfaces and  their  generalizations}, Handbook  of Differential Geometry, vol. I, 963-995, North-Holland, Amsterdam, 2000.
\bibitem{Th1} G. Thorbergsson, {\em Singular Riemannian foliations and isoparametric submanifolds}, Milan J. Math. {\bf 78}(2010), 355-370.
\bibitem{Wa1} Q.-M.\ Wang, {\em Isoparametric functions on Riemannian manifolds, I}, {Math.\ Ann.} {\bf 277}(1987), 639-646.
%\bibitem{Wa2} Q,-M.\ Wang, {\em On the topology of Clifford isoparametric hypersurfaces}, {J.\ Diffe.\ Geom.} {\bf 27}(1988), 55-66.
%\bibitem{Yo} I.\ Yokota, {Exceptional Lie groups}, arXiv: 0902.0431v1.



\end{thebibliography}
\end{document}